\documentclass[11pt]{amsart}

\usepackage{epsfig, subfig, amsmath, amssymb, latexsym, amsfonts, xcolor, easybmat, bbm}

\newcommand\ol{\ensuremath{\overline}}

\newcommand\eop{{{\hfil \ensuremath \Box}}}

\newcommand\eps{\ensuremath {\varepsilon}}

\newenvironment{cor}{\subsection{}{\textbf {Corollary.}}\em}{}
\newenvironment{defn}{\subsection{}{\textbf {Definition.}}\em}{\smallskip}
\newenvironment{eg}{\subsection{}{\textbf {Example.}}}{\smallskip}

\newenvironment{lem}{\subsection{}{\textbf {Lemma.}}\em}{\smallskip}
\newenvironment{notation}{\subsection{}{\textbf{Notation.}}}{\smallskip}

\newenvironment{prop}{\subsection{}{\textbf {Proposition.}}\em}{\smallskip}

\newenvironment{rem}{\subsection{}{\textbf {Remark.}}}{\smallskip}
\newenvironment{thm}{\subsection{}{\textbf {Theorem.}}\em}{\smallskip}

\newenvironment{pf}{\noindent{\textbf {Proof.}}} {\begin{flushright}\eop \end{flushright}\smallskip}

\newcommand\fC{\ensuremath{\mathfrak C}}

\newcommand\fE{\ensuremath{\mathfrak E}}

\newcommand\fS{\ensuremath{\mathfrak S}}

\newcommand\cB{\ensuremath{\mathcal B}}
\newcommand\cC{\ensuremath{\mathcal C}}

\newcommand\cH{\ensuremath{\mathcal H}}

\newcommand\cK{\ensuremath{\mathcal K}}
\newcommand\cL{\ensuremath{\mathcal L}}

\newcommand\cP{\ensuremath{\mathcal P}}

\newcommand\bbC{\ensuremath{\mathbb C}}
\newcommand\bbD{\ensuremath{\mathbb D}}

\newcommand\bbM{\ensuremath{\mathbb M}}
\newcommand\bbN{\ensuremath{\mathbb N}}

\newcommand\bbQ{\ensuremath{\mathbb Q}}
\newcommand\bbR{\ensuremath{\mathbb R}}

\newcommand\bbT{\ensuremath{\mathbb T}}

\newcommand\bbZ{\ensuremath{\mathbb Z}}

\newcommand\ttt{\ensuremath{\textsc}}

\newcommand\bofh{\ensuremath{\cB ( \cH)}}
\newcommand\kofh{\ensuremath{\cK ( \cH)}}

\newcommand\hilb{\ensuremath{\mathcal H}}

\newcommand\norm{\ensuremath {\Vert}}

\newcommand\ran{\ensuremath {\mathrm {ran}}}

\newcommand\nul{\ensuremath {\mathrm {nul}}}

\newcommand\rank{\ensuremath{\mathrm{rank}\, }}

\newcounter{asst}
\newcounter{lab}
\newcounter{asstAA}
\newcounter{asstBB}
\newcounter{asstCC}
\newcounter{asstDD}
\newcounter{asstEE}
\newcounter{asstFF}
\newcounter{asstGG}
\newcounter{asstHH}
\newcounter{asstII}
\newcounter{asstJJ}
\newcounter{asstKK}
\newcounter{asstLL}
\newcounter{asstMM}
\newcounter{suppAA}
\definecolor{myred}{rgb}{0.6,0,0}
\definecolor{myblue}{rgb}{0,0.2,0.4}
\definecolor{mygreen}{rgb}{0.2,0.6, 0.5}

\newcommand\jnu{\ensuremath{\textsc{j}_n^{[\textsc{u}]}}}
\newcommand\jna{\ensuremath{\textsc{j}_n^{[\textsc{au}]}}}
\newcommand\jns{\ensuremath{\textsc{j}_n^{[\textsc{s}]}}}
\newcommand\jnx{\ensuremath{\textsc{j}_n^{[\textsc{x}]}}}

\newcommand\jnmx{\ensuremath{\textsc{j}_{n^m}^{[\textsc{x}]}}}
\newcommand\jnma{\ensuremath{\textsc{j}_{n^m}^{[\textsc{au}]}}}

\usepackage{mathdots}
\usepackage{datetime}

\allowdisplaybreaks

\begin{document}


\title{Stability relations for Hilbert space operators\\ and a problem of Kaplansky}



\thanks{${}^1$ Research supported in part by NSERC (Canada)}
\thanks{${}^2$ Research supported by ERC Marie Curie-Sk\l{}odowska Personal Fellowship (No.: 101064701)}
\thanks{${}^3$ Research supported in part by National Natural Science Foundation of China (No.: 12071174), Science and Technology
Development Project of Jilin Province (No.: 20190103028JH)}

\thanks{{\ifcase\month\or Jan.\or Feb.\or March\or April\or May\or
June\or
July\or Aug.\or Sept.\or Oct.\or Nov.\or Dec.\fi\space \number\day,
\number\year} \currenttime}

\author
	[L.W. Marcoux]{{Laurent W. Marcoux${}^1$}}
\address
	{Department of Pure Mathematics\\
	University of Waterloo\\
	Waterloo, Ontario \\
	Canada  \ \ \ N2L 3G1}
\email{Laurent.Marcoux@uwaterloo.ca}


\author
	[H. Radjavi]{{Heydar Radjavi}}
\address
	{Department of Pure Mathematics\\
	University of Waterloo\\
	Waterloo, Ontario \\
	Canada  \ \ \ N2L 3G1}
\email{hradjavi@uwaterloo.ca}


\author
	[S. Troscheit]{{Sascha Troscheit${}^2$}}
\address
	{Research Unit Mathematical Sciences \\
	University of Oulu\\
	PO Box 8000, 90014 Oulu\\
	Finland}
\email{sascha@troscheit.eu}


\author
	[Y.H.~Zhang]{{Yuanhang~Zhang${}^3$}}
\address
	{School of Mathematics\\
	Jilin University\\
	Changchun 130012\\
	P.R. China}
\email{zhangyuanhang@jlu.edu.cn}


\begin{abstract}
In his monograph on Infinite Abelian Groups, I.~Kaplansky raised three ``test problems" concerning their structure and multiplicity.  As noted by Azoff, these problems make sense for any category admitting a direct sum operation.   Here, we are interested in the operator theoretic version of Kaplansky's second problem which asks:  if $A$ and $B$ are operators on an infinite-dimensional, separable Hilbert space and $A \oplus A$ is equivalent to $B \oplus B$ in some (precise) sense, is $A$ equivalent to $B$?   We examine this problem under a strengthening of the hypothesis, where a ``primitive" square root $J_2(A)$ of $A\oplus A$ is assumed to be equivalent to the corresponding square root $J_2(B)$ of $B \oplus B$.  When ``equivalence" refers to similarity of operators and $A$ is a compact operator, we deduce from this stronger hypothesis that $A$ and $B$ are similar.   We exhibit a counterexample (due to J. Bell) of this phenomenon in the setting of unital rings.    Also, we  exhibit an uncountable family $\{ U_\alpha\}_{\alpha \in \Omega}$ of unitary operators, no two of which are unitarily equivalent, such that each $U_\alpha$ is unitarily equivalent  to $J_n(U_\alpha)$, a ``primitive" $n^{th}$ root of $U_\alpha \oplus U_\alpha \oplus \cdots \oplus U_\alpha$.
\end{abstract}


\keywords{Kaplansky's Problem, similarity, direct sums, multiplicity, primitive operator roots}
\subjclass[2010]{Primary: 47A45, 47A65. Secondary: 28A20}

\maketitle
\markboth{\textsc{  }}{\textsc{}}


\section{Introduction}


\subsection{} \label{sec1.01}  \ \ \
Let $\hilb$ be a complex, separable, infinite-dimensional Hilbert space, and let $\bofh$ denote the algebra of bounded linear operators acting on $\hilb$.   By $\kofh$ we denote the closed, two-sided ideal of compact operators in $\bofh$.   The \textbf{spectrum} of an element $T \in \bofh$ is denoted by $\sigma(T) := \{ \lambda \in \bbC:  T - \lambda I \mbox{ is not invertible}\}$, and the \textbf{spectral radius} of $T$ is $\mathrm{spr}(T) := \max \{ | \lambda | : \lambda \in \sigma(T)\}$.   The Beurling-Gelfand Spectral Radius Formula is the statement that $\mathrm{spr}(T) = \lim_n \norm T^n \norm^{\frac{1}{n}}$.  An operator $T \in \bofh$ is said to be \textbf{Fredholm} if its range is closed, and if $\ker\, T$ and $\ker\, T^*$ are both finite-dimensional, in which case the \textbf{Fredholm index} of $T$ is defined as $\mathrm{ind}\, T := \nul\, T - \nul\, T^*$.   (Equivalently, $T$ is Fredholm if its image $\pi(T)$ in the Calkin algebra $\bofh/\kofh$ under the canonical map $\pi: \bofh \to \bofh/\kofh$ defined by $\pi(T) := T + \kofh$ is invertible.)  We say that $T \in \bofh$ is \textbf{biquasitriangular} provided that whenever $\lambda \in \bbC$ is such that the range of $T-\lambda I$ is closed and at least one of $\nul\, T$ and $\nul\, T^*$ is finite, then both are finite and $\mathrm{ind}\, T = 0$.  (Although this is not the original definition of biquasitriangularity, it is equivalent to it by a deep result of Apostol, Foia\c{s} and Voiculescu~\cite{ApostolFoiasVoiculescu1973}.).

We adopt the standard notation:  $\bbD := \{ z \in \bbC: |z| < 1\}$, $\bbT := \{ z \in \bbC: |z| = 1\}$, and for $A \in \bofh$ and $1 \le n < \infty$, we write $A^{(n)} := A \oplus A \oplus \cdots \oplus A$ ($n$ times).

 In a recent article~\cite{MarcouxRadjaviZhang2022.02}, the first, second and fourth authors characterised the biquasitriangular operators which lie in the closure $\textsc{clos}\, (\fC_\fE)$ of the set $\fC_\fE := \{ E F - F E: E, F \in \bofh, E = E^2, F = F^2 \}$ of commutators of idempotent operators in $\bofh$.   While investigating the question of which \emph{non}-biquasitriangular operators lie in $\ttt{clos}\, (\fC_\fE)$, they demonstrated that if $S \in \bofh$ is the unilateral forward shift (i.e. if there exists an orthonormal basis $\{ e_n\}_{n \in \bbN}$ for $\hilb$ with respect to which $S e_n = e_{n+1}$ for all $n \ge 1$) and $\mu \in \bbC$, then $\mu S \in \textsc{clos}(\fC_\fE)$ if and only if $|\mu| \le \frac{1}{2}$.   A key step in this proof of this result is the fact that the unilateral shift operator $S$ is unitarily equivalent to the operator $J_2(S) := \begin{bmatrix} 0 & I \\ S & 0 \end{bmatrix} \in \cB(\hilb \oplus \hilb)$.

In the present article, we investigate which other Hilbert space operators satisfy this and closely related conditions.    A principal reason for our being interested in operators of the form $J_2(T)$ (or more generally of the form $J_n(T)$ as defined below) -- other than their relation to the set of commutators of idempotent operators mentioned above -- is their connection to a problem of Kaplansky~\cite{Kaplansky1954} originally regarding abelian groups, and later adapted to operator theory by Azoff~\cite{Azoff1995}.    Kaplansky considered two abelian groups $G$ and $H$ such that $G \oplus G$ is (group) isomorphic to $H \oplus H$, and asked whether  it necessarily follows that $G$ is isomorphic to $H$?  He recognised that analogues of this problem ``can be formulated for very general mathematical systems".   Azoff's question asks:   if $A$ and $B$ are two Hilbert space operators and $A \oplus A$ is equivalent to $B \oplus B$ (in some sense to be made precise), is $A$ equivalent to $B$?  In the finite-dimensional setting, the most important notions of equivalence of operators are unitary equivalence and similarity.    As we shall now see, the answer to Azoff's question is ``yes" in both cases.


\begin{defn} \label{defn1.02}
Let $A, B \in \bofh$.   We say that $A$ and $B$ are \textbf{unitarily equivalent} if there exists a unitary operator $U \in \bofh$ such that $B = U^* A U$, in which case we write $A \simeq B$.   We say that $A$ and $B$ are \textbf{approximately unitarily equivalent} if there exists a sequence $(U_n)_n$ of unitary operators such that $B = \lim_n U_n^* A U_n$, in which case we write $A \simeq_a B$.  Finally, we say that $A$ and $B$ are \textbf{similar}, and we write $A \sim B$, if there exists an invertible operator $S \in \bofh$ such that $B = S^{-1} A S$.
\end{defn}

It is (very) well-known that these define equivalence relations on $\bofh$.   Note that  $A$ and $B$ are approximately unitarily equivalent if and only if each belongs to the closure of the unitary orbit of the other.  (We can also define a relation on $\bofh$ by asking that $A$ belong to the closure of the similarity orbit of $B$, but this relation is not symmetric and we do not discuss it here.)

\smallskip


\subsection{} \label{sec1.02.01}
When $\hilb$ is finite-dimensional, the unitary group of $\bofh$ is compact, and therefore the relations of unitary equivalence and of approximate unitary equivalence coincide.   Furthermore,  in this setting, both unitary equivalence and similarity admit a notion of  ``multiplicity", and so we can ``recover" the operator $A$ from $A \oplus A$.  For example, if $A \oplus A$ is unitarily equivalent to $B \oplus B$, then it is routine to verify the trace condition
\[
\mathrm{tr}(w((A\oplus A), (A \oplus A)^*)) =  \mathrm{tr}(w((B \oplus B), (B \oplus B)^*)) \]
holds for all words $w(x,y)$ in two non-commuting variables $x$ and $y$.  But then
\[
\mathrm{tr} (w(A, A^*)) = \frac{1}{2} \mathrm{tr}(w((A\oplus A), (A \oplus A)^*)) = \frac{1}{2} \mathrm{tr}(w((B \oplus B), (B \oplus B)^*)) = \mathrm{tr}(w(B, B^*))\]
for all words $w(x,y)$ as above, and thus a theorem of Specht~\cite{Specht1940} implies that $A$ is unitarily equivalent to $B$.
If $A \oplus A$ is similar to $B \oplus B$, then these two operators share a common Jordan form, and thus so do $A$ and $B$ (as the Jordan form of $A$ admits the same Jordan blocks as that of $A \oplus A$, only with half the multiplicities,  and an analogous statement holds regarding $B \oplus B$ and $B$).    It follows that $A$ is similar to $B$.

When $\hilb$ is infinite-dimensional and separable, the answer to Azoff's question is still ``yes"  for unitary equivalence (Kadison and Singer~\cite{KadisonSinger1957}) and for approximate unitary equivalence (c.f. Proposition~\ref{prop4.04} below).    It is not known, however, whether the same conclusion holds when the equivalence relation under consideration is \emph{similarity}.   That is,  it is not known whether the existence of an invertible operator $S \in \cB(\hilb \oplus \hilb)$ for which $S^{-1} (A \oplus A) S = (B \oplus B)$ implies the existence of an invertible operator $R \in \bofh$ such that  $R^{-1} A R = B$, \emph{even in the case where $A$ (and therefore $B$) is compact.}

Given $A \in \bofh$, let us define $J_2(A) = \begin{bmatrix} 0 & I \\ A & 0 \end{bmatrix} \in \cB(\hilb \oplus \hilb)$.  Observe that $J_2(A)$ is a square root of $A \oplus A$, as $J_2(A)^2 = A \oplus A$.   It follows that if $B \in \bofh$ and $J_2(A)$ is similar to $J_2(B)$, then clearly $A \oplus A$ is similar to $B \oplus B$.   Whereas Azoff's question remains open for compact operators $A$ and $B$, as  we shall see below, when $A$ is compact, the special structure of $J_2(A)$ as a ``primitive" square root of $A \oplus A$ will allow us to prove that if $A$ is compact and $J_2(A)$ is similar to $J_2(B)$, then $A$ is similar to $B$ (Theorem~\ref{thm4.10} below).


Given $A \in \bofh$, we may generalise the notion of the square root $J_2(A)$ of $A \oplus A$ as follows.

\begin{defn} \label{defn1.03}
Let $2\leq n \in \bbN$.  Given $T \in \bofh$, we define
\[
J_n(T) :=
\begin{bmatrix}
	0 & I & 0 & 0 & \cdots & 0 & 0 \\
	0 & 0 & I & 0 & \cdots & 0 & 0 \\
	0 & 0 & 0 & \ddots & \cdots & \vdots & \vdots \\
	\vdots & & & \ddots & \cdots & & \vdots \\
	\vdots & & & & \ddots & & \vdots \\
	0 & 0 & 0 & 0 & \cdots & 0 & I \\
	T & 0 & 0 & 0 & \cdots & 0 & 0
\end{bmatrix} \in \cB(\hilb^{(n)}), \]
and we refer to this as the \textbf{primitive $n^{th}$ root} of $T^{(n)}$.
We then say that $T \in \bofh$ is
\begin{itemize}
	\item{}
	$\jnu$-\textbf{stable} if $T \simeq J_n(T)$;
	\item{}
	$\jna$-\textbf{stable} if $T \simeq_a J_n(T)$;
	\item{}
	$\jns$-\textbf{stable} if $T \sim J_n(T)$;
\end{itemize}	

\end{defn}

The choice of the adjective ``primitive" is a nod to the fact that if $I$ denotes the identity operator in $\bofh$, then $\sigma(I) = \{ 1 \}$, whereas  $J_n(I) \simeq \oplus_{k=1}^n \omega^k I$, and thus $\sigma(J_n(I)) = \{ \omega^k: 1 \le k \le n\}$, where $\omega = e^{\frac{2 \pi i}{n}}$ is a primitive $n^{th}$ root of $1$.
We also direct the reader's attention to Lemma~\ref{lem2.15} for further justification of this terminology.
\smallskip

In Section Two below, we investigate the basic spectral properties of $\ttt{j}_n$-stable operators, and we show that if $T \in \bofh$ is $\ttt{j}_n$-stable (in any of the above three ways), then $\sigma(T)$ is either $\bbT$ or $\ol{\bbD}$.   In Section Three, we exhibit a host of examples of non-unitarily equivalent, cyclic unitary operators, each of which is $\jnu$-stable.    To do so, we must embark on a detailed analysis of the corresponding spectral measures of these unitary operators.    In Section Four, we study our variant of the Kaplansky/Azoff problem, obtaining the result mentioned above that if $A, B \in \bofh$ and $A$ is compact,  then $J_2(A)$ similar to $J_2(B)$ implies that $A$ is similar to $B$.  Finally, in Section Five, we show that this last result cannot be generalised in a purely algebraic way.     More specifically, we reproduce an example given to us by J.~Bell of a unital ring $R$ and two elements $x, y \in R$ such that $J_2(x) := \begin{bmatrix} 0 & 1 \\ x & 0 \end{bmatrix}$ is conjugate in $\bbM_2(R)$ to $J_2(y) = \begin{bmatrix} 0 & 1 \\ y & 0 \end{bmatrix}$, yet $x$ is not conjugate to $y$.

\bigskip

Some results we shall mention below apply equally well to all notions of $\ttt{j}_n$-stability.   To avoid repetition, when proving such results we shall write that $T$ is $\jnx$-stable, meaning that $[\ttt{x}] \in \{ [\ttt{u}], [\ttt{au}], [\ttt{s}]\}$.


\begin{rem} \label{rem1.04}
The study of $\ttt{j}_n$-stability of operators extends the work -- initiated by Conway, Pr\v{a}jitur\v{a} and Rodriguez-Mart\'{i}nez~\cite{ConwayPrajituraRodriguezMartinez2014} -- of those operators $T \in \bofh$ such that $T^2$ is similar (or unitarily equivalent) to $T^{(2)} $.     Since $J_n(T)^n = T^{(n)}$, it follows that if $T$ is $\jnu$-stable, then $T^n \simeq T^{(n)}$, and analogous statements hold regarding $\jna$-stability and $\jns$-stability.

Moreover, it is clear that if $T$ is  $\jnu$-stable, then $T$ is both $\jna$-stable and $\jns$-stable.    In general, neither $\jna$-stability nor $\jns$-stability implies $\jnu$-stability, as we shall see.
\end{rem}


\begin{eg} \label{eg1.05}
Having defined various notions of ``stability" of operators, it is worthwhile verifying that the various sets of $\textsc{j}_n$-stable operators are non-empty.

\smallskip

If $U \in \bofh$ is any unitary operator satisfying $\sigma(U) = \bbT$, then $U$ is $\jna$-stable for all $n \ge 2$.   This is because $J_n(U)$ is again a unitary operator with $\sigma(J_n(U)) = \bbT$, and any two unitary operators with spectrum equal to $\bbT$ are approximately unitarily equivalent by the Weyl-von Neumann-Berg Theorem~\cite[Theorem~II.4.4]{Davidson1996}.

To see that $\sigma(J_n(U)) = \bbT$, note that if $A \simeq_a B$, then a routine calculation shows that $J_n(A) \simeq_a J_n(B)$ for any $n \ge 2$.   As such, using the Weyl-von Neumann-Berg Theorem, we may assume without loss of generality that $U = \ttt{diag}(u_m)_m$, where $\{ u_m\}_m$ is dense in $\bbT$.   With $U$ diagonal as above,
\[
J_n(U) \simeq \oplus_m
\begin{bmatrix}
	0 & 1 & 0 & 0 & \cdots & 0 & 0 \\
	0 & 0 & 1 & 0 & \cdots & 0 & 0 \\
	0 & 0 & 0 & \ddots & \cdots & \vdots & \vdots \\
	\vdots & & & \ddots & \cdots & & \vdots \\
	\vdots & & & & \ddots & & \vdots \\
	0 & 0 & 0 & 0 & \cdots & 0 & 1 \\
	u_m & 0 & 0 & 0 & \cdots & 0 & 0
\end{bmatrix}. \]
Clearly $J_n(U)$ is unitary, and  $\sigma(J_n(U))$ contains all $n^{th}$ roots of $u_m$ for each $m \ge 1$.  Moreover, the set of all such $n^{th}$ roots of  $\{ u_m\}_m$ is dense in $\bbT$, since $\{ u_m\}_m$ itself is.  Thus $J_n(U)$ is a unitary operator whose spectrum is the circle, whence $J_n(U) \simeq_a U$.
\end{eg}


\begin{eg} \label{eg1.06}
It is routine to verify that the unilateral forward shift operator $S$ is $\jnu$-stable (and thus both $\jna$-stable and $\jns$-stable) for all $n \ge 2$.  The computation is omitted.
\end{eg}


\begin{eg} \label{eg1.07}
There exists an operator $T\in \bofh$ such that $T$ is  $\jns$-stable for all $n \in 2 \bbN$, but $T$ is not $\jnu$-stable for any $n \in 2 \bbN$.

Let $\{e_k\}_{k=1}^\infty$ be an orthonormal basis of $\mathcal{H}$, and $S$ be the unilateral forward shift operator with respect to $\{ e_k\}_k$ (i.e. $S e_k = e_{k+1}$ for all $k \ge 1$).  Define an invertible diagonal operator $D\in \bofh$ by
\[
D e_k := \begin{cases} \frac{1}{2} e_k & \mbox{ if } k \mbox{ is odd} \\
	2 e_k & \mbox{ if } k \mbox{ is even.}
	\end{cases}\]
Set $T := D^{-1}SD$.  It is straightforward to check that $T^2 \simeq S \oplus S$ is an isometry.
Since $S$ is $\jnu$-stable and $J_{2m}(T)=(D^{-1})^{(2m)}J_{2m}(S) D^{(2m)}$ for all $m \ge 1$, it follows that $T$ is $\jns$-stable for all $n \in 2 \bbN$.

If $T \simeq J_{2n}(T)$ for some $n\in \mathbb{N}$, then $T^{2n}\simeq T^{(2n)}$.
This is impossible as $T^{2n}$ is an isometry,  whereas $T^{(2n)}$ is not.
Hence, $T$ is not $\jnu$-stable for any $n \in 2 \bbN$.
\end{eg}



\vskip 1 cm

\section{$\ttt{j}_n$-stability} \label{section2}


\subsection{} \label{sec2.01}
We begin with an examination of the spectral properties shared by $\ttt{j}_n$-stable operators.  It is worth noting that these properties may often be derived by observing that if $n \in \bbN$ and $T \in \bofh$ is $\jnu$-stable (resp. $\jna$-stable, $\jns$-stable), then $T^n \simeq T^{(n)}$ (resp. $T^n \simeq_a T^{(n)}$, $T^n \sim T^{(n)}$).   This last property (for unitary equivalence and similarity) was originally investigated for the case where $n=2$ in the paper~\cite{ConwayPrajituraRodriguezMartinez2014} by Conway, Pr\v{a}jitur\v{a} and Rodr\'{i}guez-Mart\'{i}nez, who referred to these as $\ttt{Property }U$ and $\ttt{Property }S$ respectively.  Amongst (many) other things, they proved that the only compact operator $K \in \kofh$ which satisfies  $K^2 \sim K^{(2)}$ is the zero operator.


\smallskip

\begin{prop} \label{prop2.02}
Let $T \in \bofh$ and $2 \le n \in \bbN$.   If $T^{n} \equiv T^{(n)}$ \emph{(}where $\equiv$ denotes one of $\simeq, \simeq_a$ or $\sim$\emph{)}, then $T^{n^m} \equiv T^{(n^m)}$ for all $m \ge 1$.    Consequently,
\[
\sigma(T) = \sigma(T)^{n^m} \mbox{\ \ \ \ \ for all } m \ge 1. \]
In particular, if $T$ is $\jnx$-stable, then $\sigma(T) = \sigma(T)^{n^m}$ for all $m \ge 1$.
\end{prop}

\begin{pf}
We argue the case where $\equiv$ denotes approximate unitary equivalence.   The proof  is an easy adaptation of Proposition~2.7 of~\cite{ConwayPrajituraRodriguezMartinez2014}, as are the corresponding proofs for the other two cases.

In this case we see that
	\[
	T^{n^2} = (T^n)^n \simeq_a (T^{(n)})^n = (T^n)^{(n)} \simeq_a (T^{(n)})^{(n)} = T^{(n^2)}, \]
	and by an easy induction argument,
	\[
	T^{n^m} \simeq_a T^{(n^m)} \mbox{ for all } m \ge 1. \]
	From this it follows that for all $m \ge 1$,
	\[
	\sigma(T)^{n^m} = \sigma(T^{n^m}) = \sigma(T^{(n^m)}) = \sigma(T). \]
That the last statement of the proposition holds is clear from Paragraph~\ref{sec2.01}.
\end{pf}


\begin{prop} \label{prop2.03}
Let $T \in \bofh$ and $2 \le n \in \bbN$, and suppose that  $T^{n} \equiv T^{(n)}$ \emph{(}where $\equiv$ denotes one of $\simeq, \simeq_a$ or $\sim$\emph{)}.   Then
\begin{enumerate}
	\item[(a)]
	$\mathrm{spr}(T) \in \{ 0, 1\}$; and
	\item[(b)]
	if $T$ is invertible,  then $\sigma(T) \subseteq \bbT$.
\end{enumerate}
In particular, these properties hold if $T$ is $\jnx$-stable, where $[\ttt{x}] \in \{ [\ttt{u}], [\ttt{au}], [\ttt{s}]\}$.
\end{prop}

\begin{pf}
\begin{enumerate}
	\item[(a)]
	Using Proposition~\ref{prop2.02}, we find that
	\[
	0 \le \mathrm{spr}(T) =  \mathrm{spr}(T^{n}) = \mathrm{spr}(T)^n,\]
	from which the result clearly follows.
	\item[(b)]
	By part (a) above, $\sigma(T) \subseteq \ol{\bbD}$.   Suppose that $\alpha \in \sigma(T)$ and that $0 < |\alpha| < 1$.   By Proposition~\ref{prop2.02},
	\[
	\alpha^{n^m} \in \sigma(T)^{n^m} = \sigma(T), \]
	and since $\lim_m \alpha^{n^m} = 0$ and $\sigma(T)$ is closed, we have that $0 \in \sigma(T)$.    Thus the invertibility of $T$, combined with the hypothesis that $T^n \equiv T^{(n)}$, implies that $\sigma(T) \cap \bbD = \varnothing$.   Hence $\sigma(T) \subseteq \bbT$.
\end{enumerate}
\end{pf}


\begin{lem} \label{lem2.03.01}
Let $A, B \in \bofh$ and $n \in \bbN$.
\begin{enumerate}
	\item[(a)]
	If $A \simeq B$, then $J_n(A) \simeq J_n(B)$.
	\item[(b)]
	If $A \simeq_a B$, then $J_n(A) \simeq_a J_n(B)$.
	\item[(c)]
	If $A \sim B$, then $J_n(A) \sim J_n(B)$.
\end{enumerate}
\end{lem}

\begin{pf}
Since all three proofs are essentially the same, we shall only prove (c).

\smallskip

If $A \sim B$, say $B = R^{-1} A R$ for some $R \in \bofh$ invertible, we set $S = R^{(n)}$ and note that $S$ is invertible in $\cB(\hilb^{(n)})$.   A simple computation then shows that
\[
J_n(B) = S^{-1}  J_n(A) S, \]
whence $J_n(B) \sim J_n(A)$.
\end{pf}


\smallskip

We next prove the analogue of Proposition~\ref{prop2.02} for $\ttt{j}_n$-stable operators.
\smallskip

\begin{prop} \label{prop2.04}
Let $R, T \in \bofh$, and let $[\ttt{x}] \in \{ [\ttt{u}], [\ttt{au}], [\ttt{s}]\}$.
\begin{enumerate}
	\item[(a)]
	Let $m, n \ge 1$.   Then $J_m(J_n(R)) \simeq J_{mn} (R)$.
	\item[(b)]
	If $T$ is $\jnx$-stable for some $n \ge 2$, then $T$ is $\jnmx$-stable for all $m \ge 1$.
\end{enumerate}	
\end{prop}

\begin{pf}
\begin{enumerate}
	\item[(a)]
	This is a routine computation which is omitted.
	\item[(b)]
	We argue by induction on $m$.   Let $[\ttt{x}] = [\ttt{s}]$.  By hypothesis, the result holds for $m =1$.   Suppose that it holds for $m = k \ge 1$.   Then by part (a) and Lemma~\ref{lem2.03.01} (c),
	\[
	T \sim J_{n^{k}}(T) \sim J_{n^{k}}( J_n(T)) \simeq J_{n^{k + 1}}(T). \]
	That is, the result holds for $m = k +1$.  By induction, the result holds for all $m \ge 1$.
	
	The cases where $[\ttt{x}] = [\ttt{u}]$ or $[\ttt{x}] = [\textsc{au}]$ are handled in a similar manner.
\end{enumerate}	
\end{pf}


\begin{prop} \label{prop2.05}
Let $2 \le n \in \bbN$ and suppose that $T \in \bofh$ is $\jna$-stable. Then
\begin{enumerate}
	\item[(a)]
	$\norm T^m \norm =\norm T \norm$, for all $m \in \mathbb{N}$.
	\item[(b)]
	The spectral radius of $T$ is $\mathrm{spr}(T) = 1$.	
\end{enumerate}	
\end{prop}

\begin{pf}
\begin{enumerate}
	\item[(a)]
	First observe that $\norm J_n(T)\norm \ge 1$ for any $T \in \bofh$ and any $n \ge 2$.    Thus if $T$ is $\jna$-stable, then $\norm T \norm \ge 1$.
	\smallskip
	
	We shall argue by induction on $m$.   For each $m \in \bbN$, let $P(m)$ be the statement that $\norm T^m \norm = \norm T \norm$.   Obviously $P(1)$ is true.  Now fix $2 \le m \in \bbN$, and suppose that $P(j)$ is true, $1 \le j \le m-1$.
	\smallskip
	
	Choose $k \in \bbN$ such that $(k-1) n < m \le k n$.   A routine calculation shows that  $(J_n(T))^m$ admits an operator-matrix form $[X_{i j}] \in \cB(\hilb^{(n)})$, each of whose rows and columns contains exactly one non-zero entry, and that entry is either of the form $T^{k-1}$ or $T^k$ with $T^k$ appearing $m - (k-1)n$ times amongst the non-zero entries.  Hence, whether $T$ is $\jnu$-stable or $\jna$-stable,
	\[
	\norm T^m \norm = \max( \norm T^{k-1} \norm, \norm T^k \norm). \]
	If $k = 1$, then $\norm T^{k-1} \norm = \norm I \norm = 1 \le \norm T \norm = \norm T^k \norm$, so $\norm T^m \norm = \norm T \norm$.   If $k > 1$, then we note that $k < m$, and thus by our induction step, $\norm T^m \norm = \max(\norm T^{k-1} \norm, \norm T^k \norm) = \norm T \norm$.
	
	This completes the induction and the proof of (a).
	\item[(b)]
	By part (a), $\norm T \norm \ge 1$, and so $\norm T \norm \ne 0$.  Applying Beurling's Spectral Radius Formula, we find that
	\[
	\mathrm{spr}(T) = \lim_m \norm T^m \norm^{\frac{1}{m}} = \lim_m \norm T\norm^{\frac{1}{m}} = 1. \]

\end{enumerate}
\end{pf}


\begin{rem} \label{rem2.05.01}
It is worth noting that the full hypothesis of $\jnu$-stability (or $\jna$-stability) was required to prove the above Proposition.   Indeed, suppose that $T = \begin{bmatrix} 0 & 2 I \\ \frac{1}{2} I & 0 \end{bmatrix} \in \cB(\hilb^{(2)})$.    It is straightforward to check that $T^2 = \begin{bmatrix} I & 0 \\ 0 & I \end{bmatrix}$ and thus $T^3 = T$.   Furthermore, $T \simeq T^{(3)}$, whence $T^3 = T \simeq T^{(3)}$.   Nevertheless, $\norm T \norm = 2  \ne \norm T^2 \norm = 1$.
\end{rem}

\begin{prop} \label{prop2.06}
Let $2 \le n \in \bbN$ and suppose that $T \in \bofh$.  Let $[\ttt{x}] \in \{ [\ttt{u}], [\ttt{au}], [\ttt{s}]\}$, and  set $\theta_m := e^{i \frac{2 \pi }{n^m}}$ for each $m \ge 1$.
\begin{enumerate}
	\item[(a)]
	If $T$ is  $\jnu$-stable, then $T \simeq \theta_m^k \, T$ for each $m \ge 1$ and $1 \le k \le n^m$.  	
	\item[(b)]
	If $T$ is $\jnx$-stable, then $\sigma(T)$ has circular symmetry.
	\item[(c)]
	If $T$ is $\jnx$-stable and invertible, then $\sigma(T) = \bbT$.
\end{enumerate}
\end{prop}

\begin{pf}

\begin{enumerate}
	\item[(a)]
	Since $T$ is $\jnu$-stable, for each $m \ge 1$, $T$ is unitarily equivalent to $J_{n^m}(T)$ by Proposition~\ref{prop2.04}.  Let $m \in \bbN$ and consider $\theta_m$ as in the statement of the Proposition.
	Let
	\[
	V_m := \ttt{diag} (I, \theta_m\, I, \theta_m^2 \, I, \ldots, \theta_m^{n^m-1} \, I) \in \cB(\hilb^{(n^m)}). \]
	Then $V$ is clearly unitary and
	\[
	T \simeq J_{n^m}(T) \simeq V_m^* J_{n^m}(T) V_m = \theta_m J_{n^m}(T) \simeq \theta_m \, T. \]
	Of course, $T \simeq \theta_m \, T$ implies that $T \simeq \theta_m \, T \simeq \theta_m \, (\theta_m\, T) = \theta_m^2 \, T$, and
	by a simple induction argument, $T \simeq \theta_m^k \, T$ for all $k \ge 1$.  But $\theta_m^{k} = \theta_m^{k+n^m}$
	for all $k \ge 1$, so that this reduces to the statement that $T \simeq \theta_m^k \, T$, $1 \le k \le n^m$.
	\item[(b)]
	\begin{itemize}
		\item{}
		It is routine to verify that by replacing each instance of $\simeq$ by $\simeq_a$ in the proof of (a) above, we arrive at the conclusion that if $T$ is $\jna$-stable, then $T \simeq_a \theta_m^k T$ for all $k, m \ge 1$.
	
	But the set $\{ \theta_m^k : k, m \ge 1 \}$ is dense in $\bbT$, and thus
	\[
	T \simeq_a \lambda \, T \mbox{ for all } \lambda \in \bbT. \]
	From this it easily follows that if $T$ is $\jna$-stable, then
	\[
	\sigma(T) = \sigma(\lambda \, T) = \lambda \, \sigma(T) \mbox{ for all } \lambda \in \bbT. \]
	In particular, if $T$ is $\jnu$-stable, then $\sigma(T)$ has circular symmetry.
		\item{}
		If $T$ is $\jns$-stable, then $T \sim J_{n^m}(T)$ for all $m \ge 1$, again by Proposition~\ref{prop2.04}.  Let $V_m$ denote the unitary operator from (a).   	Then
		\[
		T \sim J_{n^m}(T) \simeq V_m^* J_{n^m}(t) V_m = \theta_m J_{n^m}(T) \sim \theta_m T, \]
		whence $T \sim \theta_m T \sim \theta_m (\theta_m T) = \theta_m^2 T$, and as before, by a simple induction argument, $T \sim \theta_m^k T$ for all $k \ge 1$.
		Thus $\sigma(T)$ is invariant under rotation by $\theta_m^k$ for all $m, k \ge 1$, and since this set is dense in $\bbT$, it follows that $\sigma(T)$ has circular symmetry.
	\end{itemize}	
	
	\item[(c)]
	This follows from (b), together with Proposition~\ref{prop2.03}.
\end{enumerate}
\end{pf}


\begin{prop} \label{prop2.07}
Let $T \in \bofh$ and $n \ge 2$ be an integer.   Suppose that $T$ is $\jns$-stable.   Then
\begin{enumerate}
	\item[(a)]
	$\sigma(T)=\bbT$, if $T$ is invertible; and
	\item[(b)]
	$\sigma(T)=\overline{\mathbb{D}}$, if $T$ is not invertible.
\end{enumerate}
\end{prop}

\begin{pf}
\begin{enumerate}
	\item[(a)]
	This is Proposition~\ref{prop2.06} (c).
	\item[(b)]	
	Since $T\sim J_n(T)$, there exists an invertible operator $R$, such that $T=RJ_n(T)R^{-1}$.

	\smallskip
	
	\textsc{Step One.} Since $0\in \sigma(T)$, it follows that either $T$ is not bounded below, or that the range of $T$ is not
dense. In the latter case, observe that $T^*$ is not bounded below.
Since $T\sim J_n(T)$ if and only if $T^*\sim J_n(T^*)$, and since $\sigma(T^*)=\sigma(T)^*=\{\overline{\lambda}:\lambda\in \sigma(T)\}$, it suffices to consider the case where $T$ is not bounded below, which we now do.

Fix $0<\lambda<1$. For each $m\geq 1$, we can find $x_m\in \hilb$ such that $\|x_m\|=1$ but $\|Tx_m\|\leq \lambda^{n^m}$.
Noting that $T\sim J_n(T)$ implies that $T\sim J_{n^m}(T)$ for all $m\geq 1$, we set
\[y_m:=(\lambda x_m,\lambda^2 x_m,\cdots,\lambda^{n^m-1}x_m,\lambda^{n^m}x_m)^t\in (\hilb)^{(n^m)}.\]
Then
\[\|y_m\|^2=\sum_{k=1}^{n^m}\lambda^{2k}=\frac{\lambda^2(1-(\lambda^2)^{n^m+1})}{1-\lambda^2}.\]
Note that $\lim_m \norm y_m \norm^2 = \frac{\lambda^2}{1-\lambda^2}$, which is obviously independent of $m$.

Moreover,
\[(J_{n^m}(T)-\lambda I)y_m=(0,0,\cdots,0,\lambda(Tx_m-\lambda^{n^m}x_m))^t,\]
implying that
\[\|(J_{n^m}(T)-\lambda I)y_m\|<\lambda(\lambda^{n^m}+\lambda^{n^m})=2\lambda^{n^m+1} < 2 \lambda^{n^m}.\]

\bigskip

\textsc{Step Two.} Next, observe that as $J_n(T)=R^{-1}TR$, for each $m\geq 1$, there exists an invertible operator
$R_m$ such that
\[
	J_{n^m}(T)=R_{m}^{-1}TR_m,\]
and $\norm R_m\norm \leq \|R\|^m$, $\|R_m^{-1}\|\leq \|R^{-1}\|^m$.

It follows that
\[
	\|R_m^{-1}(T-\lambda I)R_my_m\|<2\lambda^{n^m}.\]
Now $R$ is bounded below by $\frac{1}{\|R^{-1}\|}$,
and thus $R_m$ is bounded below by $\frac{1}{\|R^{-1}\|^m}$,
implying that if $z_m:=R_my_m$, then
\[\|z_m\|\geq \frac{1}{\|R^{-1}\|^m}\|y_m\|.\]
Similarly, $R^{-1}$ is bounded below by $\frac{1}{\|R\|}$,
and thus $R_m^{-1}$ is bounded below by $\frac{1}{\|R\|^m}$.

Thus (for some $m$ sufficiently large),
\begin{align*}
2\lambda^{n^m}
	&\geq  \|R_m^{-1}(T-\lambda I)R_my_m\|  \\
	&\geq \frac{1}{\|R\|^m} \|(T-\lambda I)z_m\| \\
	&\geq \frac{1}{\|R\|^m}\frac{1}{\|R^{-1}\|^m}\ \sqrt{\frac{\lambda^2}{2(1-\lambda^2)}} \ \|(T-\lambda I)[\frac{z_m}{\|z_m\|}]\|.
\end{align*}	
In other words, (for $m$ sufficiently large),
\[\|(T-\lambda I)[\frac{z_m}{\|z_m\|}] \norm \leq \|R\|^m\|R^{-1}\|^m \ \sqrt{\frac{2(1-\lambda^2)}{\lambda^2}}\ 2\lambda^{n^m}.\]
But
\[\underset{m}{\lim}\|R\|^m\|R^{-1}\|^m\cdot\lambda^{n^m}=0.\]
This proves that $T-\lambda I$ is not bounded below, so $\lambda\in \sigma(T)$.
Thus, $\sigma(T)=\overline{\mathbb{D}}$, since $\sigma(T)$ has circular symmetry and $\mathrm{spr}(T)=1$ (by Proposition~\ref{prop2.03}).
\end{enumerate}
\end{pf}


Since $\jnu$-stability obviously implies $\jns$-stability, the next result is immediate.

\smallskip

\begin{cor} \label{cor2.08}
Let $2 \le n \in \bbN$, and suppose that  $T \in \bofh$ is $\jnu$-stable.    Either
\begin{enumerate}
	\item[(a)]
	$T$ is invertible and $\sigma(T) = \bbT$; or
	\item[(b)]
	$T$ is not invertible and $\sigma(T) = \ol{\bbD}$.
\end{enumerate}
\end{cor}
	

The proof that if $T \in \bofh$ is $\jna$-stable, then $\sigma(T) = \bbT$ or $\sigma(T) = \ol{\bbD}$ is an adaptation of that of Proposition~\ref{prop2.07}.

\smallskip

\begin{prop} \label{prop2.09}
Let $T \in \bofh$ and $n \ge 2$ be an integer.   Suppose that $T$ is $\jna$-stable.   Then
\begin{enumerate}
	\item[(a)]
	$\sigma(T)=\bbT$, if $T$ is invertible; and
	\item[(b)]
	$\sigma(T)=\overline{\mathbb{D}}$, if $T$ is not invertible.
\end{enumerate}
\end{prop}

\begin{pf}
\begin{enumerate}
	\item[(a)]
	Once again, this is Proposition~\ref{prop2.06} (c).
	\item[(b)]	
	
	\textsc{Step One.}
	The argument of Step One of the proof of Proposition~\ref{prop2.07} carries over \emph{mutatis mutandis} to our setting by replacing every occurrence of $\sim$ in that proof with the current equivalence relation $\simeq_a$.   Thus, the problem reduces to showing that if $T$ is not bounded below, then $(T-\lambda I)$ is not bounded below for any $0 < \lambda < 1$.   Choosing such a $\lambda$, we construct the vectors $(y_m)_m$ as in that proof, and recall that $\lim_m \norm y_m \norm^2 = \frac{\lambda^2}{1-\lambda^2} > 0$.  Choose $m_0 \ge 1$ such that $m \ge m_0$ implies that $\norm y_m \norm^2 > \frac{1}{2} \frac{\lambda^2}{1 - \lambda^2}$.  We still have that for all $m \ge 1$,
\[
	\|(J_{n^m}(T)-\lambda I)y_m\|<\lambda(\lambda^{n^m}+\lambda^{n^m}) < 2\lambda^{n^m}.\]

\bigskip

\textsc{Step Two.}   By Proposition~\ref{prop2.04}(b), the fact that $T$ is $\jna$-stable implies that $T$ is $\jnma$-stable for all $m \ge 1$.  Thus, for each $m \ge m_0$, we may choose a unitary operator $V_m$ such that $\norm J_{n^m}(T) - V_m^* T V_m \norm < \lambda^{n^m} \cdot \frac{1}{\norm y_m\norm}$.

From this we see that
\[
	\|V_m^*(T-\lambda I)V_m y_m\| \le \norm V_m^* T V_m - J_{n^m}(T) \norm \ \norm y_m \norm +  \|(J_{n^m}(T)-\lambda I)y_m\| < 3 \lambda^{n^m}.\]
	
	With $z_m := V_m y_m$, for all $m \ge m_0$ we have that $\norm z_m \norm^2 \ge \frac{1}{2} \frac{\lambda^2}{1 - \lambda^2}$ and that
	\[
	\norm (T-\lambda I) z_m \norm = \norm V_m^*(T-\lambda I)V_m y_m\norm < 3 \lambda^{n^m}. \]
	
	This proves that $T-\lambda I$ is not bounded below, so $\lambda\in \sigma(T)$.  Since $0 < \lambda < 1$ was arbitrary, and since $\sigma(T)$ has circular symmetry and is contained in $\ol{\bbD}$, we conclude that in fact $\sigma(T)=\overline{\mathbb{D}}$.
\end{enumerate}
\end{pf}


\begin{prop} \label{prop2.10}
Suppose that $2 \le n \in \bbN$, and that $(T_k)_{k=1}^\infty$ is a sequence of operators in $\bofh$, each of which is $\jnu$-stable.  Then for all $m \in \bbN \cup \{ \infty\}$,
\[
	\oplus_{k=1}^m T_k \]
is $\jnu$-stable.
\end{prop}

\begin{pf}
We omit the proof of the case where $m \in \bbN$, as it is an easy adaptation of the following proof.   Let $T := \oplus_{k=1}^\infty T_k$.   By hypothesis, for each $k \ge 1$, there exists a unitary operator $U_{k}$ such that $J_n(T_k) = U_{k}^* T_k U_{k}$ for all $k \ge 1$. Thus, by identifying $\oplus_{k=1}^\infty \hilb^{(n)}$ with $(\oplus_{k=1}^\infty \hilb)^{(n)}$ and defining the unitary operator $U := \oplus_{k=1}^\infty U_k$, we find that
\[
	J_n(T) = J_n(\oplus_{k=1}^\infty T_k) \simeq \oplus_{k=1}^\infty J_n(T_k) = \oplus_{k=1}^\infty U_k^* T_k U_k = U^* T U. \]

Thus, $T$ is $\jnu$-stable.
\end{pf}


\begin{rem} \label{rem2.11}
It is not hard to see that if $(T_k)_k$ is a sequence of operators in $\bofh$, each of which is $\jna$-stable, then $T := \oplus_{k=1}^m T_k$ is $\jna$-stable for all $m \in \bbN \cup \{ \infty\}$.   The proof is a routine adaptation of that of Proposition~\ref{prop2.10}.

When each $T_k$ is $\jns$-stable, replacing each $U_k$ by the appropriate invertible operator $R_k$, $k \ge 1$ shows that if $m \in \bbN$, then $\oplus_{k=1}^m T_k$ is also $\jns$-stable.   There is a problem when $m = \infty$, however, as there is no reason why $R := \oplus_{k=1}^\infty R_k$ should be invertible.
\end{rem}


\begin{eg} \label{eg2.12}
The converse of Proposition~\ref{prop2.10} also fails.

Let $\hilb$ be an infinite-dimensional, separable, complex Hilbert space and suppose that $\{ e_n\}_{n \in \bbZ}$ is an orthonormal basis for $\hilb$.   Let $B \in \bofh$ be the bilateral shift operator satisfying $B e_n = e_{n+1}$, $n \in \bbZ$.     An easy computation shows that $B$ is $\jnu$-stable for all $n\geq 2$.

Let $\Omega_+ := \{ z \in \bbT:   \mathrm{Im} (z) \ge 0\}$ and $\Omega_- := \{ z \in \bbT: \mathrm{Im} (z) \le 0\}$.   Obviously $\bbT = \Omega_1 \cup \Omega_2$.   Since $B$ is unitary and $\sigma(B) = \bbT$, we see that $\Omega_+$ and $\Omega_-$ are Borel subsets of $\sigma(B)$.   By the Borel functional calculus, we may write $B \simeq B_+ \oplus B_-$, where $B_+, B_-$ are unitary, $\sigma(B_+) = \Omega_+$ and $\sigma(B_-) = \Omega_-$.

By Proposition~\ref{prop2.06} (b), neither $B_+$ nor $B_-$ are $\jnu$-stable, as neither operator's spectrum has circular symmetry.
\end{eg}


\subsection{} \label{sec2.13}
We now prove a couple of results which will be useful in analysing the examples of the next section.

Our first Lemma is well-known;  we include its proof for the convenience of the reader.  Recall that if $\hilb_1$, $\hilb_2$ are complex Hilbert spaces, $X \in \cB(\hilb_2)$, $Y \in \cB(\hilb_1)$, then the associated \textbf{Rosenblum operator} is the operator $\tau_{X, Y} \in \cB(\cB(\hilb_2, \hilb_1))$ defined by  $\tau_{X,Y}(Z) := Z X - Y Z$, $Z \in \cB(\hilb_2, \hilb_1)$.  If, furthermore, $\sigma(X) \cap \sigma(Y) = \varnothing$, then $\tau_{X, Y}$ is invertible.   We shall refer to this result as \emph{Rosenblum's Theorem}~\cite[Corollary 3.2]{Herrero1989}.


\begin{lem} \label{lem2.16}
Let $\hilb_1$ and $\hilb_2$ be two complex Hilbert spaces.   Let $A_1, B_1 \in \cB(\hilb_1)$ and $A_2, B_2 \in \cB(\hilb_2)$.   Suppose that $\Omega_1$ and $\Omega_2$ are two non-empty, disjoint, compact subsets of  $\bbC$ and that $\sigma(A_1) = \Omega_1 = \sigma(B_1)$, while $\sigma(A_2) = \Omega_2 = \sigma(B_2)$.

If $A_1 \oplus A_2$ is similar to $B_1 \oplus B_2$, then $A_j$ is similar to $B_j$, $j = 1, 2$.
\end{lem}

\begin{pf}
Suppose that $T = [T_{ij}]$ is invertible and $T (A_1 \oplus A_2) = (B_1 \oplus B_2) T$.  Then
\begin{align*}
\begin{bmatrix} T_1 A_1 & T_2 A_2 \\ T_3 A_1 & T_4 A_2 \end{bmatrix}
	&= \begin{bmatrix} T_1 & T_2 \\ T_3 & T_4 \end{bmatrix} \begin{bmatrix} A_1 & 0 \\ 0 & A_2 \end{bmatrix} \\
	&= \begin{bmatrix} B_1 & 0 \\ 0 & B_2 \end{bmatrix}  \begin{bmatrix} T_1 & T_2 \\ T_3 & T_4 \end{bmatrix} \\
	&= \begin{bmatrix} B_1 T_1 & B_1 T_2 \\ B_2 T_3 & B_2 T_4 \end{bmatrix}.
\end{align*}	
From this it follows that $T_2 A_2 = B_1 T_2$ and that $T_3 A_1 = B_2 T_3$.   But
\[
\sigma(A_2) \cap \sigma(B_1) = \varnothing = \sigma(A_1) \cap \sigma(B_2), \]
so by Rosenblum's theorem, $T_2 = T_3 = 0$.   This implies that both $T_1$ and $T_4$ are invertible.  But then $T_1 A_1 = B_1 T_1$ and $T_4 A_2 = B_2 T_4$ implies that $A_1$ is similar to $B_1$ and $A_2$ is similar to $B_2$.
\end{pf}


\begin{lem}\label{lem2.15}
Let $2\leq n\in \mathbb{N}$ and let $U$ be a unitary operator.   Let $V$ be any $n^{th}$ root of $U$.   Then $U$ is $\jnu$-stable if and only if
\[U\simeq \textup{diag}(V,\omega V,\cdots,\omega^{n-1}V),\]
where $\omega$ is a primitive $n$-th root of unity in $\bbC$.
\end{lem}

\begin{pf}
Let $W := \textup{diag}(V, V^2,\cdots,V^n)$,
and note that
\[W^*J_n(U)W=\begin{bmatrix}0& V & 0&\cdots& 0\\
0&0&V&\cdots&0\\
\vdots&\vdots&\ddots&\ddots&\vdots\\
0&0&\cdots&0&V\\
V&0&\cdots&0&0\\
\end{bmatrix}.\]

But the right-hand side is $P\otimes V$, where $P$ is a cyclic permutation matrix of order $n$.
Since $P\simeq  \textup{diag}(1,\omega,\cdots,\omega^{n-1})$, the claim follows.

\end{pf}


\begin{rem} \label{rem2.16}
If $n \in \bbN$ and $R \in \bofh$ is invertible and admits an $n$-th root $Y$, then an analogous argument to that above shows that $R$ is $\jns$-stable if and only if
\[
R \sim \mathrm{diag} (Y, \omega Y, \omega^2 Y, \ldots, \omega^{n-1} Y), \]
where $\omega$ is a primitive $n$-th root of unity in $\bbC$.
\end{rem}




\bigskip

\section{Examples} \label{sec3}


\subsection{} \label{sec3.01}

We now turn our attention to the question of determining which normal operators are $\ttt{j}_n$-stable.

Suppose that $N \in \cB(\hilb)$ is normal and $\jnu$-stable for some $2\leq n\in \bbN$.   Then $J_n(N)$ is normal, and thus
\[
I^{(n-1)} \oplus N N^* = J_n(N)J_n(N)^*=J_n(N)^*J_n(N) = N^* N \oplus I^{(n-1)}. \]
This shows that $N N^*=N^*N=I$, implying that $N$ is in fact a unitary operator.

Unitary operators, and more generally normal operators, are described up to unitary equivalence by their spectrum and their spectral measures, and we shall investigate $\jnu$-stability of unitary operators in terms of these.    A key observation lies in the fact that any $\jnu$-stable operator $T \in \bofh$ satisfies $T^n \simeq T^{(n)}$, which imposes a great deal of structure on the spectral measure of $T$ when $T$ is unitary.

We first concentrate on the case where the spectral measure is atomic; in other words, where the unitary operator is diagonalisable with respect to some orthonormal basis. The answer here is quite satisfactory.


\begin{notation} \label{not3.02}
Let $2 \le n \in \bbN$.  Given $\alpha \in \bbT$, we define
\[
\fS_n(\alpha) :=\{ \omega \in \bbT: \omega^{n^k} = \alpha \mbox{ for some } k \in \bbZ \}. \]

\end{notation}

Observe that $\fS_n(1)$ denotes all possible $n^k$-th roots of $1$, and in general, $\fS_n(\alpha)$ contains all $n^j$-th powers and $n^j$-th roots of $\alpha$, $j \ge 1$.   Moreover, if $\alpha,\beta\in \mathbb{T}$, then $\fS_n(\alpha)=\fS_n(\beta)$ if and only if $\fS_n(\alpha)\cap \fS_n(\beta)\neq \varnothing$.
	

\begin{eg} \label{eg3.03}
It is now easy to verify the following ``minimal" examples of unitary operators satisfying the condition $T^{(n)} \simeq T^n$ for some $n \ge 2$ (resp. for all $n\in \mathbb{N}$).

Given a sequence $\theta := (\theta_n)_n$ of complex numbers of absolute value $1$, define the diagonal unitary operator $U_\theta := \ttt{diag} (\theta_n)_n \in \cB(\ell_2)$.   It is well-known that the \textbf{point spectrum} $\sigma_p(U_\theta)$ (i.e. the set of \emph{eigenvalues} of $\theta$) is just the set $\{ \theta_n\}_{n}$.    This last statement also holds for the \textbf{infinite ampliation} $U := U_\theta^{(\infty)}$ of $U_\theta$.

Let $\alpha \in \mathbb{T}$ and $m\in \mathbb{N}$. Let us denote by $\alpha^{\frac{1}{m}}$ the $m$-th root of $\alpha$ \emph{with the smallest argument} in $[0,2\pi)$.

\begin{enumerate}
	\item[(a)]
	Fix $2 \le n$. The set $\{\alpha^{n^k}:k\in \mathbb{Z}\}$ is clearly countable, and so we can find a sequence $(\theta_n)_n$ such that $\{ \theta_n\}_n = \{\alpha^{n^k}:k\in \mathbb{Z}\}$.
	Let $U := U_\theta^{(\infty)}$.   Then $U^n \simeq U^{(n)}$.
	\item[(b)]
	The set $\{\alpha^{m^{\pm 1}}:m\in \mathbb{N}\}$ is also countable, and so we can find a sequence $(\theta_n)_n$ such that $\{ \theta_n\}_n = \{\alpha^{m^{\pm 1}}:m\in \mathbb{N}\}$.
	Let $U := U_\theta^{(\infty)}$.   Then $U^n \simeq U^{(n)}$ for all $n \ge 1$.
\end{enumerate}
\end{eg}


\begin{lem} \label{lem3.04}
If $U$ is a diagonal unitary and $\jnu$-stable for a fixed  $n\geq 2$, then every $\alpha\in \sigma_p(U)$ has infinite multiplicity,
and all its $n^k$-th powers and all its $n^k$-th roots belong to $\sigma_p(U)$ for $k\in \mathbb{N}$.
\end{lem}

\begin{pf}
Since $U\simeq J_n(U)$ implies that $U^n \simeq U^{(n)}$, we need only verify the last statement about the roots.
To do this, it suffices (by induction) to show that if $\alpha \in \sigma_p(U)$, then all $n$-th roots of $\alpha$ are in $\sigma_p(U)$.

Let $V$ be any $n^{th}$ root of $U$. Then $\sigma_p(V)$ contains at least one $n$-th root $\beta$ of $\alpha$.
It follows from Lemma~\ref{lem2.15} that $\omega^j\beta\in \sigma_p(U)$ for all $j$, where $\omega$ is a primitive $n$-th root of 1.
\end{pf}


\begin{prop}\label{prop3.05}
Let $\alpha\in \mathbb{T}$ and $U$ be the diagonal unitary whose point spectrum is $\fS_n(\alpha)$ with uniform infinite multiplicity. Then $U\simeq J_n(U)$.
\end{prop}

\begin{pf}
By Lemma~\ref{lem2.15}, we must show that
\[U\simeq \textup{diag}(V,\omega V,\cdots,\omega^{n-1}V),\]
where $V$ is an $n$-th root of $U$. This follows from the structure of $\fS_n(\alpha)$ and the fact that every eigenvalue of $V$ has infinite multiplicity.
\end{pf}

The diagonal unitary given in Proposition~\ref{prop3.05} turns out to be the building block for all $\jnu$-stable \emph{diagonal} unitaries $U$.  Let us denote it by $U_n(\alpha)$.


\begin{prop}\label{prop3.06}
Let $2 \le n \in \bbN$, and let $U$ be any $\jnu$-stable diagonal unitary operator. Then there exists a (finite or countably infinite) set $\{ \alpha_j\}_j \subseteq \bbT$ such that $U$ is a direct sum
\[\oplus_{j}U_n(\alpha_j).\]
The direct summands are unique $($up to permutation of the terms$)$.

\end{prop}

\begin{pf}
For any $\alpha\in \sigma_p(U)$, it follows from Lemma~\ref{lem3.04} that $\mathfrak{S}_n(\alpha)\subset \sigma_p(U)$. Thus $U_n(\alpha)$
is a summand of $U$.

If $\sigma_p(U)= \mathfrak{S}_n(\alpha)$, we are done. Otherwise, pick a $\beta\in
\sigma_p(U)\setminus \mathfrak{S}_n(\alpha)$. We have seen that
$\mathfrak{S}_n(\alpha)\cap \mathfrak{S}_n(\beta)=\varnothing$, and $U_n(\beta)$ is another summand of $U$.
Now it is easy to take a maximal subset $\{\alpha_1,\alpha_2,\cdots\}$ of $\sigma_p(U)$
such that the sets $\mathfrak{S}_n(\alpha_j)$ are disjoint and their union is $\sigma_p(U)$. Since
every member $\alpha$ of $\sigma_p(U)$ uniquely determines $\mathfrak{S}_n(\alpha)$, we obtain the desired decomposition.
\end{pf}


\subsection{} \label{sec3.07}
We now consider a stronger assumption on a diagonal unitary operator $U$, namely:   when is $U$ $\jnu$-stable for all $n\in \mathbb{N}$? To answer this question, we need only a
slight modification of our notation. If $\alpha$ is in $\sigma_p(U)$ for such a diagonal unitary, then all $n$-th powers and all $n$-th roots of $\alpha$ are in $\sigma_p(U)$. Furthermore if $\omega$ is any primitive $n$-th root of 1, $\omega\alpha\in \sigma_p(U)$. So, for $\alpha\in \mathbb{T}$, let us denote by $\fS(\alpha)$ the set of all $\omega\lambda$, where $\omega$ is any
$n$-th root of 1, and $\lambda$ is any $m$-th power or any $j$-th root of $\alpha$.

With a reasoning similar to what we used in the case of fixed $n$, we can verify that $\fS(\alpha)=\fS(\beta)$ if and only if
$\fS(\alpha)$ and $\fS(\beta)$ intersect. Now the following proposition is easy to prove.


\begin{prop}\label{prop3.08}
Let $U$ be a diagonal unitary operator and suppose that $U$ is $\jnu$-stable for all $n\geq 2$. Then $U$ is a $($finite or infinite$)$ direct
sum $\oplus_{j}U_j$, where for each $j$, $\sigma_p(U_j)=\fS(\alpha_j)$ for some $\alpha_j\in \mathbb{T}$ with uniform infinite multiplicity.
\end{prop}


\bigskip

\begin{eg} \label{eg3.09}
As a corollary to Proposition~\ref{prop3.06}, we have the following:   there exist $U \in \bofh$ unitary and $\theta \in \bbT$ such that $U$ is $\ensuremath{\textsc{j}_2^{[\textsc{u}]}}$-stable, but $\theta U$ is not $\ensuremath{\textsc{j}_2^{[\textsc{u}]}}$-stable.

Let $U = I \oplus D$, where $D = \ttt{diag} (d_n)_n$, and where $\{ d_n\}_n$ consists of all $2^k$-roots of unity repeated with infinite multiplicity.   By Proposition~\ref{prop3.06}, $U$ is $\ensuremath{\textsc{j}_2^{[\textsc{u}]}}$-stable.  Let $\theta = e^{2 \pi i/3}$, so that $\theta$ is a third root of unity.

Then $V := \theta U = \theta I \oplus \theta D$ is diagonal, and $V$ has $\theta$ as an eigenvalue.  As such, for $V$ to be $\ensuremath{\textsc{j}_2^{[\textsc{u}]}}$-stable we require that $\theta^{\frac{1}{2}} = e^{2 \pi i/6}$ should also be an eigenvalue of $V$.  But then
\[
e^{2 \pi i/ 6} = \omega_k^{j} e^{2 \pi i/2^k} e^{2 \pi i/3} \mbox{     for some }  k \ge 1\mbox{ and } 1 \le j \le 2^k, \]
where $\omega_k$ is a primitive $2^k$-th root of unity.

Raising both terms to the power $2^k$, we see that
\[
e^{2^k (2\pi i)/6} = e^{2^k (2 \pi i)/3}, \]
so that $e^{2^k (2 \pi i)/6}  = 1$, whence $2^k (2 \pi i)/6 \in 2 \pi i \, \bbZ$.

Thus $2^k/6 \in  \bbZ$, or equivalently $2^k \in 6 \bbZ$, which is clearly impossible.
\end{eg}

%
%


\subsection{} \label{sec3.11}
Having described which \emph{diagonal} unitary operators $U$ are $\jnu$-stable for some $2 \le n \in \bbN$, we now investigate the $\jnu$-stability of unitary operators whose spectral measures are non-atomic.

\smallskip

An operator $T\in \cB(\hilb)$ is said to be \textbf{cyclic} if there exists $0\neq x\in \hilb$ such that $\hilb=\overline{\textup{span}}\{T^nx:n\geq 0\}$, where $T^0:=I$.
If $U\in \cB(\hilb)$ is a \emph{cyclic} unitary operator, then --
by the Spectral Theorem for unitary operators  -- there exists a Borel probability measure $\nu$ with support $\sigma(U)$ such that
$U$ is unitarily equivalent to the multiplication operator $M_{z,\nu}$ on $L^2(\sigma(U),\nu)$, where
\[M_{z,\nu}f(z)=zf(z),~~z\in \mathbb{T},f\in L^2(\sigma(U),\nu).\]
If, furthermore, $U$ is   $\jnu$-stable for some $n\geq 2$, then  by Proposition~\ref{prop2.07}, $\sigma(U)$ is equal to $\bbT = \{ z \in \bbC: |z| = 1\}$.


\begin{thm}\label{measure}
Let $U$ be a cyclic unitary operator with $\sigma(U) = \bbT$ and $\nu$ be the Borel probability measure corresponding to $U$ as above.  Suppose that $\nu$ is a non-atomic  measure. Let  $2\leq n\in \bbN$.  Set $\bbT_0 := \{e^{2\pi i\theta}: \theta\in [0,\frac{1}{n}]\}$, and for $1 \le k \le n-1$, set $\bbT_k = e^{2 \pi i \frac{k}{n}} \bbT_0$.

Given a Borel set $X \subseteq \mathbb{T}_0$, define $\nu_k(X) :=\nu(e^{2\pi i \frac{k}{n}}\cdot X), k=0, 1,\cdots,n-1$, where $\alpha \cdot X := \{ \alpha x: x \in X\}$.

Then $U$ is $\jnu$-stable if and only if $\nu$ satisfies the following conditions.
\begin{enumerate}
	\item[(a)]	
	The Borel measures $\nu_0, \nu_1, \ldots, \nu_{n-1}$ on $\mathbb{T}_0$ are  mutually absolutely continuous with each other.
	\item[(b)]
	The measure $\nu_0$ is mutually absolutely continuous with $\varphi$, where $\varphi$  is a continuous measure with $\mathrm{supp}\, \varphi = \mathbb{T}_0$ such that for any Borel set $Y \subseteq \mathbb{T}_0$,
	\[
	\varphi(Y) := \nu (Y^n), \]
	where $Y^n := \{ y^n: y \in Y\}$.
\end{enumerate}
\end{thm}

\begin{pf}
Since $\jnu$-stability is clearly invariant under unitary equivalence, we may assume without loss of generality that $U=M_{z, \nu}$ acting on $L^2(\mathbb{T},\nu)$.
Let us denote by $\alpha^{\frac{1}{n}}$ the $n$-th root of $\alpha$ \emph{with the smallest argument} in $[0,2\pi)$.
Then $V:=M_{z^{\frac{1}{n}}, \nu}$ is a $n^{th}$ root of $U$.
Consider the map $W: L^2(\mathbb{T},\nu)\to L^2(\mathbb{T}_0,\varphi)$ which sends $g(z)$ to $g(z^n)$. It is straightforward to check that $W$
is a unitary operator which implements the unitary equivalence of $V$ and $M_{z,\varphi}$, where $M_{z,\varphi}$ is the multiplication operator  on $L^2(\mathbb{T}_0,\varphi)$.

By Lemma~\ref{lem2.15}, $U$ is $\jnu$-stable if and only if
\[U\simeq \textup{diag}(V,\omega V,\cdots,\omega^{n-1}V),\]
where $\omega = e^{\frac{2 \pi i}{n}}$ is a primitive $n$-th root of 1. Note that $\mathbb{T}_k:=\omega^k \mathbb{T}_0$, $k=0,1,\cdots, n-1$.
Let $\varphi_k$ be the Borel measure on $\mathbb{T}$ defined by $\varphi_k(Y)=\varphi(\omega^{-k}\cdot (Y\cap \mathbb{T}_k))$, where $Y\subseteq \mathbb{T}$ is a Borel set, and define $\eta=\sum_{k=0}^{n-1}\varphi_k$.
Then $\eta$ is a Borel measure on $\mathbb{T}$.
Note that $\textup{diag}(V,\omega V,\cdots,\omega^{n-1}V)$ is
unitarily equivalent to the multiplication operator $M_{z,\eta}$ on $L^2(\mathbb{T},\eta)$, where $z\in \mathbb{T}$.
So $U$ is $\jnu$-stable if and only if $M_{z,\nu}$ is unitarily equivalent to $M_{z,\eta}$; that is, if and only if $\nu$ is mutually absolutely continuous with $\eta$.

\end{pf}



\subsection{A measure-theoretic interlude.}\label{sec3.01.01}\ \ \ The above shows that the problem of characterising $\jnu$-stable cyclic non-atomic unitary operators is purely a measure-theoretic one.
Clearly, normalised Lebesgue measure on $\mathbb{T}$ satisfies the measure-theoretic assumptions of Theorem \ref{measure} (and this, incidentally, yields an alternative way of showing that the bilateral shift operator is $\jnu$-stable, $n\geq 2$).   Our next goal is to show that Lebesgue measure is far from being unique in this respect.  We next  describe a construction of a continuum of \textbf{mutually singular} such measures based on self-similar measure theory (c.f. Theorem \ref{thm:main}).  As a consequence, we conclude that there exists a continuum of unitary operators, no two of which are unitarily equivalent, yet each of which is $\jnu$-stable.

We first need to recall some measure-theoretic notions.
Let $X$ (resp. $Y$) be a metric space, and denote by $\mathfrak{Bor}(X)$ (resp. $\mathfrak{Bor}(Y)$) the $\sigma$-algebra of Borel sets of $X$ (resp. of $Y$). Given two Borel measures  $\mu$ and $\nu$, we write $\mu\asymp \nu$ to indicate that $\mu$ and $\nu$ are mutually absolutely continuous with respect to each other.
Let $f$ be a Borel measurable map from $X$ to $Y$.
The \textbf{pushforward} $f_* \mu$ of a Borel measure $\mu$ on $X$ under the map $f$ is the Borel measure on $Y$ defined by $(f_*\mu)(B)=
\mu(f^{-1}(B))$ for all $B\in \mathfrak{Bor}(Y)$.

We also bring to the attention of the reader the fact that in measure-theory parlance, one says that a sequence $(\mu_n)_n$  of Borel measures on a metric space $X$ is said to \textbf{converge weakly} to the Borel measure $\mu$ on $X$ if
\[
\lim_n \mu_n(f) := \lim_n \int_X f \, d\mu_n = \int_X f\, d\mu =: \mu(f) \]
for all $f \in \cC_b(X, \bbC) := \{ f: X \to \bbC \mid f \mbox{ is continuous and bounded}\}$.   (We mention this because if $X$ is compact, this is what functional analysts normally refer to as \textbf{weak}${}^*$-\textbf{convergence}.)


\bigskip

The following result is known as the \textbf{Portmanteau Theorem}.

\smallskip

\begin{thm}\cite[Theorem 2.2.5,~Corollary 2.2.6,~Theorem 2.4.1]{Bogachev2018}\label{equivalent conditions for weak convergence}
Suppose that we are given a sequence of finite Borel measures $(\mu_n)_n$ and a finite Borel measure $\mu$ on a metric space $X$.
The following conditions are equivalent:
\begin{enumerate}
	\item[(i)]
	the sequence $(\mu_n)_n$ converges weakly to $\mu$;
	\item[(ii)]
	for every closed set $F$ in $X$, one has
	\[
	\underset{n\to \infty}{\lim\sup}\, \mu_n(F)\leq \mu(F),\]
	and $\underset{n\to \infty}{\lim}\mu_n(X)=\mu(X)$;
	\item[(iii)]
	for every open set $U$, one has
	\[
	\underset{n\to \infty}{\lim\inf}\mu_n(U)\geq \mu(U),\]
	and $\underset{n\to \infty}{\lim}\mu_n(X)=\mu(X)$;
	\item[(iv)]
	for every bounded, upper semicontinuous function $f:X \to \bbR$, one has
	\[
	\underset{n\to \infty}{\lim\sup}\int_{X}fd\mu_n\leq \int_X fd\mu,\]
	and $\underset{n\to \infty}{\lim}\mu_n(X)=\mu(X)$;
	\item[(v)]
	for every Borel set $E$ of $X$ satisfying $\mu(\overline{E}\setminus (\mathrm{int}\, E))=0$, one has
\[\underset{n\to \infty}{\lim}\mu_n(E)= \mu(E).\]
\end{enumerate}
\end{thm}


Our next goal is to obtain a continuum of non-atomic measures on $\mathbb{T}$ which fulfill the assumptions of Theorem \ref{thm:main},  no two of which are  mutually absolutely continuous with respect to one another.   In so doing, we obtain a continuum of $\textsc{j}_N^{[\textsc{u}]}$-stable cyclic non-atomic unitary operators.  In fact, there is no loss of generality in replacing the circle $\bbT$ by the interval $[0,1]$, since one then simply considers the associated family of pushforwards of the measures on $[0,1]$ by the continuous function  $g: [0,1] \to \bbT$ defined by $g(t)=e^{2\pi i t}$.


\begin{thm}\label{thm:main}
  Let $N \geq 3$.   There exist an uncountable set $\Omega$ and a continuum of non-atomic Borel probability
  measures $(\mu_\alpha)_{\alpha\in \Omega}$, each with support equal to $[0,1]$, such that
\begin{itemize}
	\item[(i)]
    	$\mu_\alpha(A+\tfrac iN) = \mu_\alpha(A)$ for all $A\in \cB([0,\frac{1}{N}])$ and $1\leq i\leq N-1$.
    	\item[(ii)]
	$\mu_\alpha|_{[0,\tfrac 1N]}$ is mutually absolutely continuous to $f_*\mu_\alpha$, where $f(x) = \tfrac x N$.
    	\item[(iii)]
	If $\alpha \ne \beta \in \Omega$, then  $\mu_\alpha$ and  $\mu_\beta$ are  mutually singular.
\end{itemize}
\end{thm}

The proof is constructive and we will exhibit a method for constructing  a  measure satisfying conditions (i) and (ii) of the above theorem as the weak limit of a  sequence of specifically defined measures.   First, we shall require a number of lemmas, definitions and notations.



\subsection{} \label{paragraph2.5}
Fix $N\geq 3$.  We start by constructing a probability measure $\nu$ on $[0,1]$ as follows.

\bigskip

Given $0 \le K \in \bbZ$, we write  $\mathtt{i}=(i_1 \, i_2 \, \cdots \, i_K)$ to denote a word of length $K$ over $\{0,1,\dots,N-1\}$, with $\varnothing$ being the unique (empty) word of length zero.
Write $\Sigma_K =
\{0,1,\dots,N-1\}^K$ for the set of \emph{all} words of length equal to $K$ and $\Sigma_* =
\bigcup_{K=0}^\infty \Sigma_K$ for all words of finite length over $\{ 0, 1, \ldots, N-1\}$.
Given $\mathtt{i}$, we denote by $| \mathtt{i} |$ the \textbf{length} of $\mathtt{i}$;  that is, $|\mathtt{i}| = K$ precisely when $\mathtt{i} \in \Sigma_K$.
If $\mathtt{i} = (i_1\, i_2\, \cdots\, i_K)$ and $\mathtt{j} = (j_1\, j_2\, \cdots\, j_L)$ for some $K, L \ge 0$, we write $\mathtt{i j}$ to denote the concatenation of these two words, namely:   $\mathtt{ij} = (i_1\, i_2\, \cdots\, i_K\, j_1\, j_2\, \cdots\, j_L) \in \Sigma_{K+L}$.
If $\mathtt{j} = (j)$ is a word of length one, we may also write $\mathtt{i}j$ instead of $\mathtt{i j}$.
Let $\mathtt{b} \in \Sigma_*$, and suppose that we may write $\mathtt{b} = \mathtt{i} \mathtt{j}$, where $\mathtt{i}, \mathtt{j} \in \Sigma_*$.
We shall refer to $\mathtt{i}$ (resp.\ to $\mathtt{j}$) as a \textbf{prefix} (resp.\ as a \textbf{suffix}) of $\mathtt{b}$.
A \textbf{subword} of $\mathtt{b}$ is any $\mathtt{k} \in \Sigma_*$   for which there exist $\mathtt{i}, \mathtt{j} \in \Sigma_*$ (possibly of length zero) such that $\mathtt{b} = \mathtt{i} \mathtt{k} \mathtt{j}$.

If $\mathtt{i}= (i_1 \, i_2 \, \cdots i_K)  \in \Sigma_K$ and $0 \le n \le K$, we denote by $\mathtt{i}|_n = (i_1 \, i_2 \, \cdots i_n) \in \Sigma_n$, the truncation of $\mathtt{i}$ to the first $n$ coordinates.   If $n = K-1$, we also concisely write $\mathtt{i}^-$ in place of $\mathtt{i}|_{K-1}$.

To each $\mathtt{i} \in \Sigma_*$ we may associate an $N$-ary rational number via the natural projection
\[
\begin{array}{rccc}
  \pi: & \Sigma_* & \to & [0,1] \\[.5em]
  &\mathtt{i}=(i_1 \, i_2 \, \cdots i_K) &\mapsto &\sum_{k=1}^K \frac{i_k}{N^k}.
\end{array} \]		
Note that
\[
  \pi(\Sigma_K) = \bigg\{\frac{m}{N^K} : m \in\{0,1,\dots,N^K -1\}\bigg\}
\]
consists of all rational numbers in $[0,1)$  with
denominator $N^K$.

The projection $\pi:\Sigma_*\to[0,1]$ has a natural extension to infinite
sequences $\mathtt{i}\in \Sigma_\infty = \{0,1,\dots,N-1\}^{\bbN}$ by setting
\[
  \pi( (i_1\, i_2 \, \dots \, i_n \, \dots)) = \sum_{k=1}^{\infty} \frac{i_k}{N^k}.
\]
This extension is a surjection from $\Sigma_\infty$ to $[0,1]$ that fails to be injective.
However,
there are only countably many $x\in [0,1]$ such that $\pi^{-1}(x)$ is not a singleton.
These points
are of the form $(i_1 \, \dots \, i_{n-1} \,i_n\, (N-1)\,(N-1)\,(N-1)\,\dots)$ and $(i_1\,\dots\,i_{n-1}\,(i_n+1)\,0\,0\,\dots)$.

\bigskip
For $i \in \bbN$, define $\mathtt{b}^i = (1\overbrace{0\,0\,\cdots 0}^{i \text{ times }}2) \in \Sigma_{i+2}$, and let $B_0 := \{ \mathtt{b}^1, \mathtt{b}^2, \mathtt{b}^3, \ldots \}$, so that
\[
B_0 := \{ (1\,0\,2), (1\,0\,0\,2), (1\,0\,0\,0\,2), \ldots \} \subseteq \Sigma_*. \]
The following key properties\footnote{We note that the arguments that follow on the next pages can be generalised to any $B$ that satisfies these four criteria with little effort, though substantially more notation. The argument can also be extend to much more general sets of words $B$, but we will not give details here. The interested reader may peruse texts on sub-self-similar sets and their measures, e.g. \cite{Falconer1995}.}  of $B_0$ will be of use to us below:

\begin{itemize}
	\item{}
	$3\leq |\mathtt{b}^i| < |\mathtt{b}^{i+1}|$ for all $i \ge 1$.
	\item{}
	The first digit of each $\mathtt{b}^i$ is $1$, the last digit of each $\mathtt{b}^i$ is $2$, and every other digit of each $\mathtt{b}^i$ is $0$.
	\item{}
	  Given $1 \le i \ne j$, no non-trivial prefix of $\mathtt{b}^i$ is a suffix of $\mathtt{b}^j$.
	\item{}
	  Given $1 \le i$, no non-trivial proper prefix of $\mathtt{b}^i$ is a suffix of $\mathtt{b}^i$.
      \end{itemize}
We shall refer to elements of $B_0$ as  ``forbidden" words, and we note that no word $\mathtt{b}^i$ is a subword of $\mathtt{b}^j$ (unless $i = j$).

Of course, if $B \subseteq B_0$, and if we order the elements of $B$ by length,  then the same four properties hold for the elements of $B$.

\bigskip

We now consider an arbitrary but fixed subset $B \subseteq B_0$.
Our goal is to construct a non-atomic probability measure $\nu$ (which is uniquely determined by $B$) with support $[0,1]$, and to associate to $\nu$ a push-forward measure $\mu_0 := f_* \, \nu$ under the function $f(x) = x/N$, $x \in [0,1]$.   The measure $\mu_0$ is then used to construct a measure satisfying conditions (i) and (ii) of Theorem~\ref{thm:main}.  The dependence of this construction upon the initial choice of $B \subseteq B_0$ is what will guarantee that condition (iii) of Theorem~\ref{thm:main} is met.

The construction of $\nu$ is, alas, somewhat elaborate.   It will be defined as a weak limit of convex combinations of certain non-atomic probability measures $L_{\pi([\mathtt{i}])}$ defined on subintervals of $[0,1]$.   However, in order to define the probability measures $L_{\pi([\mathtt{i}])}$, we shall  first construct a non-atomic measure $\widetilde{\nu}$ whose support is a Cantor-like subset $\cC$ (also uniquely determined by $B$) of $[0,1]$ which avoids the forbidden words appearing in $B$.   As we shall soon see, this measure $\widetilde{\nu}$ is itself a weak limit of a sequence of weighted combinations of non-atomic probability measures $\widetilde{L}_{\pi([\mathtt{i}])}$ defined on subintervals of $[0,1]$.   In fact, the  measures $\widetilde{L}_{\pi([\mathtt{i}])}$  are nothing more than the usual Lebesgue measure $\cL$ restricted to the subinterval $\pi([\mathtt{i}])$ and then normalised.  We will then define $L_{\pi([\mathtt{i}])} = {F_{\mathtt{i}}} _* \widetilde{\nu}$ as the push-forward measure of $\widetilde{\nu}$ under the linear mapping $F_{\mathtt{i}}(x) = N^{-|\mathtt{i}|} x + \pi(\mathtt{i})$.

In order to improve the readability of the arguments below, we shall omit the subscript $B$ from our sets, functions and measures until the proof of Theorem~\ref{thm:main}.

\bigskip

We begin by describing the weights that will appear in the linear combinations of probability measures we shall require to obtain $\widetilde{\nu}$ and subsequently $\nu$.

\bigskip

Given $\mathtt{i} = (i_1 \, i_2\, \cdots \, i_K) \in\Sigma_*$,  we define
\[
B(\mathtt{i})
	= 	\begin{cases} \{ (i_1 \, i_2\, \cdots \, i_K \, 2) \} & \text{ if some } \mathtt{b} \in B  \text{  is a suffix of } (i_1 \, i_2\, \cdots \, i_K \, 2) \\
		\varnothing  & \text{ otherwise.}
		\end{cases} \]
Observe that $B(\mathtt{i}^-) = \{ \mathtt{i}\}$ precisely if $\mathtt{b}$ is a suffix of $\mathtt{i}$ for some $\mathtt{b} \in B$.
A necessary (but very insufficient) condition to have $B(\mathtt{i}^-) = \{ \mathtt{i}\}$ is that $i_K = 2$.

\bigskip

Given a set $Z$, let $\# Z$ denote the cardinality of $Z$.
We define weight functions $p, \widetilde{p}:\Sigma_*\to\bbR^+_0:=\{t:0\leq t\in \mathbb{R}\}$  (which again depends upon $B$) inductively by
\begin{align*}
  p(\varnothing) = \widetilde{p}(\varnothing)= 1
  \quad\text{,}\quad
  p(\mathtt{i}) &= 	
  \begin{cases}
    p(\mathtt{i}^-)2^{-|\mathtt{i}|} &\text{if }  B(\mathtt{i}^-) = \{ \mathtt{i} \} \\
    p(\mathtt{i}^-)\frac{1-\#B(\mathtt{i}^-)\cdot2^{-|\mathtt{i}|}}{N-\#B(\mathtt{i}^-)}&\text{otherwise.}
  \end{cases}
  \intertext{and}
  \widetilde{p}(\mathtt{i}) &= 	
  \begin{cases}
    0 &\text{if }  B(\mathtt{i}^-) = \{ \mathtt{i} \} \\
    \widetilde{p}(\mathtt{i}^-)\frac{1}{N-\#B(\mathtt{i}^-)}&\text{otherwise.}
  \end{cases}
\end{align*}
Writing
\[
p_{n}(\mathtt{i}|_n)
	= \begin{cases} 2^{-n}&\text{if } B(\mathtt{i}|_{n-1}) = \{ \mathtt{i}|_n \}\\
    	\frac{1-\#B(\mathtt{i}|_{n-1})\cdot2^{-n}}{N-\#B(\mathtt{i}|_{n-1})}&\text{otherwise},
 	 \end{cases} \]
we obtain
\[
p(\mathtt{i}) = \prod_{k=1}^{|\mathtt{i}|}p_{k}(\mathtt{i}|_k).
\]
A similar formula holds for $\widetilde{p}(\mathtt{i})$, replacing $2^{-n}$ with $0$.

The weight functions are defined in such a way to allow for the following consistency relation.
Let $\mathtt{i} \in \Sigma_K$.
If no element of $B$ is a suffix of $\mathtt{i}2$, then $B((\mathtt{i}j)^-) = B(\mathtt{i}) = \varnothing$ for all $0 \le j \le N-1$, whence
\begin{equation}\label{eq:splittingup}
	\sum_{j=0}^{N-1} p( \mathtt{i} \, j) =\sum_{j=0}^{N-1} p(\mathtt{i})\frac{1}{N} = p(\mathtt{i}),
\end{equation}
whereas if there exists $\mathtt{b} \in B$ which is a suffix of $\mathtt{i}2$, then
\begin{equation}\label{eq:splittingup2}
	\sum_{j=0}^{N-1} p( (\mathtt{i} \, j))
		=\sum_{j\in\Sigma_1\setminus\{ 2 \}} p(\mathtt{i})\frac{1-2^{-(K+1)}}{N-1} + p(\mathtt{i})2^{-(K+1)} = p(\mathtt{i}).
\end{equation}
A similar consistency relation can be derived for $\widetilde{p}$.
Heuristically, the weights $p(\mathtt{i})$ an $\widetilde{p}(\mathtt{i})$ of any
word $\mathtt{i}\in\Sigma_*$ get split amongst its successor words $\mathtt{i}j$.

Let $[\mathtt{i} ]$ be the set of (infinite) codings with initial sequence $\mathtt{i}$ followed by arbitrary codings,
\[
  [ \mathtt{i} ] =\big\{ \mathtt{j}\in\Sigma_\infty : i_k = j_k \text{ for all }k\leq
  |\mathtt{i}| \big\}.
\]
The set $[ \mathtt{i} ]$ is clearly independent of $B$ and
it is useful to keep in mind that $\pi([ \mathtt{i}])$ is a closed subinterval of $[0,1]$.

\bigskip


\subsection{} \label{sec3.16.new}
Suppose that we are given non-atomic probability measures $\lambda_{\pi([\mathtt{i}])}$ on the the interval $\pi([\mathtt{i}]):=[\pi(\mathtt{i}),\pi(\mathtt{i})+N^{-|\mathtt{i}|}]$, $\mathtt{i} \in \Sigma_*$.   Then  for each $K \ge 1$, one may define the probability measures
\[
\nu_{K, \lambda} := \sum_{\mathtt{i} \in \Sigma_K} p(\mathtt{i}) \lambda_{\pi{([\mathtt{i}])}}
\mbox{\ \ \ \ \ \ \ \ \ \  and \ \ \ \ \ \ \ \ \ }
\widetilde{\nu}_{K, \lambda} := \sum_{\mathtt{i} \in \Sigma_K} \widetilde{p}(\mathtt{i}) \lambda_{\pi{([\mathtt{i}])}}  \]

Heuristically, these measures split the measure at coding $\mathtt{i}\in\Sigma_K$ into its ``children''
$\mathtt{i} \mathtt{j}$ ($\mathtt{j} \in \Sigma_1$) evenly, giving any forbidden word (of which there is at most one) exponentially small weight (in the case of $\nu_{K, \lambda}$) or zero weight (in the case of $\widetilde{\nu}_{K, \lambda})$.
Since $\nu_{K, \lambda}$ and $\widetilde{\nu}_{K, \lambda}$ are  non-atomic and $\pi([\mathtt{i}])\cap \pi([\mathtt{j}])$
contains at most one point if $\mathtt{i},\mathtt{j}\in\Sigma_K$ are distinct, we have
$\nu_{K, \lambda}(\pi([\mathtt{i}]))= p(\mathtt{i})$ and $\widetilde{\nu}_{K, \lambda} (\pi([\mathtt{i}]))= \widetilde{p}(\mathtt{i})$
for all $K\geq |\mathtt{i}|$, independent of our choice of $\lambda_{\pi([\mathtt{i}])}$.
Perhaps surprisingly, the weak limit of $(\nu_{K,\lambda})_K$  and of $(\widetilde{\nu}_{K, \lambda})_K$ not only exist, but can be shown to be independent of the actual choice of $\lambda_{\pi([\mathtt{i}])}$ (even without the assumption that each $\lambda_{\pi([\mathtt{i}])}$ is non-atomic).  We shall not require this, however, and so the proof is omitted.
	

\bigskip

The following estimates on  the size of $p(\mathtt{i})$ and $\widetilde{p}(\mathtt{i})$ will prove useful below.

\begin{lem}\label{thm:weightDecay}
  Let $\mathtt{i}\in\Sigma_*$ and let $p(\mathtt{i})$ and $\widetilde{p}(\mathtt{i})$ be as above.
  Then, there exists a constant $C_N$ only depending on $N$ such that
  \[
    p(\mathtt{i}) \leq C_N \cdot(N-1)^{-|\mathtt{i}|}
    \quad
    \text{and}
    \quad
    \widetilde{p}(\mathtt{i}) \leq (N-1)^{-|\mathtt{i}|}.
  \]
  In particular, $p(\mathtt{i})\to 0$ and $\widetilde{p}(\mathtt{i})\to0$ as $|\mathtt{i}|\to\infty$.

\end{lem}
\begin{pf}
  By definition,
  $p(\mathtt{i}) = \prod_{k=1}^{|\mathtt{i}|} p_k (\mathtt{i}|_k)$ and
  for $1\leq k \leq |\mathtt{i}|$ we have
  $p_{k}(\mathtt{i}|_{k})\leq \max \left( 2^{-k}, \tfrac{1}{N-1} \right)$.
  Notice that $2^{-k} < \tfrac{1}{N-1}$ for $k> \log_2(N-1)$ and
  $2^{-k} \geq \tfrac{1}{N-1}$ for $k \leq \log_2(N-1)$.
  Now $2^{-k}\leq \tfrac12$ and so $2^{-k}\leq \tfrac{N-1}{2}\tfrac{1}{N-1}$.
  Thus,
  \[
    \prod_{k=1}^{|\mathtt{i}|} p_k (\mathtt{i}|_k)
    \leq \left( \tfrac{N-1}{2} \right)^{\log_2(N-1)} \prod_{k=1}^{|\mathtt{i}|}\tfrac{1}{N-1}
    =\left( \tfrac{N-1}{2} \right)^{\log_2(N-1)} \left(\tfrac{1}{N-1}\right)^{|\mathtt{i}|}.
  \]
  The conclusion for $\widetilde{p}$ follows as $p_k(\mathtt{i}|_k)\leq (N-1)^{-1}$ for all $k$.
\end{pf}
	

\subsection{}
Keeping $B \subseteq B_0$ and the family $\{ \lambda_{\pi([\mathtt{i]})} : \mathtt{i} \in \Sigma_*\}$ fixed as above, our next goal is to show that the sequences $(\nu_{n, \lambda})_n$ and $(\widetilde{\nu}_{n,  \lambda})_n$
of probability measures converge weakly to non-atomic probability measures $\nu_{\lambda}$ and $\widetilde{\nu}_{\lambda}$, respectively.

The Kantorovich-Rubinstein metric (also known as the Wasserstein metric) $W_1(\cdot, \cdot)$ is a metric on the set of probability measures on a metric space~\cite{Wasserstein1969}.
We shall not need the general definition of this metric.
For our purposes, we only need to know the following two things.
\begin{itemize}
	\item{}
	Convergence of the sequence $(\nu_{n})_n$ with respect to the Kantorovich-Rubinstein distance implies weak convergence of that sequence (see e.g.\ \cite[Theorem 7.12]{Villani2003}).
	\item{}
	If $\mu$ and $\nu$ are probability measures, there is the following equivalent formulation of the Kantorovich-Rubinstein metric  in terms of Lipschitz functions  $\textup{Lip}_1([0,1],\bbR)$ with constant $1$ (see \cite[Equation (7.1)]{Villani2003}):
	\[
 	 W_1(\mu,\nu) =
	 \sup_{h\in\textup{Lip}_1([0,1],\bbR)}\left( \int h(x) d\mu(x) - \int h(x) d\nu(x)\right).
		      \]
\end{itemize}			
	
In the next few results, we shall require the Cantor-like set:
\[
	\cC := \{ \mathtt{i} \in \Sigma_\infty:   \text{there does not exist } \mathtt{b} \in B \text{ such that } \mathtt{b} \text{ is a subword of } \mathtt{i}\}. \]


\begin{lem} \label{lemmasmallpnew}
  Let $B \subseteq B_0$.  Given non-atomic probability measures $\lambda_{\pi([\mathtt{i}])}$ on the  interval $\pi([\mathtt{i}]):=[\pi(\mathtt{i}),\pi(\mathtt{i})+N^{-|\mathtt{i}|}]$, $\mathtt{i} \in \Sigma_*$, let $(\nu_{n, \lambda})_{n=1}^\infty$ and $(\widetilde{\nu}_{n, \lambda})_{n=1}^\infty$  denote the corresponding sequences of probability measures constructed as above.
  Then the weak limits $\nu_\lambda  =\lim_n \nu_{n, \lambda}$ and $\widetilde{\nu}_\lambda =\lim_n \widetilde{\nu}_{n, \lambda}$ exist
  and are non-atomic probability measures.
  Further, $\nu_\lambda$ is fully supported on $[0,1]$, whereas $\widetilde{\nu}_\lambda$ is fully supported on $\pi(\mathcal{C})\subseteq[0,1]$.
\end{lem}

\begin{pf}
  The convergence proofs for $\nu_\lambda$ and $\widetilde{\nu}_\lambda$ are identical and we will only show the former.  To aid the exposition, for the remainder of this Lemma we shall abbreviate the notation $\nu_{n, \lambda}$ to $\nu_n$, $n \ge 1$.

Recall that the space of probability measures on $[0,1]$ is itself compact with respect to the
weak-$*$ topology it inherits as a subset of the dual space $(\cC([0,1], \bbC))^*$ of the space $\cC([0, 1], \bbC)$ of continuous, complex-valued functions on $[0,1]$.
In particular, it is complete in that topology.
Furthermore, since $[0,1]$ is compact, measure-theoretic weak convergence of these measures (as defined in paragraph~\ref{sec3.01.01}) coincides with weak-$*$ convergence in the functional-analytic sense.

Thus to prove that $(\nu_{n})_n$ converges weakly to some measure $\nu_\lambda (= \nu_{B, \lambda})$ on $[0,1]$, we see from above that it suffices to show that it is Cauchy with respect to the Kantorovich-Rubinstein metric.
To that end, we first calculate $W_1(\nu_{n},  \nu_{n+1})$, $n \ge 1$.

Fix $K \ge 1$ and consider the two measures $\nu_K,\nu_{K+1}$.
Let $h\in\textup{Lip}_1([0,1],\bbR)$.
The Kantorovich-Rubinstein distance between $\nu_K$ and $\nu_{K+1}$ becomes
\begin{align*}
  W_1(\nu_K,\nu_{K+1})
  &=
  \sup_{h\in\textup{Lip}_1([0,1],\bbR)}\left( \int h(x)d\nu_{K}(x)- \int h(x)d\nu_{K+1}(x) \right)
  \\
  &=
  \sup_{h\in\textup{Lip}_1([0,1],\bbR)}\left( \sum_{\mathtt{i}\in\Sigma_K}\int_{\pi([\mathtt{i}])} h(x)d\lambda_{\pi([\mathtt{i}])}(x)
    -  \sum_{\mathtt{i}\in\Sigma_{K+1}}\int_{\pi([\mathtt{i}])} h(x)d\lambda_{\pi([\mathtt{i}])}(x)
 \right)
  \\
  &=
  \sup_{h\in\textup{Lip}_1([0,1],\bbR)}\sum_{\mathtt{i}\in\Sigma_K}
  \left( \int_{\pi([\mathtt{i}])} h(x)d\lambda_{\pi([\mathtt{i}])}(x)
    -  \sum_{\mathtt{j}\in\Sigma_1}\int_{\pi([\mathtt{ij}])} h(x)d\lambda_{\pi([\mathtt{ij}])}(x)
 \right)
  \\
  &\leq
  \sum_{\mathtt{i}\in\Sigma_K}
  \sup_{h\in\textup{Lip}_1([0,1],\bbR)}
  \left( \int_{\pi([\mathtt{i}])} h(x)d\lambda_{\pi([\mathtt{i}])}(x)
    -  \sum_{\mathtt{j}\in\Sigma_1}\int_{\pi([\mathtt{ij}])} h(x)d\lambda_{\pi([\mathtt{ij}])}(x)
  \right)
  \\
  &\leq
  \sum_{\mathtt{i}\in\Sigma_K}\sup_{h\in\textup{Lip}_1([0,1],\bbR)} \left( h(x(\mathtt{i}))p(\mathtt{i})
  - \sum_{\mathtt{j}\in\Sigma_1}^{N-1}h(x'(\mathtt{ij})) p(\mathtt{ij} ) \right)
  \intertext{for some $x(\mathtt{i}) \in \pi([\mathtt{i}])$ and $x'(\mathtt{ij})\in \pi([\mathtt{ij}])$.}
  \end{align*}
  Since $\sum_{\mathtt{j}\in\Sigma_1} p(\mathtt{ij} ) = p(\mathtt{i})$,
 \begin{align*}
 W_1(\nu_K, \nu_{K+1})
 	&\leq \sum_{\mathtt{i}\in\Sigma_K} \sum_{\mathtt{j} \in \Sigma_1} p(\mathtt{ij}) | x(\mathtt{i}) - x'(\mathtt{ij}) | \\
  	&\le  \sum_{\mathtt{i}\in\Sigma_K} \sum_{\mathtt{j} \in \Sigma_1} \sup_{x,y\in\pi([\mathtt{ij}])}  p(\mathtt{ij}) |x -y|   \\
  	&\le  \sum_{\mathtt{i}\in\Sigma_K} \sum_{\mathtt{j} \in \Sigma_1}  p(\mathtt{ij}) N^{-K-1}   \\	
  	&\leq \sum_{\mathtt{i}\in\Sigma_K}p(\mathtt{i}) N^{-K}= N^{-K}.
\end{align*}
Thus, for $m>K$, we have
\[
  W_1(\nu_K,\nu_{m}) \leq W_1(\nu_K,\nu_{K+1})+\dots+W_1(\nu_{m-1},\nu_m) \leq N^{-K}+ N^{-(K +1)} + \dots + N^{-m} \leq N^{-(K-1)},
\]
and so $(\nu_n)_n$ is a Cauchy sequence with respect to $W_1$.
As argued above, $(\nu_n)_n$ converges weakly to some measure $\nu_\lambda$ in the space of probability measures on $[0,1]$.
As noted above, the argument for $\widetilde{\nu}_\lambda$ is identical.

We proceed to show that $\nu_\lambda$ and $\widetilde{\nu}_\lambda$ are non-atomic.
Fix $x\in  [0,1]$ and $n \ge 1$.
As remarked above, there are at most two $\mathtt{i}, \mathtt{j}\in\Sigma_{\infty}$ with $\pi(\mathtt{i})=\pi(\mathtt{j})=x$.
Note also that for all $k>n$, we have $\nu_k(\mathrm{int}(\pi([\mathtt{i}|_{n}])))=\nu_k(\pi([\mathtt{i}|_{n}]))=p(\mathtt{i}|_n)$.
By Theorem \ref{equivalent conditions for weak convergence} and Lemma \ref{thm:weightDecay},
\begin{align*}
  \nu_\lambda(\left\{ x \right\})
  &\leq \sum_{\substack{\mathtt{i}\in\Sigma_{\infty}:\\\pi(\mathtt{i})=x}}\nu_\lambda(\mathrm{int}(\pi([\mathtt{i}|_n])))
  \\
  &\leq \liminf_{k\to\infty} \sum_{\substack{\mathtt{i}\in\Sigma_{\infty}:\\\pi(\mathtt{i})=x}}\nu_k(\mathrm{int}(\pi([\mathtt{i}|_n])))
  \\
  &=\sum_{\substack{\mathtt{i}\in\Sigma_{\infty}:\\\pi(\mathtt{i})=x}}p(\mathtt{i}|_n)
  \leq 2C_N (N-1)^{-n}.
\end{align*}
The claim now follows for $\nu_\lambda$.
The case for $\widetilde{\nu}_\lambda$ is similar and left to the reader.

Let $\mathcal{O}\subset[0,1]$ be open. Since $\pi(\Sigma_{\infty})=[0,1]$, there exists an infinite word $\mathtt{i}\in\Sigma_{\infty}$ such that $\pi(\mathtt{i})\in \mathcal{O}$.
Since $\pi([\mathtt{i}|_n])$ is an interval of length $N^{-n}$, we must have $\pi([\mathtt{i}|_n])\subseteq \mathcal{O}$ for large enough $n$.
Hence, by Theorem~\ref{equivalent conditions for weak convergence},
\[
  \nu_\lambda(\mathcal{O}) \geq \nu_\lambda(\pi([\mathtt{i}|_n])) \geq \limsup_{k\to\infty} \nu_k(\pi([\mathtt{i}|_n])) = p(\mathtt{i}|_n) >0.
\]
Thus, $\nu_\lambda$
must have support $[0,1]$.

By considering open sets $\mathcal{O}$ that intersect $\pi(\mathcal{C})$, it can similarly shown
that the support of $\widetilde{\nu}_\lambda$ must contain $\pi(\mathcal{C})$.
Details are left to the reader.
To see that the support is exactly $\pi(\mathcal{C})$, consider the countable decomposition
\[
  \widetilde{\nu}_\lambda(\pi(\mathcal{C})^c) \leq \sum_{n=0}^{\infty}\sum_{\mathtt{i}\in\Sigma_n}\sum_{\mathtt{b}\in B}
  \widetilde{\nu}_\lambda(\pi([\mathtt{i}\mathtt{b}]))
  \leq \sum_{n=0}^{\infty}\sum_{\mathtt{i}\in\Sigma_n}\sum_{\mathtt{b}\in B}
  \widetilde{p}_{n+|\mathtt{b}|}(\mathtt{i}\mathtt{b})
  =0.
\]
Hence our claim follows.
\end{pf}


We now move from the abstract non-atomic measures $\lambda_{\pi([\mathtt{i}])}$ above to the specific measures $\widetilde{L}_{\pi([\mathtt{i}])}$ and $L_{\pi([\mathtt{i}])}$ that we require.   We begin by considering  the measures
\[
\widetilde{L}_{\pi([\mathtt{i}])} :=  \frac{\chi_{\pi([\mathtt{i}])}(x)}{\mathcal{L}(\pi([\mathtt{i}]))} \cdot d\mathcal{L}(x), \ \ \ \ \ \mathtt{i} \in \Sigma_*, \]
where $\cL$ denotes Lebesgue measure on $[0,1]$.  By Lemma~\ref{lemmasmallpnew}, we find that the sequence $(\nu_{n, \widetilde{L}})_{n=1}^\infty$ converges to a weak limit $\widetilde{\nu}_{\widetilde{L}}$ which is defined on $\pi(\mathcal{C})$, the Cantor-like set avoiding forbidden words from $B$.  Since we have now fixed $\widetilde{L}$, let us abbreviate  $\widetilde{\nu}_{\widetilde{L}}$ to $\widetilde{\nu}$, and $\widetilde{\nu}_{n,\widetilde{L}}$ to $\widetilde{\nu}_n$, $n\geq 1$.   Note that by Lemma~\ref{lemmasmallpnew}, $\widetilde{\nu}$ is a non-atomic probability measure.

Next, we set $L_{\pi([\mathtt{i}])} := F_{\mathtt{i}}{}_*\widetilde{\nu}$,
where $F_{\mathtt{i}}(x) = N^{-|\mathtt{i}|}x+\pi(\mathtt{i})$ is the linear mapping that maps $[0,1]$ into $\pi([\mathtt{i}])$.
A second application of Lemma~\ref{lemmasmallpnew} yields a non-atomic probability measure $\nu_L$  with support $[0,1]$ as the weak limit of the corresponding sequence $(\nu_{n, L})_n$.  Again, since $L$ is now fixed (since we fixed $\widetilde{L}$ and thus $\widetilde{\nu}$), we abbreviate $\nu_L$ to $\nu$, and $\nu_{n,L}$ to $\nu_n$, $n\geq 1$.  We then define $\mu_{0} := f_* \nu$ on $[0,\tfrac1N]$ as the push-forward measure of $\nu$ under $f(x)=x/N$, that is:
\[
  \mu_{0}(E) = \nu(N \cdot E)\text{ for all measurable }E.
\]
(As always, we should keep in mind that $\nu$ and thus that $\mu_0$ depend upon our original choice of $B$.)

Now $\mu_0$ is a non-atomic probability measure supported on $[0,\tfrac1N]$.
It will be convenient to think of $\mu_{0}$ as a measure defined on $\bbR$, taking the value $0$ for sets contained in $\bbR \setminus [0, \frac{1}{N}]$.
We construct $\mu$ ($=\mu_B$) satisfying the translational invariance in condition (i) of Theorem~\ref{thm:main}  by setting
\[
  \mu := \sum_{i=0}^{N-1}  g_i{}_*\mu_{0} = \sum_{i=0}^{N-1} g_i{}_*(f_*\nu)
  =\sum_{i=0}^{N-1}(g_i\circ f)_{*}\nu,
\]
where $g_i(x)=x+\frac{i}{N}$, $i=0,1,\cdots,N-1$.
Observe that $\mu$ is by construction a non-atomic Borel measure, supported on $[0,1]$.

\bigskip


Next we set up some additional notation.
For finite sequences $\mathtt{i}\in\Sigma_*$, we write
\begin{multline*}
  \langle\mathtt{i}\rangle
  := \big\{ \mathtt{j}\in\Sigma_\infty : i_k = j_k \text{ for all }k\leq |\mathtt{i}|\\
    \text{ and if  $\mathtt{j}=\mathtt{k}\mathtt{b}\mathtt{l}$ for some $\mathtt{b}\in B,
  \mathtt{k}\in\Sigma_*,\mathtt{l}\in\Sigma_\infty$ then }|\mathtt{k}\mathtt{b}|\leq |\mathtt{i}|\big\}.
\end{multline*}
Thus if $\mathtt{j} \in \langle\mathtt{i}\rangle$, then $\mathtt{j}$ admits $\mathtt{i}$ as a prefix, and any element $\mathtt{b} \in B$ which appears as a subword of $\mathtt{j}$ appears only as a subword of the prefix $\mathtt{i}$.
To conform to this new notation, we may also write $\langle.\rangle = \mathcal{C}$ for all codings not containing any subword $\mathtt{b}\in B$.

We write $\Lambda$ for the set of codings that contain only finitely many forbidden words from $B$.
That is,
\[
  \Lambda =\langle.\rangle\cup\Bigg(\bigcup_{k=0}^{\infty}\bigcup_{\mathtt{i}\in\Sigma_k}\bigcup_{\mathtt{b}\in B}
  \langle\mathtt{i}\mathtt{b}\rangle\Bigg).
\]
We note here that the union above is disjoint by our careful choice of $B$.

For $k\in \mathbb{N}_0:=\mathbb{N}\cup\{0\}$ and $\mathtt{b}\in B$, define 
\[E_{k,\mathtt{b}}=\{\mathtt{i}\mathtt{b}\mathtt{j}\in
\Sigma_{\infty}: \mathtt{i}\in\Sigma_k, \mathtt{j}\in\Sigma_\infty\}=\bigcup_{\mathtt{i}\in\Sigma_k}[\mathtt{i}\mathtt{b}].\]
Clearly the boundary of $\pi(E_{k,\mathtt{b}})$ is a finite set for each $k\in \bbN$ and
$\mathtt{b}\in B$ as $\pi([ \mathtt{i}\mathtt{b}])$ is an interval.

We claim that $\underset{n \to \infty}{\lim} \nu_n(\pi(E_{k, \mathtt{b}})) \le 2^{-(k + |\mathtt{b}|)}.$  To see this, let  $\mathtt{i} \in \Sigma_k$, and consider $[\mathtt{i}\mathtt{b}]$.
For $n \ge k+|\mathtt{b}|$, we obtain from the definition of $\nu_n$ that
\[
  \nu_n(\pi([\mathtt{i}\mathtt{b}]))
  = p(\mathtt{ib}) = \prod_{l=1}^{|\mathtt{i}| + |\mathtt{b}|} p_l (\mathtt{ib}|_l)
  = \left( \prod_{l=1}^k p_l (\mathtt{ib}|_l) \right) \, \left( \prod_{l=k+1}^{k + |\mathtt{b}|} p_l(\mathtt{ib}|_l) \right).
\]
However, from paragraph~\ref{paragraph2.5}, we see that
\begin{itemize}
  \item{}
    for $1 \le l \le k=|\mathtt{i}|$, $p_l(\mathtt{ib}|_l) = p_l(\mathtt{i}|_l)$;
  \item{}
    for $k + 1 \le l \le k + |\mathtt{b}| -1$, $p_l(\mathtt{ib}|_l) < 1$; and
  \item{}
    for $l = k+|\mathtt{b}|$, since $\mathtt{b} \in B$, $p_l(\mathtt{ib}|_l) = 2^{-(k+|\mathtt{b}|)}$.
\end{itemize}
Therefore, by the definition of $p$,
\[
\nu_n(\pi([\mathtt{i}\mathtt{b}])) \le \left( \prod_{l=1}^k p_l(\mathtt{i}|_l) \right) 2^{-{(k + |\mathtt{b}|)}} = p(\mathtt{i}) 2^{-(k+|\mathtt{b}|)}, \]
for $n\geq k+|\mathtt{b}|$ and thus $\nu(\pi(E_{\mathtt{i},k,\mathtt{b}})) \leq p(\mathtt{i}) 2^{-(k+|\mathtt{b}|)}$.
Consequently,
\[
  \nu(\pi(E_{k,\mathtt{b}})) \leq  \sum_{\mathtt{i} \in \Sigma_k} \nu(\pi(E_{\mathtt{i}, k, \mathtt{b}})) \le \left (\sum_{\mathtt{i} \in \Sigma_k} p(\mathtt{i}) \right) 2^{-(k+|\mathtt{b}| )}= 2^{-(k+|\mathtt{b}|)}.
\]
Next, note that $\sum_{k=0}^\infty \sum_{\mathtt{b}\in B}\nu(\pi(E_{k,\mathtt{b}}))< \infty$, and so by the Borel-Contelli Lemma,
\[
  \nu(\pi(\Lambda^c))= \nu\Bigg( \pi\Big(\Big\{\mathtt{l}\in \Sigma_\infty: \#\{(k,\mathtt{b}): \mathtt{l}\in E_{k,\mathtt{b}}\}=\infty\Big \}\Big)\Bigg)=0.
\]
Consequently, $\mu_{0} ([0,\frac{1}{N})\cap \pi(\Lambda^c))=\nu(\pi(\Lambda^c))=0$.


\begin{lem}\label{restriction}
  For $\mathtt{j}\in \Sigma_*$, the sequence $(\nu_{n}|_{\pi(\langle\mathtt{j}\rangle)})_n$ converges weakly to ${\nu}|_{\pi(\langle\mathtt{j}\rangle)}$.
\end{lem}

\begin{pf}
  Let $\mathtt{j}\in \Sigma_*$ and $K=|\mathtt{j}|$.
  The set
  \[
    \mathcal{O} = \bigcup_{n=0}^\infty \bigcup_{\mathtt{i}\in\Sigma_n}\bigcup_{\mathtt{b}\in B}\mathrm{int}\;\pi([\mathtt{j}\mathtt{i}\mathtt{b}])
  \]
  is a countable union of open sets and thus itself open.
  Hence $\mathcal{O}^c$ is closed in $[0,1]$, and we further have $\pi(\langle\mathtt{j}\rangle) \subseteq \mathcal{O}^c$.
  This implies that $\chi_{\pi(\langle\mathtt{j}\rangle)}\leq \chi_{\mathcal{O}^c\cap \pi([\mathtt{j}])}$.
  Note that
  \[\pi([ \mathtt{j} ])=\pi(\langle\mathtt{j}\rangle)\cup (\bigcup_{n=0}^\infty \bigcup_{\mathtt{i}\in\Sigma_n}\bigcup_{\mathtt{b}\in B}\pi([\mathtt{j}\mathtt{i}\mathtt{b}])).
  \]
  Hence
  \[(\mathcal{O}^c \cap \pi([ \mathtt{j} ]))\setminus \pi(\langle\mathtt{j}\rangle)
    =\mathcal{O}^c \cap (\bigcup_{n=0}^\infty \bigcup_{\mathtt{i}\in\Sigma_n}\bigcup_{\mathtt{b}\in B}\pi([\mathtt{j}\mathtt{i}\mathtt{b}]))=\bigcup_{n=0}^\infty \bigcup_{\mathtt{i}\in\Sigma_n}\bigcup_{\mathtt{b}\in B}(\pi([\mathtt{j}\mathtt{i}\mathtt{b}])\setminus \mathrm{int}\;\pi([\mathtt{j}\mathtt{i}\mathtt{b}])).
  \]
  Therefore, $(\mathcal{O}^c \cap \pi([ \mathtt{j} ]))\setminus \pi(\langle\mathtt{j}\rangle)$ is countable.
  Recalling that $\nu$ is non-atomic, we have
  \[
    \chi_{\pi(\langle\mathtt{j}\rangle)}(x) = \chi_{\mathcal{O}^c\cap  \pi([ \mathtt{j} ])}(x) \ \ \ \ \nu\text{-a.e.}\ x\in[0,1].
  \]

  Observe that we can partition $[ \mathtt{j}]$ as a disjoint union of $\langle\mathtt{j}\rangle$ and a countable union of intervals (that are not necessarily disjoint) as follows:
  \[
    [ \mathtt{j} ] = \langle\mathtt{j}\rangle\sqcup \bigcup_{n=0}^\infty \bigcup_{\mathtt{i}\in\Sigma_n}\bigcup_{\mathtt{b}\in B}[\mathtt{j}\mathtt{i}\mathtt{b}].
  \]
  Since $\pi([ \mathtt{j} ])$ is a closed interval and hence the boundary of $\pi ([\mathtt{j} ])$ consists of two points, by Theorem \ref{equivalent conditions for weak convergence}, $\nu(\pi([\mathtt{j}])) = \lim_{k \to \infty} \nu_k(\pi([\mathtt{j}]))$.
  Thus, to show that  $ \nu(\pi(\langle\mathtt{j}\rangle)) = \lim_{k \to \infty} \nu_k(\pi(\langle\mathtt{j}\rangle))$,
  it suffices to show that
  \[
    \nu(\pi(\bigcup_{n=0}^\infty \bigcup_{\mathtt{i}\in\Sigma_n}\bigcup_{\mathtt{b}\in B}[\mathtt{j}\mathtt{i}\mathtt{b}])) = \lim_{k \to \infty} \nu_k(\pi(\bigcup_{n=0}^\infty \bigcup_{\mathtt{i}\in\Sigma_n}\bigcup_{\mathtt{b}\in B}[\mathtt{j}\mathtt{i}\mathtt{b}])).
  \]
  Fix a word $\mathtt{i}\in \Sigma_n$ and a forbidden word $\mathtt{b}\in B$.
  We then have
  \[
    \nu_{m}(\pi([\mathtt{j}\mathtt{i}\mathtt{b}]))=\nu_{K+n+|\mathtt{b}|}(\pi([\mathtt{j}\mathtt{i}\mathtt{b}])),
  \]
  for $m\geq K+n+|\mathtt{b}|$.
  Assume now that $K+1 \leq m< K+n+|\mathtt{b}|$.
  Write $\sigma\mathtt{j}$ for the left-shift map $\mathtt{j}= (j_1 \, j_2 \cdots j_K) \mapsto (j_2 \, j_3\, \cdots \, j_K)$.
  We get
  \begin{align*}
    \nu_{m}(\pi([\mathtt{j}\mathtt{i}\mathtt{b}]))
    &=
    \int_{\pi([\mathtt{j}\mathtt{i}\mathtt{b}])}
    p(\mathtt{j}\mathtt{i}\mathtt{b}|_m)\, d L_{\pi([\mathtt{j}\mathtt{i}\mathtt{b}|_m])}
    \\
    &=
    p(\mathtt{j}\mathtt{i}\mathtt{b}|_m)\cdot
    \widetilde{\nu}(\pi([\sigma^m(\mathtt{j}\mathtt{i}\mathtt{b})]))
    \\
    &=
    p(\mathtt{j}\mathtt{i}\mathtt{b}|_{m-1})\cdot
    p_{m}(\mathtt{j}\mathtt{i}\mathtt{b}|_{m})\cdot
    \widetilde{\nu}(\pi([\sigma^m(\mathtt{j}\mathtt{i}\mathtt{b})]))
    \\
    &\geq
    \frac{1}{N-1}\cdot
    p(\mathtt{j}\mathtt{i}\mathtt{b}|_{m-1})\cdot
    \widetilde{\nu}(\pi([\sigma^m(\mathtt{j}\mathtt{i}\mathtt{b})]))
    \\
    &\geq
    p(\mathtt{j}\mathtt{i}\mathtt{b}|_{m-1})\cdot
    \widetilde{\nu}(\pi([\sigma^{m-1}(\mathtt{j}\mathtt{i}\mathtt{b})]))
    \\
    &=\nu_{m-1}(\pi([\mathtt{j}\mathtt{i}\mathtt{b}])).
  \end{align*}
  We may apply Theorem \ref{equivalent conditions for weak convergence} to conclude that
  \[
    \left( \nu_{m}\left( \pi\left(  \bigcup_{n=0}^\infty \bigcup_{\mathtt{i}\in\Sigma_n}\bigcup_{\mathtt{b}\in B}[\mathtt{j}\mathtt{i}\mathtt{b}]\right) \right) \right)_{m=1}^{\infty}
  \]
  \bigskip
  is an (eventually) increasing sequence which converges to
  \[
    \nu(\pi(\bigcup_{n=0}^\infty \bigcup_{\mathtt{i}\in\Sigma_n}\bigcup_{\mathtt{b}\in B}[\mathtt{j}\mathtt{i}\mathtt{b}])).
  \]
  Finally, we note that $\chi_{\mathcal{O}^c\cap {\pi([ \mathtt{j} ])}}$ is upper semi-continuous as
  $(\mathcal{O}^c \cap {\pi([ \mathtt{j} ])})$ is closed.
  Let $h$ be any bounded upper semi-continuous functions on $[0,1]$.
  Let $t\geq 0$ be s.t. $h(x)+t\geq 0$.
  Then $(h(x)+t)\cdot \chi_{\mathcal{O}^c\cap {\pi([ \mathtt{j} ])}}$ is upper semi-continuous.
  By Theorem \ref{equivalent conditions for weak convergence}, we have
  \begin{align*}
    \underset{n\to \infty}{\lim\sup}\int h \ d\nu_{n}|_{\pi(\langle\mathtt{j}\rangle)}
    &=\underset{n\to \infty}{\lim\sup}\int h \chi_{\pi(\langle\mathtt{j}\rangle)} \ d\nu_{n}\\
    &=\underset{n\to \infty}{\lim\sup}\left(\int(h+t)\chi_{\pi(\langle\mathtt{j}\rangle)}\ d\nu_{n} - t \nu_{n}(\pi(\langle\mathtt{j}\rangle))\right)\\
    &= \underset{n\to \infty}{\lim\sup}\int(h+t)\chi_{\mathcal{O}^c\cap {\pi([ \mathtt{j} ])}}\ d\nu_n-t\nu(\pi(\langle\mathtt{j}\rangle))\\
    &\leq \int (h+t)\chi_{\mathcal{O}^c\cap {\pi([ \mathtt{j} ])}}\ d \nu-t\nu(\pi(\langle\mathtt{j}\rangle))\\
    &=\int (h+t)\chi_{\pi(\langle\mathtt{j}\rangle)}\ d \nu -t\nu (\pi(\langle\mathtt{j}\rangle))\\
    &=\int h \ d \nu|_{\pi(\langle\mathtt{j}\rangle)}.
  \end{align*}
  By one last application of Theorem \ref{equivalent conditions for weak convergence}, we conclude that $(\nu_{n}|_{\pi(\langle\mathtt{j}\rangle)})_n$ converges weakly to $\nu|_{\pi(\langle\mathtt{j}\rangle)}$.
\end{pf}


\begin{lem}\label{thm:numucontinuous}
  The measures $\nu|_{[0,\frac{1}{N}]}$ and $\mu_0$ are mutually absolutely continuous; in other words,  $\nu|_{[0,\frac{1}{N}]}\asymp \mu_0$.
\end{lem}

\begin{pf}
  Note that
  \begin{equation}\label{eq:decomposition}
    \big[0,\tfrac{1}{N}\big)
    =\pi(\langle0\rangle)\cup\Bigg(\bigcup_{k=0}^\infty \bigcup_{\mathtt{j}\in\Sigma_k}\bigcup_{\mathtt{b}\in B}
    \pi(\langle0\mathtt{j}\mathtt{b}\rangle)\Bigg)\cup
    \big([0,\tfrac{1}{N})\cap \pi(\Lambda^c)\big)
  \end{equation}
  and observe that
  \[
  \mu_0([0,\tfrac{1}{N})\cap \pi(\Lambda^c))=\nu([0,\tfrac{1}{N})\cap \pi(\Lambda^c))=0. \]
  Further, note that the unions in \eqref{eq:decomposition} are disjoint symbolically, i.e.\ all distinct $\langle0\rangle, \langle0\mathtt{j}\mathtt{b}\rangle$ have empty intersection.
  While this is no longer true for the images under $\pi$, the images of any distinct cylinders still only have countably many intersections.
  Hence, it suffices to show that $\mu_0|_{\pi(\langle\mathtt{j}\rangle)}\asymp \nu|_{\pi(\langle\mathtt{j}\rangle)}$, $\forall \mathtt{j}\in \Sigma_*$ with $j_1=0$.

  For $n\in \bbN$, define $\mu_n=f_*\nu_n$, where $f(x)=\frac{x}{N}$, $x\in [0,1]$.
  Clearly $(\mu_n)_n$ converges weakly to $\mu_0$.
  Fix $\mathtt{j}\in \Sigma_*$ with $j_1=0$.
  Set $K:=|\mathtt{j}|$.
  Consider arbitrary $M\in\bbN$ and  $\mathtt{i}\in\Sigma_M$ such that if $\mathtt{b} \in B$ is not
  a subword of $\mathtt{j}$ then it is also
  not a subword of $\mathtt{ji}$.
  Let $A_{\mathtt{ji}}$ be such that
  $\mathrm{int}\,\pi([\mathtt{ji}])\subseteq A_{\mathtt{ji}}\subseteq\pi([\mathtt{ji}])$.
  That is, $A_{\mathtt{ji}}$ is either the closed interval $\pi([\mathtt{ji}])$, its interior, or one of
  the two half-open subintervals.

  For $n\geq |\mathtt{ji}|=K+M$,
  \[
    \nu_{n+1}|_{\pi(\langle\mathtt{j}\rangle)}(A_{\mathtt{ji}})
    =\prod_{k=1}^{K+M}p_k(\mathtt{j}\mathtt{i}|_{k})
  \]
  and, as $j_1=0$,
  \begin{align*}
    \mu_n|_{\pi(\langle\mathtt{j}\rangle)}(A_{\mathtt{ji}})
    &=\nu_n|_{\pi\left( \langle\mathtt{j}\rangle \right)}(N\cdot A_{\mathtt{ji}})
    =
    \prod_{k=1}^{K+M-1}p_k(\sigma(\mathtt{j}\mathtt{i})|_{k})
  \end{align*}

  Note that no forbidden word gets completed in $\mathtt{ji}|_k$ for $k>K$. Therefore
  \[
    p_k(\sigma(\mathtt{ji})|_k) =
    \begin{cases}
      \frac{1-2^{-k}}{N-1} &\text{if }\#B(\sigma(\mathtt{ji})|_{k-1})=1
      \\
      \frac{1}{N} & \text{if }\#B(\sigma(\mathtt{ji})|_{k-1})=0
    \end{cases}
  \]
  and
  \[
    p_{k+1}(\mathtt{ji}|_{k+1}) =
    \begin{cases}
      \frac{1-2^{-{k+1}}}{N-1} &\text{if }\#B(\mathtt{ji}|_{k})=1
      \\
      \frac{1}{N} & \text{if }\#B(\mathtt{ji}|_{k})=0
    \end{cases}
  \]
  for $k\geq K$.
  Further observe that $\#B(\mathtt{ji}|_{k}) = \#B(\sigma(\mathtt{ji})|_{k-1})$.
  Hence
  \begin{align*}
    \mu_{n}|_{\pi(\langle\mathtt{j}\rangle)}(A_{\mathtt{ji}})
    &=
    \frac{\prod_{k=1}^{K+M-1}p_k(\sigma(\mathtt{ji})|_k)}{\prod_{k=1}^{K+M}p_{k}(\mathtt{j}\mathtt{i}|_k)}
    \cdot\nu_{n+1}|_{\pi(\langle\mathtt{j}\rangle)}(A_{\mathtt{ji}})
    \\[0.5em]
    &=
    \frac{\prod_{k=1}^{K-1}p_k(\sigma(\mathtt{j})|_k)}{\prod_{k=1}^{K}p_{k}(\mathtt{j}|_k)}
    \cdot
    \prod_{k=K}^{K+M-1}
    \left(
      \frac{p_k(\sigma(\mathtt{ji})|_k)}{p_{k+1}(\mathtt{j}\mathtt{i}|_{k+1})}
    \right)
    \cdot\nu_{n+1}|_{\pi(\langle\mathtt{j}\rangle)}(A_{\mathtt{ji}})
    \\[0.5em]
    &=
    \underbrace{
      \vphantom{\prod_{\substack{K\leq k \leq K+ M-1\\ \#B(\mathtt{ji}|_{k})=1}}}
      \frac{\prod_{k=1}^{K-1}p_k(\sigma(\mathtt{j})|_k)}{\prod_{k=1}^{K}p_{k}(\mathtt{j}|_k)}
    }_{R_{\mathtt{j}}}
    \cdot
    \underbrace{
    \prod_{\substack{K\leq k \leq K+M-1\\ \#B(\mathtt{ji}|_{k})=1}}
    \left(
      \frac{1-2^{-k}}{1-2^{-(k+1)}}
    \right)
  }_{r_{\mathtt{ji}}}
    \cdot\nu_{n+1}|_{\pi(\langle\mathtt{j}\rangle)}(A_{\mathtt{ji}})
  \end{align*}
  Note that $R_{\mathtt{j}}$ depends only on $\mathtt{j}$ and that
  $r_{\mathtt{ji}}$ is bounded independently of $\mathtt{i}$:
  \[
    \frac12
    =
    \prod_{k=1}^{\infty}\frac{1-2^{-k}}{1-2^{-(k+1)}}
    \leq
    \prod_{\substack{K\leq k \leq K+M-1\\ \#B(\mathtt{ji}|_{k})=1}}
    \left(
      \frac{1-2^{-k}}{1-2^{-(k+1)}}
    \right)
    =r_{\mathtt{ji}}
    \leq
    1.
  \]
  By Lemma \ref{restriction}, $(\nu_{n}|_{\pi(\langle\mathtt{j}\rangle)})_n$ converges weakly
  to $\nu|_{\pi(\langle\mathtt{j}\rangle)}$, and similarly
  $(\mu_{n}|_{\pi(\langle\mathtt{j}\rangle)})_n$ converges weakly to $\mu_0|_{\pi(\langle\mathtt{j}\rangle)}$.
  Since $\mathrm{int}\,\pi([\mathtt{ji}]) \subseteq A_{\mathtt{ji}} \subseteq\mathrm{int}\,\pi([\mathtt{ji}])$ we may apply Theorem \ref{equivalent conditions for weak convergence} to conclude that
  $\nu|_{\pi(\langle\mathtt{j}\rangle)}$ and $\mu_0|_{\pi(\langle\mathtt{j}\rangle)}$
  are constant on all $A_{\mathtt{ji}}$ for fixed $\mathtt{j}$ and $\mathtt{i}$, respectively.
  We obtain
  \[
    \mu_0|_{\pi(\langle\mathtt{j}\rangle)}(\pi([\mathtt{ji}]))
    =
    \mu_0|_{\pi(\langle\mathtt{j}\rangle)}(A_{\mathtt{ji}})
    = R_{\mathtt{j}}\cdot r_{\mathtt{ji}}
    \cdot\nu|_{\pi(\langle\mathtt{j}\rangle)}(A_{\mathtt{ji}})
    =
     R_{\mathtt{j}}\cdot r_{\mathtt{ji}}\cdot
    \nu|_{\pi(\langle\mathtt{j}\rangle)}(\pi([\mathtt{ji}])).
  \]
  Write
  \[
    \mathcal{M}_{\mathtt{j}} = \{\mathtt{i}\in \Sigma_{*}
      \,:\, \mathtt{b} \in B  \text { not a subword of } \mathtt{j} \text{ implies  that }
    \mathtt{b} \text{ is not a subword of } \mathtt{ji} \}
  \]
  and
  \[
    \mathcal{R}_{\mathtt{j}} = \left\{ \bigcup_{i=1}^n A_i \,:\, n\in\bbN\text{ and } A_i\in \mathcal{R}'_{\mathtt{j}} \right\}
  \]
  which are the finite unions of elements in
  \[
    \mathcal{R}'_{\mathtt{j}} = \left\{ A_{\mathtt{ji}}\,:\, \mathtt{i}\in\mathcal{M}_{\mathtt{j}} \text{ and }
    \mathrm{int}\,\pi([\mathtt{ji}])\subseteq A_{\mathtt{ji}}\subseteq \pi([\mathtt{ji}])  \right\}.
  \]
  Note that $\mathcal{R}'_{\mathtt{j}}$, and so $\mathcal{R}_{\mathtt{j}}$, generates the Borel $\sigma$-algebra
  $\mathfrak{Bor}([0,\tfrac1N))|_{\pi(\langle\mathtt{j}\rangle)}$
  and that $\mathcal{R}_{\mathtt{j}}$ is stable under intersections (in fact, it is a ring of sets).
  By Carath\'eodory's extension theorem, see \cite[Theorem 5.4]{BauerMeasureIntegration},
  the finite measures $\nu|_{\pi(\langle\mathtt{j}\rangle)}$ and $\mu_0|_{\pi(\langle\mathtt{j}\rangle)}$ are uniquely determined by their values of elements in $\mathcal{R}_{\mathtt{j}}$.
  In particular, for all $E\in\mathfrak{Bor}([0,\tfrac1N))$ they can be expressed as the outer measures
  \begin{equation}\label{eq:nudecomp}
    \nu|_{\pi(\langle\mathtt{j}\rangle)}(E)
    =\inf\left\{ \sum_{i=1}^{\infty} \nu|_{\pi(\langle\mathtt{j}\rangle)}(E_i) \,:\,E_i\in\mathcal{R}_{\mathtt{j}}\text{ and }E\subseteq\bigcup_{i=1}^{\infty}E_i \right\}
  \end{equation}
  and
  \begin{equation}\label{eq:mudecomp}
    \mu_0|_{\pi(\langle\mathtt{j}\rangle)}(E)
    =\inf\left\{ \sum_{i=1}^{\infty} \mu_0|_{\pi(\langle\mathtt{j}\rangle)}(E_i) \,:\,E_i\in\mathcal{R}_{\mathtt{j}}\text{ and }E\subseteq\bigcup_{i=1}^{\infty}E_i \right\}.
  \end{equation}
  Without loss of generality we may replace $\mathcal{R}_{\mathtt{j}}$ by $\mathcal{R}'_{\mathtt{j}}$ in \eqref{eq:nudecomp} and \eqref{eq:mudecomp} since we are taking infimums over countable unions.
  We may further assume that the $E_i\in\mathcal{R}'_{\mathtt{j}}$
  are of the form $\pi([\mathtt{ji}])$ for some $\mathtt{i}\in\mathcal{M}$ as
  $\mathrm{int}\, \pi([\mathtt{ji}])\subseteq E_i\subseteq \pi([\mathtt{ji}])$ and the measures $\mu_0|_{\pi(\langle\mathtt{j}\rangle)}$ and $\nu|_{\pi(\langle\mathtt{j}\rangle)}$ of those sets agree, respectively.
  That is, we may write
  \begin{equation*}
    \nu|_{\pi(\langle\mathtt{j}\rangle)}(E)
    =\inf\left\{ \sum_{k=1}^{\infty} \nu|_{\pi(\langle\mathtt{j}\rangle)}(\pi([\mathtt{ji}^{(k)}]))
    \,:\,\mathtt{i}^{(k)}\in\mathcal{M}_{\mathtt{j}}\text{ and }E\subseteq\bigcup_{k=1}^{\infty}\pi([\mathtt{ji}^{(k)}]) \right\}
  \end{equation*}
  with an analogous expression holding for $\mu_{0}|_{\pi(\langle\mathtt{j}\rangle)}$.

  Thus, for any $E\in\mathfrak{Bor}([0,\tfrac1N))$,
  \begin{align*}
    &\mu_{0}|_{\pi(\langle\mathtt{j}\rangle)} (E)
    =
    \inf\left\{ \sum_{k=1}^{\infty} \mu_0|_{\pi(\langle\mathtt{j}\rangle)}(\pi([\mathtt{ji}^{(k)}]))
      \,:\,\mathtt{i}^{(k)}\in\mathcal{M}_{\mathtt{j}}\text{ and }
    E\subseteq\bigcup_{k=1}^{\infty}\pi([\mathtt{ji}^{(k)}]) \right\}
    \\
    &=
    \inf\left\{ \sum_{k=1}^{\infty} R_{\mathtt{j}}\cdot r_{\mathtt{ji}}\cdot
      \nu|_{\pi(\langle\mathtt{j}\rangle)}(\pi([\mathtt{ji}^{(k)}]))
      \,:\,\mathtt{i}^{(k)}\in\mathcal{M}_{\mathtt{j}}\text{ and }
    E\subseteq\bigcup_{k=1}^{\infty}\pi([\mathtt{ji}^{(k)}]) \right\}
    \\
    &\leq
    R_{\mathtt{j}}\sup_{\mathtt{i}\in\mathcal{M}_{\mathtt{j}}}r_{\mathtt{ji}}
    \cdot \inf\left\{ \sum_{k=1}^{\infty}
      \nu|_{\pi(\langle\mathtt{j}\rangle)}(\pi([\mathtt{ji}^{(k)}]))
      \,:\,\mathtt{i}^{(k)}\in\mathcal{M}_{\mathtt{j}}\text{ and }
    E\subseteq\bigcup_{k=1}^{\infty}\pi([\mathtt{ji}^{(k)}]) \right\}
    \\[.5em]
    &\leq R_{\mathtt{j}}\cdot
    \nu|_{\pi(\langle\mathtt{j}\rangle)} (E).
  \end{align*}
  Similarly,
  for any $E\in\mathfrak{Bor}([0,\tfrac1N))$,
  \begin{align*}
    &\mu_{0}|_{\pi(\langle\mathtt{j}\rangle)} (E)
    =
    \inf\left\{ \sum_{k=1}^{\infty} \mu_0|_{\pi(\langle\mathtt{j}\rangle)}(\pi([\mathtt{ji}^{(k)}]))
      \,:\,\mathtt{i}^{(k)}\in\mathcal{M}_{\mathtt{j}}\text{ and }
    E\subseteq\bigcup_{k=1}^{\infty}\pi([\mathtt{ji}^{(k)}]) \right\}
    \\[0.5em]
    &\geq
    R_{\mathtt{j}}\cdot\inf_{\mathtt{i}\in\mathcal{M}_{\mathtt{j}}}r_{\mathtt{ji}}
    \cdot
    \nu|_{\pi(\langle\mathtt{j}\rangle)} (E)
    \geq
    \tfrac12 R_{\mathtt{j}}\cdot
    \nu|_{\pi(\langle\mathtt{j}\rangle)} (E).
  \end{align*}
  Thus $\nu|_{\pi(\langle\mathtt{j}\rangle)}\asymp\mu_0|_{\pi(\langle\mathtt{j}\rangle)}$
  as required.

\end{pf}

\begin{lem}\label{thm:translation}
  We have $\nu|_{\pi([1])}\asymp g_1{}_*\nu|_{\pi([0])}$ and $\nu|_{\pi([i])} = g_i{}_*\nu|_{\pi([0])}$ for all $1<i\leq N-1$.
\end{lem}

\begin{pf}
  For $k \ge 0$ an integer, let us write $0^k$ to denote the word $(0 \, 0 \, 0 \, \cdots \, 0) \in \Sigma_k$.  The measure $\nu|_{\pi([1])}$ decomposes into countably many cylinder measures
  \begin{align*}
    \nu|_{\pi([1])}
    &= \sum_{k=1}^\infty \nu|_{\pi([10^k2])}\;+
    \sum_{i \neq 2 }
    \sum_{k=0}^{\infty}
    \;\nu|_{\pi([10^ki])}\;+ \nu|_{\pi([12])}
    \\
    &=\sum_{k=1}^\infty \frac{\prod_{j=1}^{k+2}p_{j}(00^k2|_j)}{\prod_{j=1}^{k+2}p_{j}(10^k2|_j)}g_1{}_*\nu|_{\pi([00^k2])}\;+
    \sum_{i\neq2}\sum_{k=0}^{\infty}\frac{\prod_{j=1}^{k+2}p_{j}(00^ki|_j)}{\prod_{j=1}^{k+2}p_{j}(10^ki|_j)}g_1{}_*\nu|_{\pi([00^ki])}\;
    \\
    &  \ \ \ \ \ \ \ \ \ \ \ \ \ \ \ \ \ \ \ \ \ \ \ \ \
    + \frac{\prod_{j=1}^{2}p_{j}(02|_j)}{\prod_{j=1}^{2}p_{j}(12|_j)}g_1{}_*\nu|_{\pi([02])}
    \\
    &=\sum_{k=1}^\infty \frac{p_{k+2}(00^k2)}{p_{k+2}(10^k2)}g_1{}_*\nu|_{\pi([00^k2])}\;+
    \sum_{i\neq 2}\sum_{k=0}^{\infty}g_1{}_*\nu|_{\pi([00^ki])}\;+ g_1{}_*\nu|_{\pi([02])}.
  \end{align*}
  Noting that $g_1{}_*\nu|_{\pi([0])}$ decomposes as
  \[
    g_1{}_*\nu|_{\pi([0])} =\sum_{k=1}^\infty g_1{}_*\nu|_{\pi([00^k2])}\;+
    \sum_{l \ne 2}\sum_{k=0}^{\infty}g_1{}_*\nu|_{\pi([00^kl])}\;+ g_1{}_*\nu|_{\pi([02])}
  \]
  gives $\nu|_{\pi([1])}\asymp g_1{}_*\nu|_{\pi([0])}$.

  The second conclusion follows similarly, with equality achieved as
  $\frac{\prod_{j=1}^{k+2}p_{j}(00^k2|_j)}{\prod_{j=1}^{k+2}p_{j}(i0^k2|_j)}=1$ for $i\neq 1$.
\end{pf}




\bigskip

\noindent \textbf{Proof}~(\textbf{of Theorem~\ref{thm:main}}).
Condition (ii) of Theorem~\ref{thm:main} holds, since
\begin{align*}
  \mu_0
  \;=\;
  f_*\nu
  \;&=\;
  \sum_{i=0}^{N-1} f_*\nu|_{\pi([i])}
  \\
  &\asymp
  \sum_{i=0}^{N-1} (g_i\circ f)_*\nu|_{\pi([0])}
  &&\text{(Lemma \ref{thm:translation})}
  \\
  &\asymp
  \sum_{i=0}^{N-1} (g_i\circ f)_*\mu_0
  &&\text{(Lemma \ref{thm:numucontinuous})}
  \\
  &=
  f_* \mu .
\end{align*}
Finally, to complete the proof of Theorem~\ref{thm:main}, we need only show that condition (iii) of that Theorem holds.

Let us write the power set of $B_0$ as $\cP(B_0) = \{ B_\alpha: \alpha \in \Omega\}$, and note that the cardinality of $\Omega$ is $2^{\aleph_0}$.
For each $\alpha \in \Omega$ we have constructed a measure $\mu_\alpha$ with support equal to $[0,1]$ (arising from a measure $\nu_\alpha$ with support $[0, 1]$) satisfying conditions (i) and (ii) of Theorem~\ref{thm:main}.
There remains to show that if $\alpha \ne \beta \in \Omega$, then the associated measures $\mu_\alpha$ and $\mu_\beta$ are mutually singular.
Thus (iii) will be satisfied, completing the proof of the theorem.
(We remark in passing that the measure $\mu_\theta$ associated to the empty set $B_\theta := \varnothing$ is just $N$ times Lebesgue measure.)

To this end, let us fix $\alpha \ne \beta \in \Omega$.
Without loss of generality, we may assume that $B_\alpha$ is not contained in $B_\beta$.
Let
$\mathtt{b}\in B_\alpha\setminus B_\beta$ and write $\Lambda_{\mathtt{b}} = \left\{ \mathtt{i} \in\Sigma_\infty : \mathtt{b}\text{ occurs finitely many times in }\mathtt{i}\right\}$.
First, we show that $\nu_{\alpha}(\pi(\Lambda_{\mathtt{b}})) = 1$.
Note that $\Lambda_{\mathtt{b}}^c \subseteq \bigcup_{n=k}^{\infty}\bigcup_{\mathtt{i}\in\Sigma_n} [ \mathtt{i}\mathtt{b}]$ for all $k\in\bbN$.
But then
\[
  \nu_{\alpha}(\pi(\Lambda_{\mathtt{b}}^c))
  \leq \sum_{n=k}^\infty \sum_{\mathtt{i}\in\Sigma_{n}}\nu_{\alpha}(\pi([ \mathtt{i} \mathtt{b}]))
  \leq \sum_{n=k}^{\infty}2^{-n}\sum_{\mathtt{i}\in\Sigma_n}\nu_{\alpha}(\pi ([\mathtt{i}]))=2^{-k+1}
\]
and as $k$ was arbitrary, $\nu_{\alpha}(\pi(\Lambda_{\mathtt{b}}^c))=0$ and
$\nu_{\alpha}(\pi(\Lambda_{\mathtt{b}}))=1$.

Next, we show that $\nu_{\beta}(\pi(\Lambda_{\mathtt{b}}))=0$, proving that $\nu_\alpha$ and $\nu_\beta$ are not mutually absolutely continuous.
For this we introduce the sets:
\[
  \Lambda_{\mathtt{b},n} = \left\{ \mathtt{i}\in\Sigma_\infty : \text{there are no forbidden words $\mathtt{b}$ in $\mathtt{i}$ after digit $n$} \right\}
\]
and
\[
  \Lambda_{\mathtt{b},n,k} =
  \left\{ \mathtt{i}\mathtt{j}^1\mathtt{j}^2\cdots\mathtt{j}^k \mathtt{l}: \mathtt{i}\in\Sigma_n, \mathtt{j}^m\in(\Sigma_{|\mathtt{b}|}\setminus\left\{
  \mathtt{b} \right\}) \text{ for all }1\leq m\leq k, \mathtt{l}\in \Sigma_{\infty}\right\}.
\]
Clearly, $\Lambda_{\mathtt{b}} = \bigcup_{n\in\bbN_0} \Lambda_{\mathtt{b},n}$ and
$\Lambda_{\mathtt{b},n}\subseteq \Lambda_{\mathtt{b},n,k}$ for all $k\in\bbN$.
Now,
\begin{align*}
  \nu_{\beta}(\pi(\Lambda_{\mathtt{b},n,k}))
  &\leq \sum_{\mathtt{i}\in\Sigma_{n}}\sum_{\substack{\mathtt{j}^m\in(\Sigma_{|\mathtt{b}|}\setminus\left\{
	  \mathtt{b}
  \right\})\\1\leq m \leq k}}
  \prod_{i=1}^{n+|\mathtt{b}|k}p_{i}(\mathtt{i}\mathtt{j}^1\ldots\mathtt{j}^k|_i)
  \\
  &= \sum_{\mathtt{i}\in\Sigma_n} \left( \prod_{i=1}^n p_i(\mathtt{i}|_i)\right)
  \cdot\prod_{m=1}^k \left( 1-p_{n+m|\mathtt{b}|}(\mathtt{i}\mathtt{j}^1\ldots\mathtt{j}^k|_{n+m|\mathtt{b}|}) \right)
  \\
  &\leq \prod_{m=1}^k \left( 1-\frac{1-2^{-(n+m|\mathtt{b}|)}}{N} \right) \leq \left(
  1-\frac{1}{2N} \right)^{k}.
\end{align*}
Thus $\nu_{\beta}(\pi (\Lambda_{\mathtt{b},n}))\leq (1-1/(2N))^k$ for all $k$ and hence
$\nu_{\beta}(\pi(\Lambda_{\mathtt{b},n}))=0$.
Since
$\Lambda_{\mathtt{b}}=\bigcup_{n\in\bbN_0}\Lambda_{\mathtt{b},n}$
is a countable union, we necessarily have $\nu_{\beta}(\pi(\Lambda_{\mathtt{b}}))=0$, proving that the two measures $\nu_\alpha$ and $\nu_\beta$ are in fact  singular to one another.  Correspondingly, $\mu_\alpha$ and $\mu_\beta$ are also mutually singular.

By rescaling $\mu_\alpha$ by $\frac{1}{N}$ for all $\alpha \in \Omega$, we can ensure that each $\mu_\alpha$ is a probability measure.

\hfill \eop


\begin{thm} \label{unitary:main}
  Let $n \ge 3$.    There exists a continuum $\{ U_\alpha: \alpha \in \Omega\}$ of $\jnu$-stable, cyclic, non-atomic unitary operators on $\hilb$, no two of which are unitarily equivalent.
\end{thm}

\begin{pf}
  Set $N := n$, and let $\{ \mu_\alpha: \alpha \in \Omega\}$ be the continuum of non-atomic measures on $[0,1]$ satisfying the conditions of Theorem~\ref{thm:main}.  For each $\alpha \in \Omega$, let $M_\alpha$ denote the multiplication operator on $L^2([0,1], \mu_\alpha)$ defined by $(M_\alpha (f))(x) = x f(x)$ a.e.-$\mu_\alpha$ on $[0,1]$.

  Let $g: [0,1] \to \bbT$ be the continuous function $g(t) = e^{2 \pi i t}$, and let $U_\alpha := g(M_\alpha)$, $\alpha \in \Omega$.  Then $U_\alpha$ is a cyclic unitary on $(\bbT, \varrho_\alpha)$, where $\varrho_\alpha = g_* \mu_\alpha$ was specifically designed to satisfy the conditions of Theorem~\ref{measure}.   In particular, each $U_\alpha$ is $\jnu$-stable, and the fact that $\mu_\alpha$ and $\mu_\beta$ are not mutually absolutely continuous if $\alpha \ne \beta \in \Omega$ is easily seen to imply that $\varrho_\alpha$ is not mutually absolutely continuous with $\varrho_\beta$ unless $\alpha = \beta$.

\end{pf}



  The situation surrounding the $\jna$-stability of normal operators is ``infinitely" simpler.

  \begin{prop} \label{prop3.17}
    Let $N \in \bofh$ be a normal operator.   Then the following are equivalent:
    \begin{enumerate}
      \item[(i)]
	$N$ is approximately $\jna$-stable for all $n \ge 2$;
      \item[(ii)]
	$N$ is approximately $\jna$-stable for some $n \ge 2$;
      \item[(iii)]
	$\sigma(N) = \bbT$.
    \end{enumerate}
  \end{prop}

  \begin{pf}
    That (i) implies (ii) is trivial. Suppose that (ii) holds. Note that by paragraph~\ref{sec3.01}, given $T \in \bofh$, the operator $J_n(T)$ is normal if and only if $T$ is unitary.   This, combined with Proposition~\ref{prop2.06}(c), implies that we must have $\sigma(N) = \bbT$.

    To see that (iii) implies (i), observe that if $\sigma(N) = \bbT$, then $N \simeq_a D^{(\infty)}$, where $D = \ttt{diag} (d_n)_n$ and $\{ d_n: n \ge 1\}$ consists of all $m^{th}$ roots of unity, $m \ge 1$.  By Proposition~\ref{prop3.08}, $D^{(\infty)} \simeq J_n (D^{(\infty)})$, $n \ge 2$, implying that $D^{(\infty)}$ (and therefore $N$) is approximately $J_n$-stable  for all $n \ge 2$.
  \end{pf}


  We next consider $\textsc{j}_n$-stability properties for isometries and for weighted shifts.  Isometries acting on a Hilbert space admit a very well-defined structure, thanks to the Wold Decomposition Theorem.   Indeed, by the Wold decomposition theorem \cite[Theorem V.2.1]{Davidson1996}, there is a cardinal number $\alpha$ and a unitary $U$ such that $W=S^{(\alpha)}\oplus U$, where $S$ is the unilateral shift operator.

  \bigskip

  \begin{prop} \label{prop3.18}
    Let $W\in \cB(\hilb)$ be an isometry and let $W = S^{(\alpha)} \oplus U$ denote its Wold decomposition.   Let $n\geq 2$.
    Then $W$ is $\jnu$-stable if and only if $U$ $\jnu$-stable.
  \end{prop}

  \begin{pf}
    Recall that $S$ is $\ttt{j}_n^{[\ttt{u}]}$-stable.   Indeed, if $\{ e_k\}_{k=1}^\infty$ is an orthonormal basis for $\hilb$ relative to which $S e_k = e_{k+1}$ for all $k \ge 1$, then writing $\hilb_m := \ol{\mathrm{span}} \{ e_{m +k n} \}_{k=0}^\infty$, $1 \le m \le n$ we that the operator matrix for $S$ relative to $\hilb = \oplus_{m=1}^n \hilb_m$ is nothing more than $J_n(S)$.  Since $S$ is $\ttt{j}_n^{[\ttt{u}]}$-stable, by Proposition~\ref{prop2.10}, $W$ is $\jnu$-stable if $U$ is $\jnu$-stable.

    On the other hand,  if $W$ is $\jnu$-stable,  it follows that
    \[ S^{(\alpha)}\oplus U=W\simeq J_n(W)=J_n(S^{(\alpha)}\oplus U)\simeq J_n(S^{(\alpha)})\oplus J_n(U)\simeq S^{(\alpha)}\oplus J_n(U).\]
    By~\cite[Theorem 2.1]{ChenHerreroWu1992}, $U\simeq J_n(U)$.
  \end{pf}


  Next we  examine which weighted shifts are either $\jnu$-stable or $\jna$-stable.  A standard result \cite[Corollary 1]{Shi74} states that a weighted shift with weights $(w_n)_{n \in \Omega}$ (where $\Omega$ is one of $\bbN$, $-\bbN$ or $\bbZ$) is unitarily equivalent to a shift with weights $(|w_n|)_n$, $n \in \Omega$.   As such, it suffices to consider only those weighted shifts whose weight sequences consist of non-negative real numbers.


  \begin{lem}\label{lem3.19}
    Let $n \geq 2$, and $W \in \cB(\hilb)$ be a bilateral weighted shift with non-negative weights $(w_k)_{k \in \bbZ}$;  that is,  $ W e_ k =w_k e_{k+1}$, where
    $\{e_k\}_{k\in \mathbb{Z}}$ is an orthonormal basis of $\hilb$.
    Under a suitable choice of orthonormal basis of $\hilb^{(n)}$,  $J_n(W)$ is a
    bilateral weighted shift with positive
    weights $\{\cdots, \mathbbm{1}_{n-1}, w_{i},\mathbbm{1}_{n-1}, w_{i+1}, \mathbbm{1}_{n-1}, w_{i+2},\cdots\}$, where $\mathbbm{1}_{n-1}:=\overbrace{1,\cdots,1}^{n-1}$.
  \end{lem}

  \begin{pf}
    For $x=\sum_{i\in \bbZ} a_ie_i\in \hilb$, define $U x=x_n\oplus \cdots \oplus x_1$,
    where $x_l=\sum_{i\in \bbZ} a_{ni+(l-1)}e_i$, $1\leq l\leq n$.
    Clearly $U$ is a unitary which maps $\mathcal{H}$ onto $\mathcal{H}^{(n)}$.
    Thus, $U^*=U^{-1}$.
    Then
    \begin{align*}
      U^*J_n(W) Ux
      &= U^*J_n(W)[x_n\oplus \cdots \oplus x_1]  \\
      &=U^*[(x_{n-1}\oplus \cdots \oplus x_1)\oplus (W x_n)] \\
      &= [(\sum_{i\in \bbZ} a_{ni+(n-2)}e_{ni+(n-1)})+\cdots+(\sum_{i\in \bbZ} a_{ni+1}e_{ni+2})]+(\sum_{i\in \bbZ}w_ia_{ni+(n-1)}e_{ni+n}).
    \end{align*}
    Hence, it is easy to see that $U^*J_n(W)U$ is a bilateral weighted shift with weight sequence $\{\cdots, \mathbbm{1}_{n-1}, w_{i},\mathbbm{1}_{n-1}, w_{i+1}, \mathbbm{1}_{n-1}, w_{i+2},\cdots\}$.
  \end{pf}

  \bigskip


  Let $W\in \cB(\hilb)$ be a bilateral weighted shift with positive weights $\{w_i\}_{i\in \mathbb{Z}}$.
  Set $\Omega_{W}:=\{i \in \bbZ: w_i\neq 1\}$, and $d(W)=\min\{|i-j|: i\neq j, i,j\in \Omega_W\}$ if $|\Omega_{W}|\geq 2$.
  By Lemma~\ref{lem3.19}, it is easy to check that $d(J_n(W))=n\, \cdot d(W)$.
  This observation will be used in the  proof of Proposition~\ref{prop3.21}.

  For each $k\geq 1$,
  the $k$-spectrum of $W$ (denoted by $\sum_k W$) is defined to be the closure (in the usual
  topology on $\mathbb{R}^k$) of the set
  \[\{(w_{i+1}, w_{i+2}, \cdots, w_{i+k}): i\in \mathbb{Z}\}.\]


  The following Lemma is a result of O'Donovan \cite{O'Do75} (also see Proposition 2.2.14 of \cite{Mar1991}) under the assumption that $A,B$
  are invertible.   The general case can be found in \cite[Theorem 4.10]{GuoJiZhu2015}.

  \begin{lem}\label{lem3.20}
    Let $V, W\in \mathcal{B}(\hilb)$ be injective bilateral weighted shifts with positive weights. Then $W \simeq_{a} V$ if
    and only if $\sum_{k} V =\sum_{k} W$ for all $k\geq 1$.
  \end{lem}


  \begin{prop}\label{prop3.21}
    Let $n \ge 2$, and $W\in \cB(\hilb)$ be an injective bilateral weighted shift with positive
    weights $\{\cdots, w_{i},w_{i+1},w_{i+2},\cdots\}$. Then the following are
    equivalent:
    \begin{enumerate}
      \item[(i)]
	$|\Omega_W|\leq 1$, where $\Omega_W:=\{ i \in \bbZ: w_i\neq 1\}$.
      \item[(ii)]
	$W$ is $\jnu$-stable.
      \item[(iii)]
	$W$ is $\jna$-stable.
    \end{enumerate}
  \end{prop}

  \begin{pf}
    \begin{enumerate}
      \item [(i)] implies (ii). This is from Lemma~\ref{lem3.19} and \cite[Theorem 1(a)]{Shi74}.
	\bigskip
      \item [(ii)] implies (iii). This is obvious.
	\bigskip

      \item [(iii)] implies (i).
	Suppose to the contrary that $|\Omega_W|\geq 2$.
	Pick
	$s,t\in \Omega_W$, such that $L := t-s =d(W)\geq 1$.
	Obviously, for $s+1\leq i\leq t-1$, $w_i=1$.
	Since $w_s\neq 1\neq w_t$,
	we could find a sufficient small positive number $\delta<1$, such that
	\[1\notin (w_s-\delta,w_s+\delta)\cup (w_t-\delta,w_t+\delta).\]
	Denote the weight sequence of $J_n(W)$ as $\{v_i\}_{i\in \mathbb{Z}}$.
	Then $d(J_n(W))=n \, d(W)=nL$. Since $W\simeq_a J_n(W)$, by Lemma~\ref{lem3.20},
	there exists $i^*\in \mathbb{N}$, such that
	\[\|(v_{i^*},\cdots, v_{i^*+L})-(w_s,\overbrace{1,\cdots,1}^{L-1},w_t)\|_2<\delta.\]
	In particular, $|v_{i^*}-w_s|<\delta$ and $|v_{i^*+N}-w_t|<\delta$.
	Hence, $v_{i^*}\neq 1$, and $v_{i^*+L}\neq 1$. Thus,
	from the definition of $d(J_n(W))$.
	\[nL=d(J_n(W))\leq |i^*+L-i^*|=L.\]
	This contradiction implies that $|\Omega_W|\leq 1$.
    \end{enumerate}
  \end{pf}


  \begin{prop} \label{prop3.22}
    Let $\hilb$ be a complex, separable, infinite-dimensional Hilbert space with orthonormal basis $\{e_i\}_{i=0}^\infty$.
    Let $n \ge 2$, and $W\in \cB(\hilb)$ be a unilateral weighted shift satisfying $W e_i=w_ie_{i+1}$, $i\geq 0$. Suppose that
    $w_i> 0$  for all $i\geq 0$. Then $W$ is $\jnu$-stable if and only if $w_i=1$ for all $i\geq 0$.
  \end{prop}

  \begin{pf}
    For $x=\sum_{i=0}^\infty a_ie_i\in \hilb$,  define $U x=x_n\oplus \cdots \oplus x_1$,
    where $x_l=\sum_{i=0}^\infty a_{ni+(l-1)}e_i$, $1\leq l\leq n$.
    Clearly $U$ is a unitary which maps $\mathcal{H}$ onto $\mathcal{H}^{(n)}$.
    Then by using the similar argument in the proof of Lemma~\ref{lem3.19},
    $U^*J_n(W) U$ is a unilateral weighted shift with weight sequence $\{\mathbbm{1}_{n-1},w_0, \mathbbm{1}_{n-1},w_1,\cdots\}$.

    By \cite[Theorem 1(b)]{Shi74}, $U^*J_n(W) U$ is unitarily equivalent to $W$ if and only if the two sequences agree, i.e.
    \[
    (w_0,w_1,w_2,w_3,\cdots ) = (\mathbbm{1}_{n-1},w_0,\mathbbm{1}_{n-1},w_1,\cdots). \]
    So $W$ is $\jnu$-stable if and only if $w_i=1$ for all $i\geq 0$.
  \end{pf}


 Despite these results regarding weighted shifts, the set of $\jnu$-stable operators is surprisingly large.   As we shall now see, given any operator $A\in \cB(\hilb)$, we may associate to $A$ an \emph{operator-shift} $W_A$ which is $\jnu$-stable for all $n\geq 2$.

  \begin{prop}
    Let $\hilb$ be a complex Hilbert space, and $A\in \bofh$. Consider the tensor product space $\mathcal{K}=l^2(\mathbb{Z})\otimes \hilb$
    consisting of all sequences $x=(\cdots,x_{-2},x_{-1},[x_0],x_1,x_2\cdots)$   which are square-summable in the norm, where $[\cdot]$ denotes the ``zero"-th coordinate.
    Clearly $\mathcal{K}$ is a Hilbert space with the inner product
    \[ \langle x,y \rangle =\sum_{k\in \mathbb{Z}} \langle x_k,y_k \rangle.\]
    Let $W_A \in \mathcal{B}(\mathcal{K})$ be the  operator-shift which acts on $\cK$ by the equation $W_A x=\hat{x}$,
    where
    \[x=(\cdots,x_{-2},x_{-1},[x_0],x_1,x_2,\cdots),~~\hat{x}=(\cdots,x_{-2},[x_{-1}],Ax_0,x_1,x_2,\cdots).\]
    Then $W_A$ is $\jnu$-stable for all $n\geq 2$.
  \end{prop}

  \begin{pf}
    The proof is a variant of that of Lemma~\ref{lem3.19}.   We  prove the case where $n=3$.  The general case is a routine adaptation of this one.
    Given $x=(\cdots,x_{-2},x_{-1},[x_0],x_1,x_2,\cdots)\in \mathcal{K}$, define $U x=s\oplus y\oplus z$,
    where
    \[s=(\cdots,x_{-3},[x_0],x_{3},\cdots), y=(\cdots,x_{-4},[x_{-1}],x_{2},\cdots),z=(\cdots,x_{-5},[x_{-2}],x_{1},\cdots).\]
    Clearly $U$ is a unitary which maps $\mathcal{K}$ onto $\mathcal{K}^{(3)}$.
    Thus, $U^*=U^{-1}$. Note that
    \begin{align*}
      U^* & J_3(W_A) Ux \\
      &=  U^*J_3(W_A)[(\cdots,x_{-3},[x_0],x_{3},\cdots)\oplus(\cdots,x_{-4},[x_{-1}],x_{2},\cdots)\oplus(\cdots,x_{-5},[x_{-2}],x_{1},\cdots)]  \\
      &= U^*[(\cdots,x_{-4},[x_{-1}],x_{2},\cdots)\oplus(\cdots,x_{-5},[x_{-2}],x_{1},\cdots)\oplus (W_A (\cdots,x_{-3},[x_0],x_{3},\cdots))] \\
      &= U^*[(\cdots,x_{-4},[x_{-1}],x_{2},\cdots)\oplus(\cdots,x_{-5},[x_{-2}],x_{1},\cdots)\oplus (\cdots,x_{-6},[x_{-3}],Ax_0,\cdots)]\\
      &=(\cdots,x_{-2},[x_{-1}],Ax_0,x_1,x_2,\cdots)\\
      &=W_A x.
    \end{align*}
    Hence,  $U^*J_3(W_A) U=W_A$.
  \end{pf}



  \bigskip

  \section{Around a problem of Kaplansky's}


  \subsection{} \label{sec4.01}
  In his study of infinite abelian groups~\cite{Kaplansky1954}, Kaplansky raised ``\emph{three test problems}" to test the usefulness of a structure theorem for abelian groups.
  As Kaplansky himself noted:  ``\emph{all three problems can be formulated for general mathematical systems}".

  In the context of Hilbert space operators, Kaplansky's second problem becomes:   if $A, B \in \bofh$ and $A \oplus A$ is equivalent to $B \oplus B$, is $A$ equivalent to $B$?   In essence, one is seeking for an appropriate multiplicity theory for the notion of equivalence involved.  Of course, in order to solve Kaplansky's second problem, one must first define what is meant by ``equivalent", and the three notions we have already defined -- unitary equivalence, approximate unitary equivalence, and similarity are a good place to start.     We begin with a result of Kadison and Singer~\cite[Theorem~1]{KadisonSinger1957}.


  \begin{thm} \label{thm4.02} \emph{[\textbf{Kadison-Singer}]}
    Let $A, B \in \bofh$ and suppose that $A \oplus A \simeq B \oplus B$.   Then $A \simeq B$.
  \end{thm}

As pointed out in~\cite[Proposition~A.1]{Sherman2010},  if $A^{(n)} \simeq B^{(n)}$ for some $n \ge 2$, then $A \simeq B$.
As an immediate consequence of this result, we obtain:


  \begin{cor} \label{cor4.03}
    Let $A, B \in \bofh$ and suppose that there exists $n \ge 2$ such that $J_n(A) \simeq J_n(B)$.   Then $A \simeq B$.
  \end{cor}

  \begin{pf}
    This is an immediate consequence of the above result (see also~\cite{Azoff1995,Sherman2010}), and the fact that $J_n(A) \simeq J_n(B)$ implies that
    \[
    A^{(n)} = J_n(A)^n \simeq J_n(B)^n = B^{(n)}. \]
  \end{pf}

  More generally, this shows that $J_n(A) \simeq J_n(B)$ for some $n \ge 2$ if and only if $J_n(A) \simeq J_n(B)$ for all $n \ge 2$.

  We now argue that $A \oplus A$ \emph{approximately} unitarily equivalent to $B \oplus B$ implies that $A$ is \emph{approximately} unitarily equivalent to $B$.  In fact, slightly more is true.


  \begin{prop} \label{prop4.04}
    Let $A, B \in \bofh$,  and suppose that $A^{(n)} \simeq_a B^{(n)}$ for some integer $n \ge 2$.   Then $A \simeq_a B$.
  \end{prop}

  \begin{pf}
    First note that $C^*(A) \simeq^* C^*(A^{(n)})$ via the ${}^*$-isomorphism $\theta_1 (X) = X^{(n)}$ for all $X \in C^*(A)$.  Similarly, the map $\theta_2(Y) = Y^{(n)}$ for all $Y \in C^*(B)$ defines an isometric ${}^*$-isomorphism of $C^*(B)$ onto $C^*(B^{(n)})$.

    Also, since $A^{(n)} \simeq_a B^{(n)}$, the map $\varrho(A^{(n)}) = B^{(n)}$ extends to an isometric ${}^*$-isomorphism $\varrho$ of $C^*(A^{(n)})$ onto $C^*(B^{(n)})$ which is approximately unitarily equivalent to the identity representation on $C^*(A^{(n)})$.  In particular, therefore,
    \[
    \mathrm{rank}\, (X^{(n)}) = \rank\, \varrho( X^{(n)}) \mbox{ for all } X \in C^*(A). \]

    Let $\varphi: C^*(A) \to C^*(B)$ be the map
    \[
    \varphi = \theta_2^{-1} \circ \varrho \circ \theta_1. \]
    Then $\varphi$ is an isometric ${}^*$-isomorphism and
    \[
    \varphi(A) = \theta_2^{-1} \circ \varrho \circ \theta_1(A) = \theta_2^{-1} \circ \varrho(A^{(n)}) = \theta_2^{-1} (B^{(n)}) = B. \]
    It also follows that
    \[
    (\varphi(X))^{(n)} = \theta_2 \circ \varphi  (X)  = \varrho \circ \theta_1 (X) = \varrho(X^{(n)}) \]
    for all $X \in C^*(A)$.

    Thus for all $X \in C^*(A)$, we have
    \begin{align*}
      n \ \mathrm{rank}\, (\varphi(X))
      &= \mathrm{rank}\, (\varphi(X)^{(n)})  \\
      &= \mathrm{rank}\, (\varrho(X^{(n)})) \\
      &= \mathrm{rank}\, (X^{(n)}) \\
      &= n \ \mathrm{rank}\, (X),
    \end{align*}	
    and so $\mathrm{rank}\, \varphi(X) = \mathrm{rank}\, X$.

    By Voiculescu's Theorem~\cite[Theorem II.5.8]{Davidson1996}, $\varphi$ is approximately unitarily equivalent to the identity map on $C^*(A)$.   In particular,
    \[
    B = \varphi(A) \simeq_a A. \]
  \end{pf}


  \begin{rem} \label{rem4.05}
    Let $\hilb$ be a separable Hilbert space.  We recall that an operator $T \in \bofh$ is said to be \textbf{quasidiagonal} if there exists an increasing sequence $(P_n)_n$ of finite-rank projections tending in the strong operator topology to the identity operator such that
    \[
    \lim_n \norm P_n T - T P_n \norm = 0. \]
    L.~Brown~\cite{Brown1984} has demonstrated the remarkable result that there exist non-quasidiagonal operators $T_n \in \bofh, \ n \ge 2$ such that $T_n^{(n)}$ is quasidiagonal while $T_n^{(k)}$ is not quasidiagonal for any $2 \le k < n$.

    It is an easy consequence of Proposition~\ref{prop4.04} that we cannot write any such $T_n^{(n)}$ in the form $Q^{(n)}$ where $Q$ is quasidiagonal, for then $T_n$ would be approximately unitarily equivalent to $Q$, and any operator approximately unitarily equivalent to a quasidiagonal operator is itself quasidiagonal.
  \end{rem}


  \begin{cor} \label{cor4.06}
    Suppose that $A, B \in \bofh$.   The following statements are equivalent.
    \begin{enumerate}
      \item[(a)]
	$A \simeq_a B$.
      \item[(b)]
	For all $n \ge 2$, $J_n(A) \simeq_a J_n(B)$.
      \item[(c)]
	For some $n \ge 2$, $J_n(A) \simeq_a J_n(B)$.
    \end{enumerate}	
  \end{cor}

  \begin{pf}
    The only implication that is not completely routine is that (c) implies (a), and so we suppose that for some $n \ge 2$, $J_n(A) \simeq_a J_n(B)$.

    Since $J_n(A) \simeq_a J_n(B)$, it follows that $A^{(n)} \simeq (J_n(A))^n \simeq_a (J_n(B))^n \simeq B^{(n)}$.   The result now follows from Proposition~\ref{prop4.04}.
  \end{pf}

  \subsection{} \label{sec4.07}
  Kaplansky's second test problem in the context of similarity of Hilbert space operators reads as follows:   suppose that $A, B \in \bofh$ and $A \oplus A$ is similar to $B \oplus B$.   Is $A$ similar to $B$?

  At this level of generality, the answer is not yet known (unless $\hilb$ is finite-dimensional, in which case the Jordan form is easily seen to provide an affirmative solution to the problem).   For certain classes of operators, the problem does indeed admit a positive solution.   For example,  if both $A$ and $B$ are normal, then similarity of $A \oplus A$ to $B \oplus B$ implies that they are unitarily equivalent as well~\cite[Corollary IX.6.11]{Conway1985}, and the result follows from the previously cited result of Kadison and Singer~\cite{KadisonSinger1957}, or from general multiplicity theory for normal operators as outlined in ~\cite[Chapter II]{Davidson1996} or ~\cite[Chapter IX.10]{Conway1985}.

  We now consider a variant of Kaplansky's second problem which imposes a somewhat stronger hypothesis on the operators $A$ and $B$.    Recalling that $J_2(A)$ is a square root of $A \oplus A$ and $J_2(B)$ is a square root of $B \oplus B$, one may ask whether or not $J_2(A)$ similar to $J_2(B)$ implies that $A$ is similar to $B$.   The main result of this section of the paper asserts that the answer is yes if one of $A$ and $B$ is compact (in which case both are).

  We first seek a structure theorem for invertible operators which implement the similarity of $J_2(A)$ to $J_2(B)$.

  \bigskip


  \begin{lem} \label{lem4.08}
    Let $X, Y \in \bofh$ and suppose that $Z := \begin{bmatrix} X & Y \\ 0 & X \end{bmatrix}$ is a Fredholm operator of index zero in $\cB(\hilb \oplus \hilb)$.    Then $X$ is Fredholm of index zero in $\bofh$.
  \end{lem}

  \begin{pf}
    Let $z := \pi(Z) = \begin{bmatrix} x & y \\ 0 & x \end{bmatrix}$ denote the image of $Z$ in the Calkin algebra $\cB(\hilb \oplus \hilb)/\cK(\hilb \oplus \hilb)$, and observe that $z$ is invertible.   Write  $z^{-1}  = \begin{bmatrix} a & b \\ c & d \end{bmatrix}$.

    Then
    \[
      \begin{bmatrix} x a + y c & x b + y d \\ x c & x d \end{bmatrix} = z z^{-1} = \begin{bmatrix} 1 & 0 \\ 0 & 1 \end{bmatrix}
    = z^{-1} z = \begin{bmatrix} a x & a y + b x \\ c x & c y + d x \end{bmatrix}. \]
    In particular, $a x = 1 = x d$, and so $x$ is invertible in $\bofh/\kofh$, i.e. $X \in \bofh$ is Fredholm.   Define $ m:= \mathrm{ind}\, X \in \bbZ$.    Then $X \oplus X \in \cB(\hilb \oplus \hilb)$ is Fredholm of index $2 m$.    Since the Fredholm operators of any fixed index form an open set, there exists $\delta > 0$ such that if $W \in \cB(\hilb \oplus \hilb)$ and $\norm W - (X \oplus X) \norm < \delta$, then $W$ is Fredholm of index $2 m$.

    Note that for all $\eps > 0$, we have that $Z$ is similar to
    \[
    Z_\eps := \begin{bmatrix} \eps I & 0 \\ 0 & I \end{bmatrix} \cdot \begin{bmatrix} X & Y \\ 0 & X \end{bmatrix} \cdot \begin{bmatrix} \eps^{-1} I & 0 \\ 0 & I \end{bmatrix} = \begin{bmatrix} X & \eps Y \\ 0 & X \end{bmatrix}. \]

    If $0 < \eps < \frac{\delta}{2 \norm Y \norm + 1}$, then $\norm Z_\eps - (X \oplus X) \norm < \delta$, and so
    \[
    0 = \mathrm{ind}\, Z_\eps = \mathrm{ind}\, (X \oplus X) = 2m. \]
    In particular, $m = \mathrm{ind}\, X = 0$.
  \end{pf}


  \begin{lem} \label{lem4.09}
    Let $A, B \in \bofh$.   Suppose  that $S = \begin{bmatrix} S_1 & S_2 \\ S_3 & S_4 \end{bmatrix} \in \cB(\hilb \oplus \hilb)$ is invertible and that
    \[
    J_2(B) = S J_2(A) S^{-1}. \]
    Then
    \begin{enumerate}
      \item[(a)]
	$S_1 = S_4$,  $S_3 = S_2 A  = B S_2$, and $S_1 A = B S_1$.
      \item[(b)]
	Writing $S^{-1} = \begin{bmatrix} T_1 & T_2 \\ T_3 & T_4\end{bmatrix}$, we also have that $T_1 = T_4$ and $T_3 = T_2 B = A T_2$.
      \item[(c)]
	$S_i T_j$ commutes with $B$ and $T_j S_i$ commutes with $A$ for all $1 \le i, j \le 4$.
      \item[(d)]
	If $A$ is compact, then $S_1$ and $T_1$ are Fredholm operators of index zero.
    \end{enumerate}
  \end{lem}	

  \begin{pf}
    \begin{enumerate}
      \item[(a)]
	We have $S J_2(A) = J_2(B) S$, which means
	\[
	\begin{bmatrix} S_2 A & S_1 \\ S_4 A & S_3 \end{bmatrix} = \begin{bmatrix} S_3 & S_4 \\ B S_1 & B S_2 \end{bmatrix}. \]
	From this (a) immediately follows.
      \item[(b)]
	The proof is similar to that of (a) and is omitted.
      \item[(c)]
	The equation $S J_2(A)  = J_2(B) S$ implies that $S (A \oplus A) = S J_2(A)^2 = J_2(B)^2 S = (B \oplus B) S$, from which we find that $S_i A = B S_i$, $1 \le i \le 4$.   Similarly, $A T_j = T_j B$, $1 \le j \le 4$.

	But then for $1 \le i, j \le 4$, we get
	\[
	(T_j S_i) A = T_j (S_i A) = T_j (B S_i) = (T_j B) S_i = (A T_j) S_i = A (T_j S_i), \]
	and
	\[
	(S_i T_j) B = S_i (T_j B) = S_i (A T_j) = (S_i A) T_j = (B S_i) T_j = B (S_i T_j). \]
      \item[(d)]
	From (a), we have that $S = \begin{bmatrix} S_1 & S_2 \\ S_2 A & S_1\end{bmatrix}$.   Clearly $S$ is Fredholm of index zero, since it is invertible in $\cB(\hilb \oplus \hilb)$.   But $S_2 A$ is compact, and therefore
	\[
	Z := \begin{bmatrix} S_1 & S_2 \\ 0 & S_1 \end{bmatrix} \]
	is also a Fredholm operator of index zero.    By Lemma~\ref{lem4.08}, $S_1$ is Fredholm of index zero.    A similar argument applied to $S^{-1} = \begin{bmatrix} T_1 & T_2 \\ A T_2 & T_1 \end{bmatrix}$ yields that $T_1$ is Fredholm of index zero.
    \end{enumerate}
  \end{pf}


  The next result is the main result of this section.

  \bigskip

  \begin{thm} \label{thm4.10}
    Let $A \in \kofh$, $B \in \bofh$ and suppose that $J_2(A)$ is similar to $J_2(B)$.   Then $A$ is similar to $B$.
  \end{thm}

  \begin{pf}
    Under the given hypotheses, first observe that
    \[
    B \oplus B = J_2(B)^2 \sim J_2(A)^2 = A \oplus A \in \kofh, \]
    from which it immediately follows that $B$ is also compact.

    By Lemma~\ref{lem4.09}, if $S \in \bofh$ is invertible and $J_2(B) = S J_2(A) S^{-1}$, then
    \[
    S = \begin{bmatrix} S_1 & S_2 \\ S_2 A & S_1 \end{bmatrix} \]
    for an appropriate $S_1$ and $S_2$.  Furthermore, by Lemma~\ref{lem4.08}, $S_1$ has Fredholm index zero, and by Lemma~\ref{lem4.09},
    \[
    \tag{$\ast$} S_1 A =  B S_1. \]
    This has a number of interesting consequences.
    For example, the equation $(\ast)$ implies that $\ker\, S_1$ is invariant for $A$.  As a result of this, relative to $\hilb = (\ker\, S_1) \oplus (\ker\, S_1)^\perp$, we may write
    \[
    A = \begin{bmatrix} A_1 & A_2 \\ 0 & A_4 \end{bmatrix}. \]
One crucial observation for the argument below is that $0 \not \in \sigma(A_1)$.   Indeed, since $A_1$ acts on the finite-dimensional space $\ker\, S_1$, if $0 \in \sigma(A_1)$, then $0$ must be an eigenvalue of $A_1$.   If $0 \ne x_0 \in \ker\, A_1$, then $A x_0 = A_1 x_0 = 0$, implying that $x_0 \in \ker\, S_1$ and $x_0 \in \ker\, S_2 A$.   But then $\begin{bmatrix} x_0 \\ 0 \end{bmatrix} \in \ker\, S$, a contradiction.
\smallskip

    Similarly,  equation $(\ast)$ implies that $\ran\, S_1$ is invariant for $B$.    Hence relative to $\hilb = (\ran\, S_1)^\perp \oplus (\ran\, S_1)$, we may write
    \[
    B = \begin{bmatrix} B_1 & 0 \\ B_3 & B_4 \end{bmatrix}. \]

    As a map from $\hilb = (\ker\, S_1) \oplus (\ker\, S_1)^\perp$ to $\hilb = (\ran\, S_1)^\perp \oplus (\ran\, S_1)$, we may write
    \[
    S_1 = \begin{bmatrix} 0 & 0 \\ 0 & S_1^\circ \end{bmatrix}, \]
    where $S_1^\circ : (\ker\, S_1)^\perp \to \ran\, S_1$ is invertible.

    Then $S_1 A = B S_1$, as a map from  $\hilb = (\ker\, S_1) \oplus (\ker\, S_1)^\perp$ to $\hilb = (\ran\, S_1)^\perp \oplus (\ran\, S_1)$, becomes
    \[
    \begin{bmatrix} 0 & 0 \\ 0 & S_1^\circ \end{bmatrix}  \begin{bmatrix} A_1 & A_2 \\ 0 & A_4 \end{bmatrix} = \begin{bmatrix} B_1 & 0 \\ B_3 & B_4 \end{bmatrix} \begin{bmatrix} 0 & 0 \\ 0 & S_1^\circ \end{bmatrix} , \]
    or equivalently
    \[
    \begin{bmatrix} 0 & 0 \\ 0 & S_1^\circ A_4 \end{bmatrix} = \begin{bmatrix} 0 & 0 \\ 0 &  B_4 S_1^\circ \end{bmatrix}. \]

    In other words, $A_4 \in \cB((\ker\, S_1)^\perp)$ is similar to $B_4 \in \cB(\ran\, S_1)$.   Since $S_1$ is Fredholm of index zero, we know that $\dim (\ker\, S_1) = \dim\, (\ran\, S_1)^\perp$ and $\dim\, (\ker\, S_1)^\perp = \dim (\ran\, S_1)$, and so we can find a unitary operator $U$ such that $U (\ker\, S_1) = (\ran\, S_1)^\perp$ and $U (\ker\, S_1)^\perp = (\ran\, S_1)$.

    Then $U B U^* $ and $A$ are of the form
    \[
    U B U^*= \begin{bmatrix} B_1 & 0 \\ B_3 & B_4 \end{bmatrix}, \ \ \ A = \begin{bmatrix} A_1 & A_2 \\ 0 & A_4 \end{bmatrix} \in \cB((\ker\, S_1) \oplus (\ker\, S_1)^\perp),  \]
    and $A_4$ is similar to $B_4$ in $\cB( (\ker\, S_1)^\perp)$.  (Note:  technically, these are not the same $B_1, B_3$ and $B_4$ as above, but copies of these acting between the spaces $\ker\, S_1$ and $(\ker\, S_1)^\perp$ via the unitary $U$.)

    Let $\Omega = \sigma(A_1)$, and set $\Gamma := \Omega \cap \sigma(A_4)$.  Note that $0 \not \in \sigma(A_1)$ implies that $\Omega$ consists of non-zero eigenvalues of $A_1$, and hence of $A$.   Then $A_4 \sim A_{41} \oplus A_{44}$, where $\sigma(A_{41}) = \Gamma$ and $\sigma(A_{44}) \cap \Gamma = \varnothing$.   Since $B_4 \sim A_4$, we may also write $B_4 = B_{41} \oplus B_{44}$, where $\sigma(B_{41}) = \Gamma$ and $\sigma(B_{44}) \cap \Gamma = \varnothing$.  (We point out the fact that the fact that $A \oplus A \sim B \oplus B$ implies that the non-zero eigenvalues of $A$ coincide with the non-zero eigenvalues of $B$, and that their algebraic multiplicities coincide.   Since $A_4 \sim B_4$, the same statement holds for the non-zero eigenvalues of $A_4$ and of $B_4$ respectively.   From this it follows that $\sigma(A_1) = \sigma(B_1)$, including the algebraic multiplicities of these (necessarily) non-zero eigenvalues.)

    Moreover, by Lemma~\ref{lem2.16} (i.e. as deduced from an application of Rosenblum's operator),  the fact that $A_4$ and $B_4$ are similar implies that $A_{41}$ is similar to $B_{41}$ and $A_{44}$ is similar to $B_{44}$.

    Let $Z \in \bofh$ be the invertible operator such that
    \[
    Z^{-1} A Z = \begin{bmatrix} A_1 & A_{21} & A_{22} \\ 0 & A_{41}& 0 \\ 0 & 0 & A_{44} \end{bmatrix}. \]
    Then $\sigma(A_{44}) \cap \Omega = \varnothing$ implies that
    \[
    \widehat{A} := Z^{-1} A Z \sim \begin{bmatrix} A_1 & A_{21} &0 \\ 0 & A_{41}& 0 \\ 0 & 0 & A_{44} \end{bmatrix}. \]

    Similarly, $B$ is similar to
    \[
    \widehat{B} := \begin{bmatrix} B_1 & 0 & 0 \\ B_{21} & B_{41} & 0 \\ 0 & 0 & B_{44} \end{bmatrix}, \]
    where, as we noted above, $A_{44}$ is similar to $B_{44}$.

    Now
    \begin{enumerate}
      \item[(i)]
	$\widehat{A} \oplus \widehat{A}$ is similar to $\widehat{B} \oplus \widehat{B}$;
      \item[(ii)]
	\[
	\sigma(\begin{bmatrix} A_1 & A_{21} \\ 0 & A_{41} \end{bmatrix}) = \Omega= \sigma(\begin{bmatrix} B_1 & 0 \\ B_{21} & B_{41} 	\end{bmatrix}); \mbox{\ \ \ \ \ and }\]
      \item[(iii)]
	$\Omega \cap \sigma(A_{44}) = \varnothing$, $\Omega \cap \sigma(B_{44})=\varnothing$.
    \end{enumerate}
    It follows that $	(\begin{bmatrix} A_1 & A_{21} \\ 0 & A_{41} \end{bmatrix})^{(2)} \sim (\begin{bmatrix} B_1 & 0 \\ B_{21} & B_{41} 	\end{bmatrix})^{(2)}$.   But these are finite matrices, so by considering their Jordan forms, we immediately conclude that
    \[
    \begin{bmatrix} A_1 & A_{21} \\ 0 & A_{41} \end{bmatrix} \sim \begin{bmatrix} B_1 & 0 \\ B_{21} & B_{41} 	\end{bmatrix}, \]
    from which we deduce that $\widehat{A} \sim \widehat{B}$, and thus $A$ is similar to $B$.

  \end{pf}


  When $A$ and $B$ are not assumed to be compact, we may no longer conclude that $S_1$ is Fredholm of index zero, and the above proof immediately breaks down.   Indeed, the question of whether or not $A \oplus A$ similar to $B \oplus B$ implies that $A$ is similar to $B$ seems to be very difficult.   There exist some cases, however, where we can assert positive results.  First, we propose a possible ``line of attack" to solving the problem, in the hope that the reader may have more luck than we did.

  \begin{prop}\label{prop4.11}
    Let $A, B \in \bofh$.

    \begin{enumerate}
      \item[(a)] 	
	Every operator $Z \in \bofh$ in the commutant of $J_2(A)$ is of the form
	\[
	  \begin{bmatrix}
	    Z_1 &  Z_2\\
	    Z_2 A & Z_1 \\
	\end{bmatrix},\]
	where $Z_1$ and $Z_2$ commute with $A$.
      \item[(b)]
	Suppose that $S \in \cB(\hilb \oplus \hilb)$ is invertible and that $J_2(B) = S J_2(A) S^{-1}$.   By Lemma~\ref{lem4.09}(a), we may write $S = \begin{bmatrix} S_1 & S_2 \\ S_2 A & S_1 \end{bmatrix}$.

	Then $A$ is similar to $B$ if and only if there exists an invertible operator $Z = \begin{bmatrix}  Z_1 &  Z_2 \\ Z_3 & Z_4 \end{bmatrix}$ in the commutant of $J_2(A)$ such that
	\[
	W_1 := S_2  A Z_2  + S_1 Z_1 \]
	is invertible.
    \end{enumerate}
  \end{prop}

  \begin{pf}
    \begin{enumerate}
      \item[(a)]
	This follows immediately from Lemma~\ref{lem4.09} (a), setting $B = A$ in that Lemma.
      \item[(b)]
	First, suppose that there exists an invertible operator $Z = \begin{bmatrix}  Z_1 &  Z_2 \\ Z_3 & Z_4 \end{bmatrix}$ in the commutant of $J_2(A)$ such that
	$W_1 := S_2  A Z_2  + S_1 Z_1$
	is invertible.   Since $Z$ commutes with $J_2(A)$, we observe that
	\[
	J_2(B) = S J_2(A) S^{-1} = S Z J_2(A) Z^{-1} S^{-1}. \]
	By part (a) above, we have that $Z_3 = Z_2 A$ and $Z_4 = Z_1$, so that  $Z = \begin{bmatrix} Z_1 & Z_2 \\ Z_2 A & Z_1 \end{bmatrix}$.       Now, by Lemma~\ref{lem4.09}(a),
	\[
	S Z = \begin{bmatrix} S_1 & S_2 \\ S_2 A & S_1 \end{bmatrix} \begin{bmatrix} Z_1 & Z_2 \\ Z_2 A & Z_1 \end{bmatrix} = \begin{bmatrix} W_1 & W_2 \\ W_2 A & W_1 \end{bmatrix} \]
	where $W_2 = S_1 Z_2 + S_2 Z_1$.   By considering the $(1,1)$-entry of $(SZ) J_2(A) = J_2(B) (SZ)$, we see that $W_1 A = B W_1$.   Since $W_1$ is invertible by hypothesis, $B$ is similar to $A$.

	\bigskip
	Conversely, suppose that $B = W_1 A W_1^{-1}$, where $W_1 \in \bofh$ is invertible.    With $W := W_1 \oplus W_1$, it follows
	that $J_2(B) = W J_2(A) W^{-1}$.

	Next, suppose that $S \in \cB(\hilb \oplus \hilb)$ and
	\[
	S J_2(A) S^{-1} = J_2(B) = W J_2(A) W^{-1}. \]
	Then $S^{-1} W J_2(A) = J_2(A) S^{-1} W$, and so $Z := S^{-1} W$ is invertible and lies in the commutant of $J_2(A)$.  From (a) above, we may write $Z = \begin{bmatrix} Z_1 & Z_2 \\ Z_2 A & Z_1 \end{bmatrix}$ where $Z_1, Z_2 \in \bofh$ commute with $A$.  Using Lemma~\ref{lem4.09}(a), we see that we may also write $S = \begin{bmatrix} S_1 & S_2 \\ S_2 A & S_1 \end{bmatrix}$, and therefore
	\[
	\begin{bmatrix} W_1 & 0 \\ 0 & W_1 \end{bmatrix} = W = S Z =  \begin{bmatrix} S_1 Z_1 + S_2 Z_2 A &  S_1 Z_2 + S_2 Z_1 \\ S_2 A Z_1 + S_1 Z_2 A & S_2 A Z_2 + S_1 Z_1 \end{bmatrix}. \]
	In particular, $W_1$ is invertible and
	\[
	W_1 = S_1 Z_1 + S_2 Z_2 A =  S_2 A Z_2 + S_1 Z_1. \]
    \end{enumerate}

  \end{pf}


  If the spectrum of $A$ (and hence of $B$) is sufficiently nice, we can again obtain  positive results.

  \bigskip

  \begin{prop} \label{prop4.12}
    Suppose that $A, B \in \bofh$ are invertible, and  that there exists a net $(p_\alpha)_\alpha$ of polynomials such that  $X := \ttt{wot}-\lim_\alpha p_\alpha(A)$ and $Y := \ttt{wot}-\lim_\alpha p_\alpha(B)$ both exist,  $X^2 = A$, and $Y^2 = B$.   If
    \[
    J_2(A) \sim J_2(B), \]
    then $A \sim B$.
  \end{prop}

  \begin{pf}
    Note that $(\sigma(X))^2 = \sigma(X^2) = \sigma(A)$ implies that $X$ is invertible.  Now
    \[
    J_2(A) \sim \begin{bmatrix} I & 0 \\ 0 & X^{-1} \end{bmatrix} \begin{bmatrix} 0 & I \\ A & 0 \end{bmatrix} \ \begin{bmatrix} I & 0 \\ 0 & X \end{bmatrix} = \begin{bmatrix} 0 & X \\ X^{-1} A & 0 \end{bmatrix} = \begin{bmatrix} 0 & X \\ X & 0 \end{bmatrix} \sim X \oplus - X. \]
    Similarly, $J_2(B) \sim Y \oplus -Y$.
    Thus
    \[
    X \oplus -X \sim J_2(A) \sim J_2(B) \sim Y \oplus -Y. \]
    Choose $R \in \cB(\hilb \oplus \hilb)$ invertible such that
    \[
    R^{-1} (X \oplus -X) R =  Y \oplus -Y. \]
    Squaring both sides, we obtain
    \[
    R^{-1} (A \oplus A) R = B \oplus B. \]

    Let $(p_\alpha)_\alpha$ be a net of polynomials such that $\ttt{wot}-\lim_\alpha p_\alpha( A) = X$.    Then
    \begin{align*}
      Y \oplus Y
      &=\ttt{wot}-\lim_\alpha p_\alpha(B \oplus B) \\
      &= \ttt{wot}-\lim_\alpha p_\alpha( R^{-1} (A \oplus A) R)  \\
      &= R^{-1} (\ttt{wot}-\lim_\alpha (p_\alpha(A \oplus A))) R \\
      &= R^{-1} (X \oplus X) R.
    \end{align*}

    Thus
    \begin{align*}
      R^{-1} (2 X \oplus 0) R
      &= R^{-1} ((X \oplus -X) + (X \oplus X))R \\
      &= R^{-1} (X \oplus -X) R + R^{-1} (X \oplus X) R \\
      &= (Y \oplus -Y) + (Y \oplus Y) \\
      &= 2 Y \oplus 0.
    \end{align*}
    Since $X$ and $Y$ are invertible, this implies (by Lemma~\ref{lem2.16}) that $X$ is similar to $Y$, and therefore that $A = X^2$ is similar to $Y^2=B$.
  \end{pf}


  \begin{cor} \label{cor4.13}
    Let $A, B \in \bofh$ be invertible operators, and suppose that there exists a simply connected domain $\Omega \subseteq \bbC \setminus \{ 0\}$ such that $\sigma(A) \subseteq \Omega$.    If $J_2(A) \sim J_2(B)$, then $A$ is similar to $B$.
  \end{cor}

  \begin{pf}
    First note that -- as we have already seen -- $J_2(A) \sim J_2(B)$ implies that $\sigma(B) = \sigma(A)$.

    Since $0 \not \in \Omega$, it follows that the identity function $\iota(z) = z$, $z \in \Omega$ admits a holomorphic square root, say $q$.   Thus $q(A)$ is a norm-limit of polynomials in $A$, and $q(B)$ is a norm-limit of the same polynomials in $B$.
    The result now follows from Proposition~\ref{prop4.12}.
  \end{pf}


  In the case where one of the operators $A, B$ is an isometry, we can dispense with the stronger hypothesis that $J_2(A)$ is similar to $J_2(B)$ and we can obtain a positive solution to Kaplansky's original problem in the context of similarity (see Theorem~\ref{thm4.17} below).  Before doing so, however, we need to recall a few definitions and results from the literature.   The monograph~\cite{Paulsen2002} is an excellent reference.

  \begin{defn} \label{defn4.14}
    Recall that an operator $T \in \bofh$ is said to be \textbf{polynomially bounded} if there exists a constant $\kappa > 0$ such that for all polynomials $p \in \bbC[x]$, we have
    \[
    \norm p(T) \norm \le \kappa \norm p \norm_\infty, \]
    where $\norm p \norm_\infty := \sup \{ |p(z)| : |z| \le 1\}.$
    We say that $T$ is \textbf{completely polynomially bounded} if  there exists a constant $\kappa > 0$ such that
    for all $q = [q_{ij}] \in \bbM_n(\bbC[x])$, we have
    \[
    \norm [q_{ij}(T)] \norm \le \kappa \norm q \norm_\infty,  \]
    where $\norm q \norm_\infty := \sup \{ \norm q_{ij}(z) \norm: |z| \le 1\}$.
  \end{defn}


  \smallskip

  \subsection{} \label{sec4.15}
  We note that von Neumann's Inequality~\cite[Corollary~1.2]{Paulsen2002} is the assertion that every contraction $T \in \bofh$ is polynomially bounded, while Arveson~\cite{Arveson1969, Arveson1972} observed that the Sz.-Nagy Dilation theorem implies that every operator $X \in \bofh$ which is similar to a contraction is in fact completely polynomially bounded.

  The converse of this, namely that every completely polynomially bounded operator is similar to a contraction, is due to Paulsen~\cite{Paulsen1984.02}.

  \smallskip

  We also remind the reader that an operator $T \in \bofh$ is said to be \textbf{quasinormal} if $T$ commutes with $T^*T$.   It follows from the work of Brown~\cite{Brown1953} that $T$ is quasinormal if and only if $T$ is unitarily equivalent to an operator of the form $N \oplus (S \otimes P)$, where $S$ denotes the unilateral forward shift, $N$ is normal and $P$ is a positive definite operator.

  If $W \in \bofh$ is an isometry, then $W^* W = I$, and so clearly $W$ is quasinormal.    The Wold Decomposition for isometries acting on infinite-dimensional, separable Hilbert spaces asserts that there exists a unitary operator $U_W$ (perhaps acting on a space of dimensional zero) and a unique cardinal number $\alpha \in \bbN \cup \{ 0, \aleph_0\}$ such that $W \simeq U_W \oplus S^{(\alpha)}$.   Observing that $S^{(\alpha)} \simeq (S \otimes I_\alpha)$, where $\alpha$ denotes the identity operator acting on a Hilbert space of dimension $\alpha$, we see for isometries, the corresponding positive definite operator in Brown's structure theorem is nothing more than $I_\alpha$.

  \bigskip


  We shall need one more result, due to Popescu~\cite{Popescu1992}.

  \begin{thm} \label{thm4.16} \emph{\textbf{[Popescu]}}
    Let $A \in \bofh$ be such that
    \[
    P_A := \ttt{sot}-\lim_m (A^*)^m A^m \]
    exists.  Then $A$ is similar to an isometry if and only if $P_A$ is invertible.
  \end{thm}


  \begin{thm} \label{thm4.17}
    Let $T \in \bofh$ and $n \ge 1$ be an integer and let $W \in \bofh$ be an isometry.   If $T^{(n)}$ is similar to $W^{(n)}$, then $T$ is similar to $W$.  In particular, if $J_n(T)$ is similar to $J_n(W)$, then $T$ is similar to $W$.
  \end{thm}

  \begin{pf}
    Note that $\norm W \norm = 1 = \norm W^{(n)} \norm$.   It follows from Arveson's results that $T^{(n)}$ is completely polynomially bounded, which is easily seen to imply that $T$ is completely polynomially bounded, and thus by Paulsen's Theorem, that $T$ is similar to a contraction $X$.

    \smallskip

    If $A \in \bofh$ is a contraction, then the sequence $((A^*)^m A^m)_m$ is a decreasing sequence of positive operators, and so $P_A := \ttt{sot}-\lim_m (A^*)^m A^m$ exists.   In our case, we conclude that
    \[
    P_X := \ttt{sot}-\lim_m (X^*)^m X^m \mbox{\ \ \ \ and that \ \ \ \ } P_{X^{(n)}} := \ttt{sot}-\lim_m ((X^{(n)})^*)^m (X^{(n)})^m \]
    both exist, and clearly we must have
    \[
    P_{X^{(n)}} = (P_X)^{(n)}. \]
    But $X^{(n)}$ is similar to the isometry $W^{(n)}$, and so from Theorem~\ref{thm4.16}, $P_{X^{(n)}}$ is invertible.   This in turn implies that $P_X$ is invertible.   Once again applying Theorem~\ref{thm4.16}, we conclude that $X$ is similar to an isometry $V$, and consequently that $T$ is similar to $V$.
    By the Wold decomposition, $V$ is unitarily equivalent to $U_V \oplus S^{(\alpha)} \simeq U_V \oplus (S \otimes I_\alpha)$, where $U_V$ is unitary $\alpha \in \bbN \cup \{ 0, \aleph_0\}$, and $W$ is unitarily equivalent to $U_W \oplus S^{(\beta)} \simeq U_W \oplus (S \otimes I_\beta)$, where $U_W$ is unitary and $\beta \in \bbN \cup \{ 0, \beta\}$.
    In particular,
    \[
    U_V^{(n)} \oplus (S \otimes I_{n \cdot \alpha}) \simeq V^{(n)} \sim T^{(n)} \sim W^{(n)} \simeq U_W^{(n)} \oplus (S \otimes I_{n \cdot \beta}). \]

    Now $V^{(n)}$ and $W^{(n)}$ are quasinormal, and so by~\cite[Theorem~2.1]{ChenHerreroWu1992}, $U_V^{(n)} \simeq U_W^{(n)}$ and for all $\lambda \in \bbC$,
    \[
    \nul (I_{n \cdot \alpha} - \lambda I_{n \cdot \alpha}) = \nul(I_{n \cdot \beta} - \lambda I_{n \cdot \beta}). \]

    By the results of Kadison and Singer~\cite{KadisonSinger1957} (see Theorem~\ref{thm4.02} above), $U_V \simeq U_W$, and the condition on the nullity for $\lambda = 1$ implies that $\alpha = \beta$.   Thus $V$ is unitarily equivalent to $W$, and so $T$ is similar to $W$.

    The last statement easily follows from the above argument and the fact that if $J_n(T) \sim J_n(W)$, then $T^{(n)} = J_n(T)^n \sim J_n(W)^n = W^{(n)}$.
  \end{pf}


  \section{Remarks and observations.}

  \bigskip

  \subsection{} \label{sec5.01}
  Let $T \in \bofh$ and $2 \le n \in \bbN$.   As we have seen multiple times, $J_n(T)^n = T^{(n)}$, and thus if $T$ is $\jnu$-stable (resp. $\jns$-stable), then $T^n \simeq T^{(n)}$ (resp. $T^n \sim T^{(n)}$).    In their paper~\cite{ConwayPrajituraRodriguezMartinez2014}, the authors studied operators satisfying this last relation in the case where $n = 2$, referring to these phenomena as $\ttt{condition U}$ (i.e. $T^2 \simeq T^{(2)}$) and as $\ttt{condition S}$ (i.e. $T^2 \sim T^{(2)}$) respectively.   Unfortunately, they mistakenly state that ``a scalar multiple of an operator satisfying $\ttt{condition S}$ (or $\ttt{condition U}$) satisfies the same condition only if the scalar is either $0$ or $1$."

  A counterexample to this claim comes in the form of the forward unilateral shift operator $S$, as for any $\alpha \in \bbT$,
  \[
  (\alpha S)^2 = \alpha^2 S^2 \simeq S^2 \simeq S \oplus S \simeq (\alpha S \oplus \alpha S). \]


  \begin{eg} \label{eg5.02}
    There exist (positive, invertible) operators $H_1, H_2 \in \bofh$ such that $H_1 \sim H_2$ and $H_1 \simeq_a H_2$, but $H_1$ is not unitarily equivalent to $H_2$.

    Choose two non-atomic Borel measures $\mu_1, \mu_2$ on the interval $[1,2] \subseteq \bbR$, with the supports of $\mu_1$ and $\mu_2$  both equal to  $[1,2]$, such that $\mu_1$ and  $\mu_2$ are not mutually absolutely continuous. Let $H_k$ denote the (hermitian) multiplication operator $M_x$ acting on $L^2([1,2], \mu_k), \ k=1, 2$.

    Then $\sigma(H_1) = [1,2] = \sigma(H_2)$, and so $H_1$ and $H_2$ are approximately unitarily equivalent, by the  Weyl-von Neumann-Berg Theorem~\cite[Theorem II.4.4]{Davidson1996}, but they are not unitarily equivalent, since the underlying measures $\mu_1$ and $\mu_2$ are not mutually absolutely continuous~\cite[Corollary~II.3.6]{Davidson1996}.

    Let $Z_1=S \otimes H_1$, $Z_2=S \otimes H_2$, where $S$ is the usual unilateral forward shift operator. It is easy to check that $Z_1$ and $Z_2$ are approximately unitarily equivalent since $H_1$ and $H_2$ are. Furthermore, it follows from~\cite[Theorem 2.1]{ChenHerreroWu1992} that $Z_1$ and $Z_2$ are similar.

    Note, however,  that $Z_1$ and $Z_2$ are not unitarily equivalent.   For suppose otherwise.   Then
    \[
    H_1^2 \oplus 0 = (Z_1^* Z_1 - Z_1  Z_1^*) \simeq  (Z_2^* Z_2 - Z_2 Z_2^*)  = H_2^2 \oplus  0, \]
    and so by ~\cite[Theorem~1]{KadisonSinger1957}, we have that $H_1^2 \simeq H_2^2$.   Since $H_1, H_2$ are invertible and positive, we conclude that $H_1 \simeq H_2$, a contradiction.
  \end{eg}


  \begin{rem} \label{rem5.03}
    Let $A, B \in \bofh$.   The problem of whether or not $J_2(A) \sim J_2(B)$ implies that $A \sim B$ remains open.  In the hope that it may be of use to someone wishing to try their hand at this problem, we make a minor observation, which may be of some use:   if $J_2(A) \sim J_2(B)$ and $\kappa \in \bbC$, then $J_2(\kappa A) \sim J_2(\kappa B)$.

    \bigskip

    The result is obvious if $\kappa = 0$.   If $\kappa \not = 0$, then $J_2(A) \sim J_2(B)$ implies that
    \[
    \kappa^{1/2} J_2(A) \sim \kappa^{1/2} J_2(B).\]
    But then
    \[
    J_2(\kappa A) \sim \kappa^{1/2} J_2(A) \sim \kappa^{1/2} J_2(B) \sim J_2(\kappa B). \]	
  \end{rem}


  \begin{eg} (due to \textbf{J. Bell}) \label{eg5.04}
    It is natural to ask whether Kaplansky's second problem can be positively answered in the general setting of rings, which would obviously imply a positive answer in the setting of operators on a Hilbert space.  More specifically, is it true that if $R$ is a unital ring, and $x, y \in R$ are such that $x \oplus x$ is \emph{conjugate} (the ring-theoretic nomenclature for similarity) to $y \oplus y$ in $\bbM_2(R)$, then $x$ and $y$ must be conjugate?  In general, the answer is no.  In fact, in the general setting of unital rings, the analogue of Theorem~\ref{thm4.10} fails.    The following counterexample is due to our colleague Jason Bell, and we thank him for allowing us to reproduce it here.

    There exists a unital ring $R$ and elements $x, y \in R$ such that if $J_2(x) := \begin{bmatrix} 0 & 1 \\ x & 0 \end{bmatrix}$ and $J_2(y) := \begin{bmatrix} 0 & 1 \\ y & 0 \end{bmatrix}$ then $J_2(x)$ is conjugate to $J_2(y)$ in $\bbM_2(R)$, but $x$ is not conjugate to $y$ in $R$.

    Let $S = \bbQ[a,b,c]/(a^2-bc=1)$.  Notice that $S$ is an integral domain since $a^2-bc$ is irreducible.  We claim that the group of units of $S$ is $\bbQ^*$.  To see this, notice that $S$ admits a $\bbZ$-grading, by giving $a$ degree $0$, $b$ degree $1$ and $c$ degree $-1$.  Since $S$ is a $\bbZ$-graded domain, the only units of $S$ are contained in $S_0$, the homogeneous part of $S$ of degree $0$.  Indeed, since $S$ is a domain the set of units is homogeneous.   Consider $S$ as a module over $S_0$, and let $S_n$ denote the set of elements of $S$ of degree $n$, $n \in \bbZ$.   Then for $n \ge 1$, $S_n = S_0 b^n$, and $S_{-n} =S_0 c^n$.  Note further that  $b$ and $c$ are not units since $bc = a^2-1$ is not a unit in $S_0$.  Hence the group of units of $S$ lies in $S_0$.

    Observe that $S_0$ is the subring of $S$ generated by $a$ and $bc$.  Thus
    \[
    S_0 \cong \bbQ[a, bc]/(a^2-bc-1)\cong \bbQ[a]. \]
    Since the group of units of $\bbQ[a]$ is $\bbQ^*$, the claim follows.

    Next, let $T$ denote the free product (in the category of associative $\bbQ$-algebras) of $S$ and $\bbQ\{x,y\}$, the free associative algebra on $x$ and $y$.  Consider the ideal $I$ of $T$ generated by $ax-ya$, $bx-yb$, and $cx-yc$.
    If we let $R=T/I$, then in $\bbM_2(R)$, the matrix $\begin{bmatrix} a & b \\ c  & a \end{bmatrix}$ is invertible (since it lies in $\bbM_2(S)$, $S$ is commutative, and the determinant is $a^2-bc=1$).
    Furthermore,
    \[
    \begin{bmatrix} a & b \\ c  & a \end{bmatrix}(x \oplus x) = (y \oplus y)\begin{bmatrix} a & b \\ c  & a \end{bmatrix}  \in \bbM_2(R), \]
    by our definition of $I$.   Thus $x \oplus x$ and $y \oplus y$ are conjugate in $\bbM_2(R)$.

    Notice that $R$ admits an $\bbN$-grading (not compatible with the one on $S$) obtained by assigning degree $0$ to $a,b$ and $c$, and degree $1$ to $x$ and $y$, since the generators of $I$ and the relation $a^2-bc=1$ are all homogeneous.    Now $R/(x,y) \cong S$ and so $S$ is the degree $0$ component of $R$.

    We now argue by contradiction.  Suppose that $x$ and $y$ are conjugate in $R$.  Then there is some
    unit $u$ of $R$ such that $ux = yu$.   Write $u = u_0 +\cdots  + u_d$, where $u_i \in R$ is
    homogeneous of degree $i$.  Let $v=v_0+\cdots + v_e$ be the inverse of $u$.
    By considering the degree $0$ terms of both sides of $1=uv$, we see that $u_0v_0=1$ and since $u_0, v_0$ are in the degree $0$ component of $R \cong S$, we see that $u_0, v_0$ must be in $\bbQ^*$.

    Now we have that $(u_0+\dots+u_d)x = y(u_0+\dots+u_d)$ and so looking at the degree one part of both sides we see that $u_0 x = y u_0$.   But  $u_0$ is a nonzero scalar, from which it follows that $x=y$ in $R$.  Thus to show that $x$ and $y$ are not conjugate in $R$, it suffices to show that $x\ne y$ in $R$.  To do this, we shall define a homomorphism $\Phi: R \to \bbM_2(\bbQ)$ such that $\Phi(x) \ne \Phi(y)$.

    We consider the map $\Phi$ induced by setting
    \begin{itemize}
      \item{}
	$\Phi(1) = \begin{bmatrix} 1 & 0 \\ 0 & 1 \end{bmatrix}$;
      \item{}
	$\Phi(a) = \begin{bmatrix} 0 & 1 \\ 1 & 0 \end{bmatrix}$;  \ \ \ \ \ $\Phi(b) = \Phi(c) = \begin{bmatrix} 0 &  0 \\ 0 & 0 \end{bmatrix}$;
      \item{}
	$\Phi(x) = \begin{bmatrix} 1 & 0 \\ 0 & -1\end{bmatrix}$; and \ \ \ \ \ $\Phi(y) =  \begin{bmatrix}  -1 &  0 \\ 0 & 1 \end{bmatrix}$.
    \end{itemize}
    By definition of $R$ as the quotient of the free product $T$ of $S$ and $\bbQ\{ x, y \}$ by the ideal $I$,  to see that $\Phi$ extends to a ring homomorphism on $R$, it suffices to show that
    \begin{align*}
      0
      &= \Phi(a) \Phi(x) - \Phi(y) \Phi(a) \\
      &= \Phi(b) \Phi(x) - \Phi(y) \Phi(b) \\
      &= \Phi(c) \Phi(x) - \Phi(y ) \Phi(c), \mbox{ and } \\
      \Phi(1)	&=\Phi(a)^2 - \Phi(b) \Phi(c),
    \end{align*}	
    all of which are routine calculations.

    This implies that $\Phi$ indeed induces a homomorphism from $R \to \bbM_2(\bbQ)$.    But $\Phi(x) \ne \Phi(y)$, and thus $x\ne y$ in $R$, completing the argument.

    \bigskip

    As a corollary to this example, we note that $x \oplus x = J_2(x)^2$ is conjugate to $y \oplus y = J_2(y)^2$, and yet $x$ is not conjugate to $y$, and so this example also shows that this formulation of Kaplansky's second problem admits a negative answer in the setting of unital rings.
  \end{eg}

\bigskip

\centerline{\textbf{Statements and Declarations}}

\noindent \textbf{Data Availability}
 There is no dataset associated to this paper, as it is a study in pure mathematics.

\noindent \textbf{Competing Interests} There are no competing interests for the authors to disclose. The authors contribute equally to this piece of theoretical work.


\bibliographystyle{plain}

\end{document}